\documentclass[a4paper,12pt]{article}

\usepackage[utf8]{inputenc}
\usepackage[T1]{fontenc}
\usepackage[left=2.5cm,right=2.5cm,top=1.5cm,bottom=2cm]{geometry}
\usepackage{amsmath,amssymb,makeidx}
\usepackage[frenchb]{babel}

\usepackage[language=french,backend=biber]{biblatex}

\usepackage{fourier}
\usepackage{graphics}
\usepackage{graphicx}
\usepackage{amsthm}
\usepackage{thmtools}
\usepackage[hyperindex=true]{hyperref}
\usepackage{color}
\usepackage{enumitem}
\newlist{itemizeth}{itemize}{2}
\setlist[itemizeth]{label=\textbullet,noitemsep,topsep=0 mm}
\newlist{enumerateth}{enumerate}{2}
\setlist[enumerateth]{label*=\textbf{\alph*)},itemsep=0.2mm}
\setlist[itemize]{noitemsep,fullwidth, topsep=0 mm, leftmargin=0pt}
\setlist[enumerate]{label*=\textbf{\arabic*)},itemsep=0.2mm, leftmargin=0pt}

\declaretheoremstyle[
title=Démonstration,
numbered=no,
headfont=\normalfont\bfseries,
notefont=\bfseries, notebraces={}{},
postheadspace=17pt,
bodyfont=\normalfont,
headindent=0pt,
qed=\qedsymbol,
spacebelow=3,
]{demostyle}

\declaretheoremstyle[
headfont=\normalfont\bfseries,
notefont=\mdseries, notebraces={(}{)},
bodyfont=\normalfont,
thmbox=S
]{standard}

\declaretheoremstyle[
title=Théorème,
headfont=\normalfont\scshape,
notefont=\mdseries, notebraces={(}{)},
thmbox=S,
]{theo}

\declaretheoremstyle[
title=Théorème,
headfont=\scshape\bfseries,
notefont=\mdseries, notebraces={(}{)},
thmbox=L,
]{import}

\declaretheoremstyle[, 
numbered=no
]{rem}

\declaretheorem[title=Définition,style=standard]{defi}
\declaretheorem[title=Proposition,style=standard]{prop}
\declaretheorem[title=Lemme,style=standard]{lem}
\declaretheorem[title=Corollaire,style=standard]{coro}
\declaretheorem[style=theo]{Th}

\declaretheorem[title=Définition-proposition, style=standard]{defprop}
\declaretheorem[title=Exemple,style=rem]{Ex}
\declaretheorem[title=Remarque,style=rem]{rqu}
\declaretheorem[title=Remarques,style=rem]{rqus}
\declaretheorem[style=demostyle]{dem}

\makeatletter
\newcommand*{\house}[1]{%
   \mathord{%
     \mathpalette\@house{#1}%
   }%
}
\newcommand*{\@house}[2]{%
   \dimen@=\fontdimen8 %
       \ifx#1\scriptscriptstyle\scriptscriptfont
       \else\ifx#1\scriptstyle\scriptfont
       \else\textfont\fi\fi
       3 %
   \sbox0{%
     $#1%
       \vrule width\dimen@\relax
       \overline{%
         \kern2\dimen@
         \begingroup 
           #2%
         \endgroup
         \kern2\dimen@
       }%
       \vrule width\dimen@\relax
       \mathsurround=1.5\dimen@ 
     $%
   }%
   \ht0=\dimexpr\ht0-\dimen@\relax
   \dp0=\dimexpr\dp0+2\dimen@\relax
   \vbox{%
     \kern\dimen@ 
     \copy0 %
   }%
}

\newcommand{\R}{\mathbb{R}}
\newcommand{\N}{\mathbb{N}}
\newcommand{\D}{\mathbb{D}}
\newcommand{\Z}{\mathbb{Z}}
\newcommand{\Q}{\mathbb{Q}}
\newcommand{\C}{\mathbb{C}}
\newcommand{\K}{\mathbb{K}}
\newcommand{\Li}{\mathcal{L}}
\newcommand{\E}{\mathcal{E}}
\newcommand{\Oal}{\mathcal{O}}
\newcommand{\Qbar}{\overline{\mathbb{Q}}}
\newcommand{\id}{\mathrm{id}}
\newcommand{\Spec}{\mathrm{Spec} \,}
\newcommand{\Hom}{\mathrm{Hom}}
\newcommand{\Gal}{\mathrm{Gal} \,}
\newcommand{\GL}{\mathrm{GL}}
\newcommand{\Es}{\mathbb{E}}
\newcommand{\Pro}{\mathbb{P}}
\newcommand{\ord}{\mathrm{ord} \,}
\newcommand{\ssi}{\Leftrightarrow}
\newcommand{\gf}{\dfrac}
\newcommand{\ds}{\displaystyle}
\newcommand{\Vect}{\mathrm{Vect}}
\newcommand{\Tr}{\mathrm{Tr} \,}
\newcommand{\Img}{\mathrm{Im} \,}
\newcommand{\den}{\mathrm{den}}
\newcommand{\singze}{\left\lbrace 0 \right\rbrace}
\newcommand{\fonction}[5]{\begin{array}[t]{lrcl} 
#1: & #2 & \longrightarrow & #3 \\
    & #4 & \longmapsto & #5 \end{array}}

\title{Le théorème d'André-Chudnovsky-Katz}
\date{}
\author{Gabriel Lepetit}
\makeindex
\begin{document}

\makeatletter
  \begin{titlepage}
  \centering
      {\large \textsc{Université Grenoble Alpes \\ Institut Fourier}}\\
     
    \vfill
       {\LARGE \textbf{\@title}} \\
    \vspace{2em}
        {\large Gabriel \textsc{Lepetit}} \\
    \vspace{4em}
      {\large \emph{Mémoire de Master 2 réalisé sous la direction de Tanguy \textsc{Rivoal} }}\\
    \vspace{2em}
    \vfill
       {\large \textsc{janvier -- juin 2018}}
  \end{titlepage}

\makeatother

\newpage

\begin{abstract}
L'objet de ce mémoire est l'étude du théorème d'André-Chudnovsky-Katz sur la structure des solutions de l'équation différentielle d'ordre minimal non nulle à coefficients dans $\Qbar(z)$ satisfaite par une $G$-fonction. Nous commençons par présenter la théorie des opérateurs différentiels globalement nilpotents, dont le résultat principal est le théorème de Katz, qui affirme qu'ils sont fuchsiens à exposants rationnels. Puis nous donnons une preuve complète du théorème des Chudnovsky impliquant que l'opérateur différentiel minimal non nul à coefficients dans $\Qbar(z)$ d'une $G$-fonction satisfait une condition de croissance modérée sur certains dénominateurs, appelée \emph{condition de Galochkin}. Enfin, nous exposons la démonstration du théorème d'André-Bombieri établissant l'équivalence entre la condition de Galochkin et la condition de Bombieri, qui implique la nilpotence globale. Ceci nous permet de prouver le point principal du théorème d'André-Chudnovsky-Katz.
\end{abstract}

\tableofcontents
\newpage

\addcontentsline{toc}{section}{Introduction}

\section*{Introduction}

Une $G$-fonction est une série $f(z)=\sum\limits_{n=0}^{\infty} a_n z^n \in \Qbar\llbracket z\rrbracket $ telle que
\begin{enumerateth}
\item Il existe une équation différentielle linéaire à coefficients dans $\Qbar(z)$ dont $f(z)$ est solution ;
\item il existe une constante $C_1 >0$ telle que pour tout $n \in \N$, les valeurs absolues des conjugués au sens de Galois des $a_n$ sont bornées par $C_1^{n+1}$ ;
\item il existe $C_2 >0$ tel que, pour tout $n \in \N$, on puisse trouver un entier $d_n \geqslant 1$ vérifiant $$d_n a_0, \dots, d_n a_n \in \Oal_{\Qbar} \quad \mathrm{et} \quad  d_n \leqslant C_2^{n+1}.$$
\end{enumerateth}

Les $G$-fonctions constituent une classe importante de fonctions spéciales, dont l'étude s'est développée ces dernières années parallèlement à celle des $E$-fonctions, qui sont les séries $f(z)=\sum\limits_{n=0}^{\infty} \gf{a_n}{n!} z^n \in \Qbar\llbracket z\rrbracket $ solutions d'une équation différentielle sur $\Qbar(z)$ et telles que les $a_n$ vérifient les conditions \textbf{b)} et \textbf{c)} ci-dessus. En effet, un certain nombre de fonctions usuelles, telles que les polylogarithmes $\mathrm{Li}_s(z)=\sum\limits_{n=1}^{\infty} \gf{z^n}{n^s}$ -- dont les évaluations en 1 valent $\zeta(s)$ -- ou des fonctions hypergéométriques (voir section \ref{subsec:exempleGfons}), font partie de cette famille. Initialement, le terme de $G$-fonction, introduit par Siegel en 1929 dans \cite{Siegelarticle}, provient de l'exemple simple de la série géométrique $\sum\limits_{n=0}^{\infty} z^n$, de même que l'appellation $E$-fonction vient de la série exponentielle $\exp(z)=\sum\limits_{n=0}^{\infty} \gf{z^n}{n!}$.

Du point de vue de l'approximation diophantienne, le but de l'étude des $E$- et $G$-fonctions est de savoir si leurs valeurs en des points algébriques sont irrationnels, voire transcendants.

Concernant les $E$-fonctions, cette question a obtenu une réponse satisfaisante avec des raffinements successifs du théorème de Siegel-Shidlovskii (voir \cite[p. 139]{Shidlovskii}), dûs à Nesterenko et Shidlovskii (cf \cite{NesterenkoShidlovskii}), puis finalement à Beukers dans \cite{Beukers2006}, qui a utilisé la théorie des $E$-opérateurs d'André (cf \cite{AndregevreyI} et \cite{AndregevreyII}).

\begin{Th}[Siegel--Shidlovskii, 1929--1956] \label{siegelshidlovskii}
Soit $\mathbf{f}={}^t \left(
f_1(z), \dots, f_n(z)\right) \in \Qbar\llbracket z\rrbracket ^n$ un vecteur de $E$-fonctions vérifiant $\mathbf{f}'=G\mathbf{f}$, avec $G \in \mathcal{M}_n(\Qbar(z))$, soit $T(z) \in \Qbar[z]$ tel que $T(z) G(z) \in \mathcal{M}_n(\Qbar[z])$. Soit $\alpha  \in \Qbar$ tel que $\alpha T(\alpha) \neq 0$. 

Alors le degré de transcendance sur $\Qbar$ de $(f_1(\alpha),..., f_n(\alpha))$ est égal au degré de transcendance sur $\Qbar (z)$ de $(f_1(z), \dots, f_n(z))$.
\end{Th}

Le théorème \ref{siegelshidlovskii} est précisé \footnote{voir toutefois la remarque à la fin de cette introduction} par le  :

\begin{Th}[Beukers, 2006] \label{beukersefonctions}
Sous les hypothèses du théorème précédent, on a ; pour tout polynôme homogène $P \in \Q[X_1, \dots, X_n]$ tel que $P(f_1(\alpha), \dots, f_n(\alpha))=0$, il existe un polynôme $Q \in \Q[Z, X_1, \dots, X_n]$ homogène en les variables $X_1, \dots, X_n$ tel que $Q(\alpha, X_1, \dots, X_n)=P(X_1, \dots, X_n)$ et \\ $Q(z, f_1(z), \dots, f_n(z))=0$.
\end{Th}

Tout récemment, et à la suite des travaux précédents, Adamczewski et Rivoal ont trouvé dans \cite{AdamczewskiRivoal} un algorithme permettant de déterminer les points algébriques auxquels une $E$-fonction prend des valeurs algébriques.

\newpage
\begin{Th}[Adamczewski, Rivoal, 2018]\label{adamczewskirivoal}
Il existe un algorithme effectuant les tâches suivantes. Prenant une $E$-fonction $f(z)$ en entrée, il dit si $f(z)$ est transcendante sur $\Qbar(z)$ ou non. Si elle est transcendante, il donne la liste finie des nombres algébriques $\alpha$ tels que $f(\alpha)$ est algébrique, ainsi que la liste correspondante des valeurs $f(\alpha)$.
\end{Th}

Malheureusement, aucun résultat diophantien similaire n'a été prouvé concernant les $G$-fonctions. Un des théorèmes majeurs dans ce domaine est celui de Galochkin, conséquence du résultat prouvé dans \cite[p. 387]{Galochkin74}.

\begin{Th}[Galochkin, 1974] \label{thmgalochkin}
Soit $\mathbf{f}={}^t \left(f_1(z), \dots, f_n(z) \right)\in \Qbar\llbracket z\rrbracket ^n$ un vecteur de $G$-fonctions vérifiant $\mathbf{f}'=G\mathbf{f}$, avec $G \in \mathcal{M}_n(\Qbar(z))$, tel que $(f_1(z), \dots, f_n(z))$ est libre sur $\Qbar(z)$, et satisfaisant la \emph{condition de Galochkin} :

« si $T(z) \in \Qbar[z]$ est tel que $T(z) G(z) \in \mathcal{M}_n(\Qbar[z])$, et si pour tout $s \in \N$, $q_s$ est le plus petit entier tel que tous les coefficients des coefficients des matrices $TG, T^2 \gf{G_2}{2}, \dots, T^s \gf{G_s}{s!}$ sont des entiers algébriques, alors il existe une constante $C>0$ telle que pour tout $s \in \N, q_s \leqslant C^{s+1}$. »

Alors il existe une constante $c>0$ dépendant uniquement de $\mathbf{f}$ et de $G$ telle que pour tout $b \in \Z$ vérifiant $ c < |b|$, la famille $\left(f_1\left(\gf{1}{b}\right), \dots, f_n\left(\gf{1}{b}\right) \right)$ est libre sur $\Q$.
\end{Th}
Si $(a,b) \in \Z^2$, on a le même type de résultat sur les $f_i\left(\gf{a}{b}\right)$ pourvu que $0 < c_1 |a|^{c_2} < |b|$, avec $c_1, c_2$ des constantes.

La condition de Galochkin dans ce théorème illustre le fait que la nature des équations différentielles satisfaites par les $G$-fonctions influence de manière essentielle les propriétés de leurs valeurs en des points algébriques.  C'est pourquoi le but de ce mémoire est de présenter la démonstration du théorème suivant (cf \cite[pp. 718--719]{AndregevreyI}), qui étudie les caractéristiques de l'équation différentielle minimale d'une $G$-fonction.

\begin{Th}[André-Chudnovsky-Katz] \label{katzchudandre}
Soit $f(z) \in \Qbar\llbracket z\rrbracket $ une $G$-fonction, et $L \in \Qbar(z)\left[\gf{\mathrm{d}}{\mathrm{d}z}\right]$ un opérateur différentiel non nul d'ordre minimal pour $f$ tel que $L(f(z))=0$. Alors l'opérateur différentiel $L$ est globalement nilpotent. En particulier, tout point de $\mathbb{P}^{1}(\C)$ est un point singulier régulier de $L$ et les exposants de $L$ en tout point sont dans $\Q$. De plus, la condition de Galochkin est satisfaite pour $(f,f', \dots, f^{(\mu-1)})$, où $\mu$ est l'ordre de $L$. 
\end{Th}

Comme $L$ est à coefficients dans $\Qbar(z)$, on verra que les points finis non algébriques sont des points réguliers, donc à exposants entiers. Par ailleurs, comme l'explique André dans  \cite[p. 719]{AndregevreyI}, si $u \in \Qbar$, il existe une base de solutions de l'équation différentielle $L(g(z))=0$ au voisinage de $u$ s'écrivant $(f_1(z-u), \dots, f_n(z-u)) (z-u)^{C_u}$, où $C_u \in \mathcal{M}_n(\Q)$ est triangulaire supérieure, et les $f_i(z)$ sont des $G$-fonctions. Si $u=\infty$, on a une base similaire en remplaçant $z-u$ par $1/z$. Si $u$ est un point régulier, alors on peut prendre $C_u=I_n$, de sorte qu'autour de ce point, on peut trouver une base de solutions composées de $G$-fonctions uniquement.

Le théorème des Chudnovsky proprement dit, dû à G. et D. Chudnovsky dans \cite{Chudnovsky}, affirme que la condition de Galochkin introduite dans le théorème ci-dessus est automatiquement satisfaite par une famille de $G$-fonctions libre sur $\Qbar(z)$, ce qui est une avancée remarquable.

De plus, le théorème d'André-Chudnovsky-Katz est crucial dans la démonstration du théorème \ref{beukersefonctions} ci-dessus sur les $E$-fonctions. En effet la théorie d'André consiste à déduire des propriétés des équations différentielles des $G$-fonctions celles des équations différentielles des $E$-fonctions, via la transformée de Fourier-Laplace des opérateurs de $\Qbar\left[z, \gf{\mathrm{d}}{\mathrm{d}z} \right]$ :

$$ \begin{array}{cccc}
\mathcal{F} : & \Qbar\left[z,\gf{\mathrm{d}}{\mathrm{d}z} \right] & \longrightarrow & \Qbar\left[z,\gf{\mathrm{d}}{\mathrm{d}z} \right] \\ 
 & z & \longmapsto & -\gf{\mathrm{d}}{\mathrm{d}z} \\ 
 & \gf{\mathrm{d}}{\mathrm{d}z} & \longmapsto & z,
\end{array} $$ en faisant le parallèle avec le fait que la transformée de Laplace d'une $E$-fonction est une $G$-fonction en la variable $1/z$.

Dans le même article, les Chudnovsky ont obtenu le résultat diophantien suivant, qui est essentiellement le meilleur connu à ce jour (cf \cite[p. 17]{Chudnovsky}).

\begin{Th} \label{chudnovksy2}
Soit $\mathbf{f}={}^t \left(f_1(z), \dots, f_n(z) \right)\in \Qbar\llbracket z\rrbracket ^n$ un vecteur de $G$-fonctions vérifiant $\mathbf{f}'=G\mathbf{f}$, avec $G \in \mathcal{M}_n(\Qbar(z))$ telle que $1, f_1(z), \dots, f_n(z)$ sont algébriquement indépendants sur $\Qbar(z)$. Alors pour tout $t \geqslant 1$, il existe une constante $c_0 >0$ dépendant de $\mathbf{f}$ et de $t$ telle que pour tout nombre algébrique $\xi$ de degré inférieur ou égal à $t$, et vérifiant $$0< |\xi| \leqslant \exp\left(-c_0 \left(\log H(\xi) \right)^{\frac{4n}{4n+1}} \right),$$ les nombres $1, f_1(\xi), \dots, f_n(\xi)$ ne sont pas liés par une relation algébrique de degré inférieur ou égal à $t$ sur $\Q(\xi)$.

Ici, $H(\xi)$ est la hauteur naïve de $\xi$, c'est-à-dire le maximum des valeurs absolues des coefficients du polynôme minimal normalisé de $\xi$ sur $\Q$.
\end{Th}
Les preuves de ces théorèmes diophantiens ne seront pas abordées dans ce mémoire dans lequel nous nous concentrerons sur la démonstration du théorème d'André-Chudnovsky-Katz. Néanmoins, les techniques des Chudnovsky (partie \ref{subsec:thchudnovsky}) servent également à démontrer le théorème \ref{chudnovksy2}. 

\medskip

Nous prouverons le théorème \ref{katzchudandre} selon le plan suivant : on se donne $f \in \Qbar\llbracket z\rrbracket $ une $G$-fonction et $L \in \Qbar(z)\left[\gf{\mathrm{d}}{\mathrm{d}z}\right]$ un opérateur d'ordre minimal pour $f$ tel que $L(f(z))=0$.

\begin{itemize}
\item  Dans la partie \ref{sec:opdiffnilpotents}, nous nous intéresserons aux opérateurs différentiels nilpotents et au théorème de Katz (théorème \ref{katz}) affirmant que tout opérateur nilpotent sur $\Qbar(z)$ est singulier régulier en tout point de $\mathbb{P}^{1}(\Qbar)$ et que ses exposants en tout point sont rationnels.

\item La partie \ref{sec:GfonsChud} sera consacrée à la preuve du théorème des Chudnovsky (théorème \ref{chudnovsky}) affirmant que la matrice compagnon associée $A_L$ associée à $L$ vérifie la condition de Galochkin (définition \ref{galochkin}).

\item Dans la partie \ref{sec:nilpotencegalochkin}, nous démontrerons, à l'aide d'estimations sur le rayon de convergence $p$-adique, le théorème d'André-Bombieri (théorème \ref{bombieriandre}) montrant que la condition de Galochkin implique que $L$ est un opérateur différentiel globalement nilpotent. Ainsi, le théorème de Katz nous assure que $L$ est singulier régulier en tout point et que ses exposants en tout point sont rationnels.
\end{itemize} 

Une annexe sera consacrée au théorème de Chebotarev qui est crucial dans la preuve du théorème de Katz.

\bigskip

\textit{Remarque}

Dans son article \cite{Siegelarticle}, Siegel a défini les $E$-fonctions \og{} au sens large \fg{}, la condition de croissance géométrique sur les coefficients et sur les dénominateurs étant remplacée par une condition de domination, pour tout $\varepsilon >0$, par $(n!)^{\varepsilon}$ à partir d'un certain rang dépendant de $\varepsilon$. Il est clair qu'une $E$-fonction \og{} au sens strict \fg{} (c'est-à-dire telle que définie au début de l'introduction) est une $E$-fonction au sens large. Il est conjecturé que ces deux définitions sont en fait équivalentes, mais cela n'a pas été prouvé à ce jour. Le théorème \ref{siegelshidlovskii} a été prouvé pour des $E$-fonctions au sens large, mais pas les théorèmes \ref{beukersefonctions} et \ref{adamczewskirivoal}. En ce sens, le théorème \ref{beukersefonctions} n'est pas précisément une généralisation du théorème de Siegel-Shidlovskii.

\bigskip

\textit{Notations}

\begin{itemize}
\item On note $\Z_{-}$ l'ensemble des entiers relatifs négatifs.

\item $\Qbar$ désigne le corps des nombres algébriques, $\Oal$ l'anneau des entiers algébriques. Si $\K$ est un corps de nombres, $\Oal_\K=\Oal \cap \K$.

\item Si $\alpha \in \Qbar$, on note $\house{\alpha}$ la \emph{maison} de $\alpha$, c'est-à-dire le maximum des valeurs absolues des conjugués de $\alpha$ au sens de Galois.

\item Si $\K$ est un corps de nombres, on note $\Spec(\Oal_\K)$ le spectre de $\K$, c'est-à-dire l'ensemble des idéaux premiers de l'anneau de Dedekind $\Oal_\K$. Si $\mathfrak{p} \in \Spec(\Oal_\K)$ est tel que $\mathfrak{p} \cap \Z=(p), p \in \Spec(\Z)$, on note $\mathfrak{p} \mid p$.

\item Si $\K$ est un sous-corps de $\C$, on note $\K\left\{ z \right\}$ l'ensemble des séries entières de rayon de convergence non nul autour de 0 à coefficients dans $\K$.

\item Si $\K$ est un corps et $M \in \mathcal{M}_n(\K)$, la comatrice de $M$ sera notée $\mathrm{com}(M)=\left((-1)^{i+j} \Delta_{i,j} \right)_{1 \leqslant i, j \leqslant n}$, où $\Delta_{i,j}$ est le mineur d'ordre $i,j$ de $M$.

\end{itemize}

\newpage

\section{Opérateurs différentiels nilpotents}\label{sec:opdiffnilpotents}

\subsection{Généralités sur les opérateurs fuchsiens} 

Dans cette section, nous reprenons quelques éléments de la théorie des point singuliers réguliers des équations différentielles dans le plan complexe. Pour davantage de détails, on pourra se référer à  \cite{Hille}, \cite{Sauloy} ou \cite{Yoshida}. 

\subsubsection{Opérateurs différentiels}

\begin{defi}
Un \emph{corps différentiel} est un couple $(H, D)$, où $H$ est un corps et $D : H \rightarrow H$ est une \emph{dérivation}, c'est-à-dire une application vérifiant 

\begin{itemizeth}
\item $\forall x, y \in H, D(x+y)=D(x)+D(y)$ (additivité) ;
\item $\forall x, y \in H, D(xy)=xD(y)+yD(y)$ (règle de Leibniz).
\end{itemizeth}

On dit que $(H_2, D_2)$ est une \emph{extension différentielle} de $(H_1, D_1)$ si $H_1 \subset H_2$ et $D_{2|H_1}=D_1$.
\end{defi}

On se donne $H$ un corps différentiel muni d'une dérivation $D$. Alors $H_0 :=\left\{ x \in H : D(x)=0 \right\}$ est le \emph{corps des constantes} de $H$. Typiquement, on prendra $H=\K(z)$ muni de la dérivation usuelle des fractions rationnelles $D=\gf{\mathrm{d}}{\mathrm{d}z}$, avec $\K$ un corps fini ou un corps de nombres.

On note $H[D]=\left\{ \sum\limits_{i=0}^{N} a_i D^i, a_i \in H, N \in \N \right\}$ l'anneau non commutatif des opérateurs différentiels sur $H$.

\begin{defi}
L'\emph{ordre} de $L \in H[D]$, noté $\ord(L)$ est le degré de $L$ en tant que polynôme en $D$. 
\end{defi}

On remarque immédiatement que si $L_1, L_2 \in H[D], \ord(L_1 L_2)=\ord(L_1)+\ord(L_2)$. 

L'anneau $H[D]$ est euclidien à droite et à gauche pour le stathme $\ord$. En particulier, tout idéal à droite ou à gauche est principal. 

\begin{prop}
Si $L \in H[D]$ et $E$ est une extension différentielle de $H$, alors $L$ peut être vu comme un endomorphisme $H_0$-linéaire de $E$ via $x \mapsto L(x)$. En particulier, l'ensemble des solutions dans $E$, $\mathrm{Sol}(L,E):=\left\{ x \in H : L(x)=0 \right\}$, est un $H_0$-espace vectoriel de dimension inférieure ou égale à $\ord(L)$.
\end{prop}

\begin{defi}
On dit que $L \in H[D]$ est \emph{trivial} sur $H$ si $\dim_{H_0} \mathrm{Sol}(L,H)=\ord(L)$.
\end{defi}

\begin{prop} \label{produittriviaux}
Si $L_1, \dots, L_m \in H[D]$ et $L_1 \dots L_m$ est trivial sur $H$, alors pour tout $i \in \{ 1, \dots, m \}$, $L_i$ est trivial sur $H$.
\end{prop}

\begin{dem}[de la proposition \ref{produittriviaux}]
Montrons d'abord que si $L=L_1 L_2$, alors $$\dim_{H_0} \mathrm{Sol}(L,H) \leqslant \dim_{H_0} \mathrm{Sol}(L_1, H)+ \dim_{H_0} \mathrm{Sol}(L_2,H).$$
On a $\mathrm{Sol}(L_2) \subset \mathrm{Sol}(L)$. Fixons un supplémentaire $S$ du $H_0$-espace vectoriel $\mathrm{Sol}(L_2)$ dans $\mathrm{Sol} (L)$. Si $u \in S$, $L_1(L_2(u))=0$, donc $L_2$ envoie $S$ sur $\mathrm{Sol}(L_1)$. De plus $L_{2|S}$ est injectif car si $u \in S$ et $u \neq 0$, alors $L_2(u) \neq 0$. Par suite, 

$$ \dim_{H_0} \mathrm{Sol}(L,H) \leqslant \dim_{H_0} \mathrm{Sol}(L_1, H)+ \dim_{H_0} \mathrm{Sol}(L_2,H). $$

\medskip

Supposons que $L=L_1 \dots L_m$ est trivial. Alors en généralisant le résultat précédent par récurrence, on a $\ord L=\dim \mathrm{Sol}(L) \leqslant \sum\limits_{i=1}^m \dim \mathrm{Sol}(L_i) \leqslant \sum\limits_{i=1}^m \ord(L_i)=\ord L$, d'où pour tout $i \in \{ 1, \dots, m \}$, $\dim \mathrm{Sol}(L_i) = \ord L_i$. En d'autres termes, $L_i$ est trivial.\end{dem}

\begin{rqu}
La réciproque est fausse : si $H=\mathbb{F}_p(z)$ et $L=z \gf{\mathrm{d}}{\mathrm{d}z}$, alors $L^2$ est d'ordre 2 et a un noyau de dimension 1. En effet, on peut montrer que $H_0=\mathbb{F}_p(z^p)$, donc $1, z, \dots, z^{p-1}$ est une base de $H$ sur $H_0$.

De plus, si $1 \leqslant k \leqslant p-1$, $L^2(z^k)=k^2 z^k \neq 0$, donc $\mathrm{Sol}(L^2, H)=\mathbb{F}_p(z^p)$ est de dimension 1.
\end{rqu}

\bigskip

On rappelle le résultat suivant (voir \cite[p. 80]{Dwork}).

\begin{Th}
Si $L \in H[D]$, il existe une extension différentielle $E$ de $H$ telle que $L$ est trivial sur $E$.
\end{Th} 

Quand il n'y a pas d’ambiguïté, on note $\mathrm{Sol}(L) := \mathrm{Sol}(L,E)$.

\subsubsection{Points singuliers réguliers}

 Soit $\K$ un corps, on note $\delta=z \gf{\mathrm{d}}{\mathrm{d}z}$, qui est une dérivation sur $\K((z))$. Considérons un opérateur \begin{equation} \label{eq:operateurL} L=\delta^n+A_1(z) \delta^{n-1} + \dots + A_n(z) \in \K((z))\left[\gf{\mathrm{d}}{\mathrm{d}z}\right].\end{equation} On écrit pour tout $j \in \{ 1, \dots, n \}, A_j(z)=\sum\limits_{i=m}^{\infty} a_{i,j} z^i$, avec $m \in \Z_{-}$. 

On peut alternativement noter $L$ sous la forme  \begin{equation} \label{eq:operateurLddz} L=B_0(z)\left(\gf{\mathrm{d}}{\mathrm{d}z}\right)^n+B_1(z)\left(\gf{\mathrm{d}}{\mathrm{d}z}\right)^{n-1}+\dots+B_n(z), \quad B_j(z) \in \K((z))\left[\gf{\mathrm{d}}{\mathrm{d}z}\right].\end{equation}

Le passage de la forme \eqref{eq:operateurL} à la forme \eqref{eq:operateurLddz} se fait grâce au lemme suivant :

\begin{lem} \label{lemmepassagedzdelta}
Pour tout $m \in \N$, on a $$z^m \left(\gf{\mathrm{d}}{\mathrm{d}z}\right)^m=\delta(\delta-1)\dots (\delta-m+1) \quad \text{et} \quad \delta^m=z^m\left(\gf{\mathrm{d}}{\mathrm{d}z}\right)^m+\sum\limits_{k=1}^{m-1} a_{m,k} z^k \left(\gf{\mathrm{d}}{\mathrm{d}z}\right)^k,$$ où les $a_{m,k}$ sont des entiers naturels.
\end{lem}

\begin{Th}[Fuchs] \label{thfuchs}
Les assertions suivantes sont équivalentes :
\begin{enumerateth}
\item Pour tout $j \in \{1, \dots, n \}$, $\gf{B_j}{B_0}$ admet un pôle d'ordre au plus $j$ en $0$.
\item Pour tout $j \in \{1, \dots, n\}$, $A_j(z)$ n'a pas de pôle en $0$.
\end{enumerateth}
Dans ce cas, on dit que $0$ est un \emph{point singulier régulier} de $L$.
\end{Th}

\begin{dem}
Selon le lemme \ref{lemmepassagedzdelta}, il existe une matrice triangulaire supérieure $T \in \mathcal{M}_n(\Q)$, dont les éléments diagonaux valent tous $1$, telle que $${}^t (\delta^{n-1},\dots, \delta, 1 ) = T {}^t (z^{n-1} D^{n-1}, \dot, zD, 1) \quad \mathrm{avec} \;\; D=\gf{\mathrm{d}}{\mathrm{d}z}. $$ 

Donc $L$ est égal, à multiplication à gauche par un élément de $\K(z)$ près, à $$L=z^n D^n+z^{n-1} \widetilde{B}_1(z) D^{n-1}+\dots+\widetilde{B}_n(z), \quad \text{avec}\;\; \begin{pmatrix} \widetilde{B}_1 \\ \vdots \\ \widetilde{B}_n \end{pmatrix} = T^{-1} \begin{pmatrix} A_1 \\ \vdots \\ A_n \end{pmatrix}.$$ Ainsi, les $A_j(z)$ sont dans $\K\llbracket z\rrbracket $ si et seulement si les $\widetilde{B}_j(z)$ n'ont pas de pôle en $0$, c'est à dire si 
$\gf{z^{n-j}}{z^n} \widetilde{B}_j(z)$ a un pôle d'ordre au plus $j$ en $0$.

En comparant avec l'écriture \eqref{eq:operateurLddz}, on voit que pour tout $j$, $\gf{B_j(z)}{B_0(z)}=\gf{z^{n-j}}{z^n} \widetilde{B}_j(z)$ et le théorème \ref{thfuchs} s'ensuit.
\end{dem}

\bigskip

On suppose dans toute la suite que $L \in \K(z)\left[\gf{\mathrm{d}}{\mathrm{d}z}\right]$. 

Pour $a \in \mathbb{P}^1(\K)$, on définit l'opérateur $L_a \in \K(u)\left[\mathrm{d}/\mathrm{d}u\right]$, obtenu par changement de variable $u=z-a$, tel que pour toute fonction $f$, $L(f(z))=0$ si et seulement si $L_a(f_a(u))=0$, où $f_a(u)=f(u-a)$ si $a \neq \infty$ et $f_{\infty}(u)=g(u^{-1})$. Notons que $L_{a}$ est bien défini à multiplication à gauche par un élément de $\K(u)$ près. On a alors, en posant $\delta_u=u\gf{\mathrm{d}}{\mathrm{d}u}$, \begin{equation} \label{equationLa}L_a=B_0(u+a)\left(\gf{\mathrm{d}}{\mathrm{d}u}\right)^n+B_1(u+a)\left(\gf{\mathrm{d}}{\mathrm{d}u}\right)^{n-1}+\dots+B_n(u+a)=\delta_u^n+A_1(u+a)\delta_u^{n-1}+\dots+ A_n(u+a)\end{equation} et
\begin{equation}\label{equationLinfini} L_{\infty}=\delta_u^n-A_1\left(\gf{1}{u}\right)\delta_u^{n-1}+\dots+(-1)^n A_n\left(\gf{1}{u}\right).\end{equation}

\begin{defi}
On dit que $a \in \Pro^1(\K)$ est un point singulier régulier de $L$ si $0$ est un point singulier régulier de $L_a$.

L'opérateur $L$ est dit \emph{fuchsien} si tous les points de $\mathbb{P}^{1}(\overline{\K})$ sont singuliers réguliers.
\end{defi}

On peut donc réécrire le théorème de Fuchs en tout point de $\Pro^1(\K)$. En un point fini $a$, l'expression \eqref{equationLa} montre qu'on peut simplement remplacer $0$ par $a$ dans le théorème \ref{thfuchs}. En l'infini, on obtient le résultat suivant :

\begin{prop}
Les assertions suivantes sont équivalentes :
\begin{enumerateth}
\item L'opérateur $L$ admet $\infty$ pour point singulier régulier.
\item Pour tout $j \in \{1, \dots, n\}$, $A_j\left(u^{-1}\right) \in \K\llbracket u\rrbracket $).
\item Pour tout $j \in \{1, \dots, n \}$, $\gf{B_j(u^{-1})}{B_0(u^{-1})}$ admet un zéro d'ordre au moins $j$ en $0$.
\item En supposant que pour tout $j \in \{1, \dots, n \}$, $B_j(z) \in \K[z]$, alors $\deg(B_j) \leqslant \deg(B_0) - j$. En particulier, $\deg(B_0)=\max\limits_{1 \leqslant j \leqslant n} \deg(B_j)$.
\end{enumerateth}
\end{prop}

\begin{dem}

L'équivalence entre \textbf{a)} et \textbf{b)} découle de l'équation \eqref{equationLinfini}. 

Pour montrer l'équivalence entre \textbf{b)} et \textbf{c)}, on reprend les notations de la preuve du théorème \ref{thfuchs}. On a $$\begin{pmatrix} \widetilde{B}_1 \\ \vdots \\ \widetilde{B}_n \end{pmatrix} = T^{-1} \begin{pmatrix} A_1 \\ \vdots \\ A_n \end{pmatrix},$$ donc les $A_j(u^{-1})$ sont dans $\K\llbracket u\rrbracket $ si et seulement si les $\widetilde{B}_j(u^{-1})$ sont dans $\K\llbracket u\rrbracket $, c'est à dire si, pour tout $j$, $\gf{B_j(u^{-1})}{B_0(u^{-1})}=u^{j} \widetilde{B}_j(u^{-1})$ a un zéro d'ordre au moins $j$ en $0$ (ou autrement dit un pôle d'ordre au plus $-j$ en $0$).

Pour finir, pour prouver que \textbf{c)} et \textbf{d)} sont équivalents, on remarque que si $B_j(z) \in \K[z]$, l'ordre du pôle en $0$ de $B_j(u^{-1})$ est $\deg(B_j)$,  de sorte que $\gf{B_j(u^{-1})}{B_0(u^{-1})}$ a un pôle d'ordre $\deg(B_j)-\deg(B_0)$ en $0$. 
\end{dem}

\subsubsection{Polynôme indiciel et exposants d'un opérateur}

Soit $L$ écrit sous la forme \eqref{eq:operateurL}. Intéressons-nous maintenant à la notion de polynôme indiciel et d'exposants, qui permettent de préciser la structure de l'espace des solutions de l'équation $L(y(z))=0$ autour de $0$. On vérifie que pour tout $q \in \N$ et $s \in \Z$, $\delta^q z^s=s^q z^s$, donc, en étendant la dérivation $D$ à $\K((z))$, pour tout $s \in \Z$,

\begin{align}
L(z^s)&=s^n z^s+\sum\limits_{i=m}^{\infty} a_{i,1} z^i s^{n-1} z^s+\dots + \sum\limits_{i=m}^{\infty} a_{i,n} z^i z^s \notag \\
&= s^n z^s + z^{m+s} \left( \sum\limits_{i=m}^{\infty} \left( \sum\limits_{j=1}^n a_{i,j} s^{n-j} \right) z^{i-m} \right) \notag \\
&= z^{m+s} \left(\sum_{k=0}^{\infty} \phi_k(s) z^k \right), \label{Lzspolindiciel}
\end{align} où les $\phi_k$ sont des polynômes de degrés inférieurs à l'ordre $n$ de $L$.

\begin{defi} \label{defexposants}
\begin{itemizeth}
\item Supposons que $0$ est un point singulier régulier de $L$, alors $m=0$ et $$\phi_{0,L}(x) :=x^n+\sum\limits_{j=1}^n a_{0,j} x^{n-j}$$ est appelé \emph{polynôme indiciel} de $L$ en $0$. On le note simplement $\phi_0$ quand il n'y a pas d'ambiguité.

\item Les racines de $\phi_{0,L}(x)$ dans la clôture algébrique de $\K$ sont appelés les \emph{exposants} de $L$ en $0$.

\item Si $a \in \Pro^1(\K)$ est un point singulier régulier de $L$, le polynôme indiciel $\phi_{a,L}(x)$ de $L$ en $a$ est défini comme le polynôme indiciel de $L_a$ en $0$.
\end{itemizeth}
\end{defi}

\begin{rqu}
Si $0$ est un point régulier singulier de $L$, alors le polynôme indiciel s'écrit $$\phi_0(x)=x^n+A_1(0) x^{n-1}+\dots+A_n(0).$$

Plus généralement, les équations \eqref{equationLa} et \eqref{equationLinfini} montrent que si $a \in \K$ est un point singulier régulier, alors pour tout $j \in \{1, \dots, n \}$, $A_j(u+a) \in \K\llbracket u\rrbracket $, et le polynôme indiciel en $a$ est $$\phi_a(x)=x^n+A_1(a) x^{n-1}+\dots+A_n(a)  \;,$$ et si $\infty$ est une singularité régulière, alors $$\phi_{\infty}(x)=x^n-A_1(\infty) x^{n-1}+\dots+(-1)^n A_n(\infty),$$ où, par abus de notation, on note $A_j(\infty)=\lim\limits_{z \rightarrow \infty} A_j(z)$.
\end{rqu}

A l'aide du lemme \ref{lemmepassagedzdelta}, on peut déterminer en pratique rapidement l'équation indicielle d'une équation différentielle écrite sous la forme \eqref{eq:operateurLddz}.

\begin{prop}
Soit $L=\left(\gf{\mathrm{d}}{\mathrm{d}z}\right)^n+B_1(z)\left(\gf{\mathrm{d}}{\mathrm{d}z}\right)^{n-1}+\dots+B_n(z), \quad B_j(z) \in \K(z)$, supposons que $a \in \K$ est un point singulier régulier de $L$. Alors $$\phi_{a,L}(x)=x(x-1)...(x-n+1)+a_1 x(x-1)\dots (x-n+2)+\dots+a_n, \quad a_i=\lim_{z \rightarrow a} z^i B_i(z)$$ et si $\infty$ est un point singulier régulier de $L$,  on a $$\phi_{\infty,L}(x)=x(x+1)...(x+n-1)-a_1 x(x+1)\dots (x+n-2)+\dots+(-1)^n a_n, \quad a_i=\lim_{z \rightarrow \infty} z^i B_i(z)$$
\end{prop}

\begin{Ex}
Si $a \in \K$ est un point ordinaire de $L$, alors l'équation indicielle en $a$ s'écrit $x(x-1)\dots(x-n+1)=0$ (si $a=\infty$, $x(x+1)\dots(x+n-1)=0$), donc ses exposants sont $0,1, \dots, n-1$ (si $a=\infty$, ses exposants sont $0,-1,\dots,-(n-1)$).
\end{Ex}

Pour finir, la proposition suivante montre que le changement de variable et la notion de polynôme indiciel sont compatibles avec la multiplication dans $\K(z)\left[\gf{\mathrm{d}}{\mathrm{d}z}\right]$.

\begin{prop} \label{exposantsproduitop}
Soient $M, L \in \K(z)\left[\gf{\mathrm{d}}{\mathrm{d}z}\right]$. Soit $a \in \mathbb{P}^1(\K)$, supposons que $a$ est une singularité régulière de $M$ et de $L$. Alors 

\begin{enumerateth}
\item Les opérateurs $M_a$ et $L_a$ obtenus par changement de variable $u=z-a$ vérifient $(ML)_a=M_a L_a$. En d'autres termes, l'application $$\fonction{\Psi_a}{\K(z)\left[\gf{\mathrm{d}}{\mathrm{d}z}\right]}{\K(u)\left[\gf{\mathrm{d}}{\mathrm{d}u}\right]}{L}{L_a}$$ est un morphisme d'anneaux. 

\item Si de plus, $a$ est un point singulier régulier de $M$ et de $L$, alors c'est une singularité régulière de $ML$.

\item Sous les mêmes hypothèses que dans le point précédent, on a 
$$\phi_{a,ML}(x)=\phi_{a,M}(x) \phi_{a,L}(x).$$
Ainsi, l'ensemble des exposants en $a$ de $ML$ est l'union de l'ensemble des exposants en $a$ de $M$ et de l'ensemble des exposants en $a$ de $L$.
\end{enumerateth}
\end{prop}

Ainsi, pour connaître les exposants d'un opérateur différentiel, il peut être utile de le factoriser en produit d'opérateurs différentiels plus simples.

\begin{rqu}
Le point \textbf{c)} montre en particulier que la multiplication d'un opérateur différentiel par un élément de $\K(z)$ ne change pas les exposants. 
\end{rqu}

\begin{dem}
\begin{enumerate}[label=\textbf{\alph*)}]
    \item Si $a \neq \infty$, c'est clair au vu de la formule \eqref{equationLa}.
    
    Pour $a=\infty$, écrivons $M=\sum\limits_{j=0}^m b_j(z) \delta^j$ et $L=\sum\limits_{i=0}^n a_i(z) \delta^i$. Alors
    \begin{align}
    ML &= \left(\sum_{j=0}^m  b_j(z) \delta^j\right)\left(\sum_{i=0}^n a_i(z) \delta^i\right) = \sum_{j=0}^m \sum_{i=0}^n b_j(z) \delta^j(a_i(z) \delta^i) \notag\\ &\overset{\mathrm{Leibniz}}{=}\sum_{j=0}^m \sum_{i=0}^n b_j(z) \sum_{k=0}^j \binom{j}{k} \delta^{j-k}(a_i(z)) \delta^{i+k} \notag\\
        &= \sum_{\ell=0}^{m+n} \sum\limits_{i+j=\ell} b_j(z) \sum_{k=0}^j \binom{j}{k} \delta^{j-k}(a_i(z)) \delta^{i+k} \label{ML0singreg}
    \end{align} et de même  
    \begin{align*}
        \Psi_{\infty}(M) \Psi_{\infty}(L) &= (-1)^{m+n} \sum_{j=0}^m \sum_{i=0}^n (-1)^j b_j\left(\gf{1}{u}\right) (-1)^i \delta_u^j\left(a_i\left(\gf{1}{u}\right) \delta_u^i\right)\\
        &=(-1)^{m+n} \sum_{\ell=0}^{m+n} (-1)^\ell \sum_{i+j=\ell}^n b_j\left(\gf{1}{u}\right) \sum_{k=0}^j \binom{j}{k} \delta_u^{j-k}\left(a_i\left(\gf{1}{u}\right)\right) \delta^{i+k}.
    \end{align*}
    On montre par récurrence que si $a(z) \in \K(z)$, alors  $$\forall s \in \N, \;\; \delta^{s}(a(z))_{|z=u^{-1}}=(-1)^s \delta_u^s\left(a\left(\gf{1}{u}\right)\right).$$
    Donc \begin{align}
        \Psi_{\infty}(ML) &=(-1)^{m+n} \sum_{\ell=0}^{m+n} \sum\limits_{i+j=\ell} b_j\left(\gf{1}{u}\right) \sum_{k=0}^j \binom{j}{k} \left(\delta^{j-k}(a_i(z))\right)_{|z=u^{-1}} (-1)^{i+k} \delta_u^{i+k} \notag \\
        &= (-1)^{m+n} \sum_{\ell=0}^{m+n} (-1)^\ell \sum\limits_{i+j=\ell} b_j\left(\gf{1}{u}\right) \sum_{k=0}^j \binom{j}{k} \delta_u^{j-k}\left(a_i\left(\gf{1}{u}\right)\right) \delta_u^{i+k}=\Psi_{\infty}(M) \Psi_{\infty}(L). \label{MLinftysingreg}
    \end{align}
    \item Les formules \eqref{ML0singreg} et \eqref{MLinftysingreg} montrent que si $0$ (resp. $\infty$) est un point singulier régulier de $M$ et $L$ alors $0$ (resp. $\infty$) est un point singulier régulier de $ML$. Par translation, il en va de même en tout point de $\K$.
    \item Grâce au point \textbf{a)}, il suffit de montrer le résultat en $0$. Supposons que $0$ est un point singulier régulier de $M$ et de $L$. Selon l'équation \eqref{Lzspolindiciel}, si $s \in \Z$,
    \begin{align*}
        (ML)(z^s)&=M\left(\sum_{k=0}^{\infty} \phi_{k,L}(s) z^{k+s}\right)=\sum\limits_{k=0}^{\infty}M(\phi_{k,L}(s) z^{k+s})=\sum\limits_{k=0}^{\infty}\phi_{k,L}(s) M\left(z^{k+s}\right) \quad \text{car} \;\;  \phi_{k,L}(s) \in \K \\
        &=\sum\limits_{k=0}^{\infty}\phi_{k,L}(s)\sum_{\ell=0}^{\infty} \phi_{\ell,M}(k+s) z^{k+\ell+s} \in \phi_{0,L}(s) \phi_{0,M}(s) z^s + z^{s+1} \K\llbracket z\rrbracket,
    \end{align*} si bien que, par définition, le polynôme indiciel de $ML$ en $0$ est $\phi_{0,L}(x) \phi_{0,M}(x)$.
\end{enumerate}
\end{dem}

\subsubsection{Théorème de Frobenius}

Le théorème suivant donne la structure des solutions d'une équation différentielle autour d'un point singulier régulier en fonction des exposants de cette équation différentielle. On pourra trouver une preuve dans \cite[p. 349]{Hille}.

\begin{Th}[Frobenius]
Soit $L=\left(\gf{\mathrm{d}}{\mathrm{d}z}\right)^n+B_1(z)\left(\gf{\mathrm{d}}{\mathrm{d}z}\right)^{n-1}+\dots+B_n(z) \in \K(z)\left[\gf{\mathrm{d}}{\mathrm{d}z}\right]$ un opérateur d'ordre $n$ et $a \in \Pro^1(\K)$ un point singulier régulier de $L$. Notons $\overline{\K}$ la clôture algébrique de $\K$. Alors 

\begin{enumerateth}
\item Il existe une base de solutions au voisinage du $a$ de l'équation $L(y(z))=0$ de la forme $(f_1(z-a), \dots, f_n(z-a))(z-a)^B$, où $(z-a)^B := \exp(B \log(z-a))$, étant donnée une détermination du logarithme complexe, $f_i(z) \in \K\llbracket z\rrbracket $ et $B \in \mathcal{M}_n(\overline{\K})$.
\item Supposons que les séries de Laurent définissant les $B_j(z)$ autour de $0$ ont un rayon de convergence au moins $R$. Alors le rayon de convergence des $f_i(z)$ est au moins~$R$.
\item Si $\rho$ est une valeur propre de $B$, alors $\rho$ est un exposant de $L$ en $a$, c'est le minimum des exposants de la forme $\rho+k, k \in \Z$.
\item Si les exposants $\rho_1, \dots, \rho_{n}$ de $L$ en $a$ sont  tels que $\forall i \neq j, \rho_i - \rho_j \not\in \Z$, alors il existe une base de solutions de la forme $((z-a)^{\rho_1} f_1(z-a), \dots, (z-a)^{\rho_{n}} f_{n}(z-a))$, où $f_i(u) \in \overline{\K}\llbracket u\rrbracket $.
\end{enumerateth} 
\end{Th}

Concrètement, en remarquant que $B$ est semblable à une matrice sous forme de Jordan, on voit que les solutions de l'équation $L(y(z))$ seront des \emph{séries Nilsson-Gevrey}, c'est à dire de la forme $$y(z)=\sum\limits_{\alpha,k} z^{\alpha} (\log^{k}(z)) f_{\alpha,k}(z)\;, \quad f_{\alpha,k}(z) \in \overline{\K}\llbracket z\rrbracket $$ et les $\alpha$ apparaissant dans la somme sont des exposants de $L$ en $a$.

\subsection{Opérateurs nilpotents et théorème de Katz-Honda} \label{subsec:nilpotentkatzhonda}

La référence principale pour cette section et la suivante est \cite[pp. 77--100]{Dwork}

On se donne dans cette section un corps $\K$ (par exemple $\K=\mathbb{F}_p$) de caractéristique~$p>~0$, $\mathcal{F}=\K(z)$ muni de la dérivation $D=\gf{\mathrm{d}}{\mathrm{d}z}$. On note $\mathcal{F}_0$ le corps des constantes de $\K(z)$. Par exemple, si $\K=\mathbb{F}_p$, on vérifie que $\mathcal{F}_0=\mathbb{F}_p(z^p)$. Le but de cette section est d'étudier les exposants d'un opérateur différentiel sur $\K(z)$, le théorème de Katz-Honda (théorème \ref{katzhonda} ci-dessous) affirmant qu'ils sont nécessairement à coefficients dans $\mathbb{F}_p$.

\begin{defi} \label{definilpotent}
L'opérateur différentiel $L \in \K(z)[D]$ est dit \emph{nilpotent} lorsque $L=L_1 \dots L_m$, où les $L_i \in \K(z)[D]$ sont triviaux sur $\K(z)$.
\end{defi}

\begin{Th} \label{caracterisationnilpotent}
Soit $L \in \K(z)[D]$ d'ordre $n$, alors 
\begin{enumerateth}
\item Si $L$ est nilpotent, alors $D^{pn} \in \K(z)[D] L$.
\item S'il existe $\mu \in \N$ tel que $D^{p {\mu}} \in \K(z)[D] L$, alors $L$ est nilpotent.
\end{enumerateth}
\end{Th}

\begin{dem} [du théorème \ref{caracterisationnilpotent}]

\textbf{a)} \textbf{Premier cas} : $L$ est trivial sur $\K(z)$. Effectuons la division euclidienne de $D^p$ par $L$ : $D^p=AL+B$, où $\ord B < \ord L$. Soit $(u_1, \dots, u_n)$ une base de solutions sur $\mathcal{F}_0$ de l'équation $L(u)=0$, alors pour tout $i$, $D^{p}(u_i)=B(u_i)=0$, donc $B=0$ car $\dim_{\mathcal{F}_0} \mathrm{Sol}(B) \leqslant \ord B$. Donc $D^p \in \K(z)[D] L$.

\medskip

\textbf{Second cas} : $L$ est non trivial sur $\K(z)$. Prouvons le résultat par récurrence sur $n=\ord L$. C'est évident pour $n=1$. Si $n >1$, $L$ est nilpotent donc on peut écrire $L=L_1 L_2$, avec $L_i$ nilpotent sur $\K(z)$ et $n_i=\ord L_i < n$. Donc il existe $B_1, B_2 \in \K(z)[D]$ tels que $D^{p n_1}=B_1 L_1$ et $D^{p n_2}=B_2 L_2$. 

\medskip

Remarquons que si $M \in \K(z)[D]$, alors $D^p M=M D^p$. En effet, si $j \in \N$, $j \geqslant p$, $D^p(z^j)=p ! \times \binom{j}{p} z^{j-p}=0$, donc, par linéarité, $D^p(u)=0$ pour tout $u \in \K[z]$, d'où $D^p(u)=0$ pour $u \in \K(z)$. Donc par la règle de Leibniz, $D^p$ est $\K(z)$-linéaire, ce qui autorise la commutation annoncée.

Par conséquent, ici, $$ B_2 B_1 L=B_2 (B_1 L_1 L_2)=B_2 D^{p n_1} L_2=D^{p n_1} B_2 L_2=D^{p n_1} D^{p n_2}=D^{p n},$$ ce qui conclut la preuve de \textbf{a)}.

\bigskip

Passons au point \textbf{b)}. On procède par récurrence sur $\ord L$. Prenons $L$ tel que $\ord L>0$.

On choisit $\mu \in \N$ minimal tel que $D^{p \mu} \in \K(z)[D] L$. On considère, dans $\K(z)[D]$, l'idéal à gauche  $\left( D^p, L\right)$ qui est principal, donc on peut fixer $L_1 \in \K(z)[D]$ tel que $\left(D^p, L\right) = \left( L_1 \right)$. 

Si $1 \in \left( D^p, L\right)$, alors il existe $\alpha, \beta \in \K(z)[D]$ tels que $1=\alpha D^p+ \beta L$, donc $$D^{p(\mu-1)}=\alpha D^{p \mu}+\beta D^{p(\mu-1)} L \in \K(z)[D] L,$$ ce qui contredit la minimalité de $\mu$. 

Ainsi, $\ord L_1 >0$ et $L=A_1 L_1$ pour un certain $A_1 \in \K(z)[D]$. Comme $D^p \in \K(z)[D] L_1$ et $D^p$ est trivial sur $\K(z)$ -- une base de $\mathrm{Sol}(D^p, \K(z))$ étant donnée par $1, z, \dots, z^{p-1}$ --, $L_1$ est trivial sur $\K(z)$. Il reste à montrer que $A_1$ est nilpotent.

Soit $E$ extension de $\K(z)$ sur laquelle $L$ est triviale. Alors selon la proposition \ref{produittriviaux}, $A_1$ est triviale sur $E$ et par ailleurs, en reprenant la preuve de cette proposition, on a $L_1(\mathrm{Sol}(L,E))=\mathrm{Sol}(A_1,E)$.

Si $u \in \mathrm{Sol}(L,E)$, alors $D^{p \mu}(u)=0$ car $D^{p \mu} \in \K(z)[D]L$ et comme $D^{p\mu}$ commute avec $L_1$ (cf point \textbf{a )}), on a $0=L_1(D^{p \mu}(u))=D^{p \mu}(L_1(u))$, donc $D^{p \mu}_{|L_1(\mathrm{Sol}(L,E))}=0=D^{p \mu}_{|\mathrm{Sol}(A_1,E))}$.

Écrivons $$D^{p \mu}=A_2 A_1+B, \quad \ord B < \ord A_1$$ la division euclidienne de $D^{p \mu}$ par $A_1$. En évaluant l'égalité en $u \in \mathrm{Sol}(A_1,E)$, on voit que $B=0$. Donc $D^{p \mu} \in \K(z)[D] A_1$. De plus, $\ord A_1 < \ord L$ car $\ord L_1 >0$. Par hypothèse de récurrence, $A_1$ est nilpotent, d'où comme $L_1$ est trivial sur $\K(z)$, $L$ est nilpotent.\end{dem}

Soit $\Omega$ une extension de $\K(z)$ telle que $D^{p^s}$ est trivial sur $\Omega$ pour un certain $s \in \N$. On note $\Omega_s=\mathrm{Sol}(D^{p^s}, \Omega)$ et $\Omega_0$ le corps des constantes de $\Omega$.

\begin{prop}
L'ensemble $\Omega_s$ est un $\Omega_0$-espace vectoriel qui est aussi un corps différentiel.
\end{prop}

\begin{dem}
Soient $y, z \in \Omega_s$. Par la formule de Leibniz, $$D^{p^s}(yz)=\sum\limits_{i=0}^{p^s} \binom{p^s}{i} D^i(y) D^{p^s-i}(z).$$ Or si $1 \leqslant i \leqslant p^s-1, \binom{p^s}{i}$ est divisible par $p$ donc nul dans $\Omega$ de caractéristique $p$. D'où $D^{p^s}(yz)=y D^{p^s}(z)+z D^{p^s}(y)=0$, si bien que $yz \in \Omega_s$.

Par le même raisonnement, si $y \in \Omega_s$ et $yz=1$, alors $0=yD^{p^s}(z)+z \times 0$,  de sorte que $z \in \Omega_s$ car $y \neq 0$.\end{dem}

\begin{prop} \label{caracterisationnilpotenceomegas}
Si $L \in \K(z)[D]$ est d'ordre $n$ et $p^{s-1} \geqslant n$, alors $L$ est nilpotent si et seulement $L$ est trivial sur $\Omega_s$. 
\end{prop}

\begin{dem}
\begin{itemize}
\item Si $L$ est nilpotent, alors selon le théorème \ref{caracterisationnilpotent}, $D^{pn} \in \K(z)[D] L$, donc, pour $p^s \geqslant pn$, $D^{p^s} \in \K(z)[D] L$. On note $D^{p^s}=AL, A \in \K(z)[D]$. Sur $\Omega_s$, $D^{p^s}$ est trivial car $\mathrm{Sol}(D^{p^s}, \Omega_s)$ contient $\mathrm{Sol}(D^{p^s}, \Omega)$ qui est de dimension $p^s$ sur $\Omega_0$. Par conséquent, $L$ est trivial sur $\Omega_s$. 

\item Réciproquement, si $L$ est trivial sur $\Omega_s$ pour $p^s \geqslant pn$, on écrit $D^{p^s}=AL+B$, avec $\ord B < \ord L$.\end{itemize}
\noindent  Donc comme $L$ est trivial sur $\Omega_s$, on obtient, en évaluant sur une base de solutions de $L(u)=~0$, $B=0$,  de sorte que $D^{p^s} \in \K(z)[D]L$.\end{dem}

\bigskip

\begin{Th}[Katz-Honda] \label{katzhonda}
Soit $L=\delta^n+A_1 \delta^{n-1} + \dots + A_n \in \K(z)\left[\gf{\mathrm{d}}{\mathrm{d}z}\right]$. Si $L$ est nilpotent, alors 
\begin{enumerateth}
\item $0$ est un point singulier régulier de $L$.
\item Tous les exposants de $L$ en $0$ sont dans $\mathbb{F}_p$.
\end{enumerateth}
\end{Th}

\begin{dem}
Procédons par récurrence sur $n=\ord(L)$. 

\begin{itemize}
\item Si $n=1$, alors $L$ est trivial sur $\K(z)$, donc on peut fixer $u \in \K(z) \setminus \{ 0 \}$ tel que $L(u)=0$, de sorte que $L=\delta-\gf{\delta(u)}{u}$. On écrit $u=z^{\ell} v$, avec $v \in \K\llbracket z\rrbracket ^{\times}$ et $\ell \in \N$. Donc

$$ \gf{\delta(u)}{u}=\gf{z(\ell z^{\ell-1} v+z^\ell D(v))}{u}=\gf{\ell z^{\ell} v+ z^\ell \delta(v)}{z^\ell v}=\ell+\gf{\delta(v)}{v}. $$

Par suite, $L=\delta-\ell-\gf{\delta(v)}{v}$ et $\gf{\delta(v)}{v} \in \K\llbracket z\rrbracket $, donc $0$ est un point singulier régulier de $L$. 

De plus, si $s \in \N$, $$L(z^s)=\delta z^s - \ell z^s - \gf{\delta(v)}{v} z^s=(s-\ell)z^s-\gf{\delta(v)}{v} z^s.$$ Mais comme $v$ est inversible dans $\K\llbracket z\rrbracket $, $\gf{1}{v} \in \K\llbracket z\rrbracket $ et $\gf{\delta(v)}{v}=z\gf{D(v)}{v} \in z\K\llbracket z\rrbracket $. D'où $\phi_0(x)=x-\ell$, de racine $\ell \in \mathbb{F}_p$. 

\item Soit $n \geqslant 2$, supposons le résultat vrai pour tout $m <n$. Alors comme $L$ est nilpotent, on peut écrire $L=L_1 L_2$, avec $L_1$ trivial d'ordre 1 sur $\K(z)$ et $L_2$ nilpotent d'ordre $n-1$ (car un opérateur trivial est produit d'opérateurs triviaux d'ordre 1). 
\end{itemize}
Par hypothèse de récurrence, $L_1$ et $L_2$ admettent $0$ comme point singulier régulier, donc, selon le point \textbf{b)} de la proposition \ref{exposantsproduitop}, $0$ est une singularité régulière de $L$. De plus, le point \textbf{c)} de la même proposition affirme que tout exposant de $L$ en $0$ est un exposant de $L_1$ en $0$ ou un exposant de $L_2$ en $0$. Par hypothèse de récurrence, les exposants de $L$ en $0$ sont donc dans~$\mathbb{F}_p$.\end{dem}

\begin{rqu}
Via un changement de variable $u=z-a$ ou $u=1/z$, on voit que si $L$ est nilpotent, tous les points de $\mathbb{P}^1(\K)$ sont singuliers réguliers avec exposants dans $\mathbb{F}_p$. En effet, on vérifie que si $M \in \K(z)\left[\mathrm{d}/\mathrm{d}z\right]$ et $a \in \mathbb{P}^1(\K)$, $M$ est trivial sur $\K(z)$ si et seulement si l'opérateur $M_a \in \K(u)\left[\mathrm{d}/\mathrm{d}u\right]$, obtenu par changement de variable $u=z-a$, est trivial sur $\K(u)$. Ainsi, selon la définition \ref{definilpotent} et le point \textbf{a)} de la proposition \ref{exposantsproduitop}, $L$ est nilpotent si et seulement si $L_a$ est nilpotent. Le théorème \ref{katzhonda} s'applique donc à $L_a$, donc $0$ est un point singulier régulier de $L_a$ à exposants dans $\mathbb{F}_p$, de sorte que $a$ est un point singulier régulier régulier de $L$ à exposants dans $\mathbb{F}_p$. 
\end{rqu}

\begin{Ex}
L'opérateur $L=D-1$ n'est pas nilpotent car $\infty$ n'est pas un point régulier.
\end{Ex}
\subsection{Cas global : le théorème de Katz} \label{subsec:thkatz}

Soit $\K$ un corps de nombres. Le but de cette partie est de démontrer, à l'aide du théorème de Katz-Honda, le théorème de Katz (théorème \ref{katz}) affirmant que les exposants en tout point de $\Qbar$ d'un opérateur différentiel $L \in \K(z)\left[\mathrm{d}/\mathrm{d}z\right]$ globalement nilpotent sont rationnels. Un opérateur globalement nilpotent est un opérateur qui définit un opérateur nilpotent quand il est localisé en presque tout premier de $\K$. Pour bien définir cette opération de localisation, il nous faut étudier la valeur absolue de Gauss associée à une valeur absolue non-archimédienne sur $\K$. 

\subsubsection{Valeur absolue de Gauss}

\begin{defi}
Soient $\K$ un corps, et $|\cdot |$ une valeur absolue non archimédienne sur $\K$. La \emph{valeur absolue de Gauss} associé à $|\cdot |$ sur $\K(z)$ est définie par 

$$ \left\vert \gf{\sum\limits_{i=1}^n a_i z^i}{\sum\limits_{j=1}^m b_j z^j} \right\vert_{\mathrm{Gauss}}=\gf{\max\limits_{1 \leqslant i \leqslant n} |a_i|}{\max\limits_{1 \leqslant j \leqslant m} |b_j|}. $$
\end{defi}

\begin{prop} \label{valuationgauss}
La valeur absolue de Gauss est bien définie et est une valeur absolue au sens usuel sur $\K(z)$.
\end{prop}

\begin{dem}
Commençons par définir la valeur absolue de Gauss sur $\K[z]$. 

\begin{itemize}
\item Les conditions d'homogénéité et d'inégalité triangulaire découlent immédiatement du fait que $| \cdot |$ est une valeur absolue sur $\K$.

\item Vérifions la condition de multiplicativité. Si $f(z)=\sum\limits_{i=0}^n a_i z^i$ et $g(z)=\sum\limits_{j=0}^m b_j z^j$, alors $$h(z)=f(z)g(z)=\sum\limits_{k=0}^{n+m} c_k z^k, \quad \text{où} \;\; c_k=\sum\limits_{i=0}^k a_i b_{k-i}.$$ On a donc pour tout $k$, comme $|\cdot|$ est non archimédienne, $$|c_k|\leqslant \max\limits_{0 \leqslant i \leqslant k} (|a_i| |b_{k-i}|) \leqslant |f|_{\mathrm{Gauss}} |g|_{\mathrm{Gauss}}.$$

Soient $i_0$ et $j_0$ tels que $|a_{i_0}|=|f|_{\mathrm{Gauss}}$, $|b_{j_0}|=|g|_{\mathrm{Gauss}}$, et $$\forall i < i_0, |a_i| < |a_{i_0}|, \qquad \forall j < j_0, |b_j| < |b_{j_0}|.$$ Alors 
$$c_{i_0+j_0}=a_0 b_{i_0+j_0}+ a_1 b_{i_0+j_0-1} + \dots +a_{i_0} b_{j_0} + \dots + a_{i_0+j_0} b_0.$$

Si $i<i_0$, $|a_i b_{i_0+j_0-i}| < |a_{i_0}| |b_{i_0+j_0-i}| \leqslant |a_{i_0}| | b_{j_0} |$ et de même si $i > i_0$. Donc $|c_{i_0+j_0}|=|a_{i_0} b_{j_0}|$, de sorte que $$|fg|_{\mathrm{Gauss}}=|f|_{\mathrm{Gauss}} |g|_{\mathrm{Gauss}}.$$

\end{itemize}

La condition de multiplicativité sur $\K[z]$ implique immédiatement que si $f=\gf{A}{B}=\gf{C}{D} \in~\K(z)$, avec $A, B, C, D \in \K[z]$, alors $|A|_{\mathrm{Gauss}}|D|_{\mathrm{Gauss}}=|C|_{\mathrm{Gauss}}|B|_{\mathrm{Gauss}}$, donc la valeur absolue de Gauss est bien définie sur $\K(z)$ et vérifie évidemment les conditions d'homogénéité et de multiplicativité. 

De plus, si $f=\gf{A}{B}$ et $g=\gf{C}{D}$ sont deux éléments de $\K(z)$, alors $f+g=\gf{AD+BC}{BD}$, donc 
\begin{align*}
|f+g|_{\mathrm{Gauss}} &= \gf{|AD+BC|_{\mathrm{Gauss}}}{|B|_{\mathrm{Gauss}} |D|_{\mathrm{Gauss}}} \leqslant \gf{\max(|A|_{\mathrm{Gauss}} |D|_{\mathrm{Gauss}}, |B|_{\mathrm{Gauss}} |C|_{\mathrm{Gauss}})}{|B|_{\mathrm{Gauss}} |D|_{\mathrm{Gauss}}} \\
&\leqslant  \max\left( \gf{|A|_{\mathrm{Gauss}}}{|B|_{\mathrm{Gauss}}}, \gf{|C|_{\mathrm{Gauss}}}{|D|_{\mathrm{Gauss}}} \right),
\end{align*} si bien que l'inégalité triangulaire est vérifiée.\end{dem}

\begin{rqu}
Soient $\K$ corps de nombres, $\mathfrak{p}$ idéal premier de $\Oal_\K$ et $|\cdot |_{\mathfrak{p}}$ la valeur absolue sur $\K$ associée. 
 
Soit $R=\left\{ f \in \K(z) \mid |f|_{\mathfrak{p}, \mathrm{Gauss}} \leqslant 1 \right\}$ l'anneau de valuation associé à $| \cdot |_{\mathfrak{p}, \mathrm{Gauss}}$. Montrons que $R=\Oal_\K[z]_{\mathfrak{p}[z]}$, le localisé de $\Oal_\K[z]$ en la partie multiplicative $\Oal_\K[z] \setminus \mathfrak{p}[z]$.

D'abord, remarquons que $R \cap \K[z]=\Oal_\K[z]$ et que $\mathfrak{p}[z]$ est un idéal premier de $\Oal_\K[z]$ car $\mathfrak{p}[z]=\left\{P \in \Oal_\K[z] : |P|_{\mathfrak{p}, \mathrm{Gauss}} <1 \right\}$.

De plus, si $$f=\gf{P}{Q}=\gf{\sum\limits_{i=1}^n a_i z^i}{\sum\limits_{j=1}^m b_j z^j} \in \K(z)$$ et $|f|_{\mathfrak{p},\mathrm{Gauss}} \leqslant 1$, alors en fixant $j_0$ tel que $|b_{j_0}|_{\mathfrak{p}}=|P|_{\mathfrak{p}, \mathrm{Gauss}}$, on peut écrire $f=\gf{\tilde{P}}{\tilde{Q}}$, avec $\tilde{P}, \tilde{Q} \in R$ et $|\tilde{Q}|_{\mathfrak{p},\mathrm{Gauss}}=1$, de sorte que $f \in \Oal_\K[z]_{\mathfrak{p}[z]}$.

Par conséquent, cela donne un sens à l'opération de réduction modulo $\mathfrak{p}$ :

$$ \Oal_\K[z]_{\mathfrak{p}[z]} \twoheadrightarrow \gf{\Oal_\K[z]_{\mathfrak{p}[z]}}{\mathfrak{p}[z] \Oal_\K[z]_{\mathfrak{p}[z]}} \simeq \left(\Oal_\K/\mathfrak{p} \right)(z), $$ ce dernier isomorphisme étant issu du morphisme surjectif de réduction modulo $\mathfrak{p}$ coefficient par coefficient  $$\fonction{\varphi_{\mathfrak{p}}}{\Oal_\K[z]_{\mathfrak{p}[z]}}{\left(\Oal_\K/\mathfrak{p} \right)(z)}{\gf{P}{Q}, Q \not\in \mathfrak{p}[z]}{\gf{\overline{P}}{\overline{Q}}}$$ de noyau $\mathfrak{p}[z] \Oal_\K[z]_{\mathfrak{p}[z]}$.
\end{rqu}

\subsubsection{Opérateurs globalement nilpotents et théorème de Katz}

Définissons à présent la notion d'opérateur globalement nilpotent. Soit $\K$ un corps de nombres et $L \in \K(z)\left[\gf{\mathrm{d}}{\mathrm{d}z}\right]$ de la forme $L=\delta^n + A_1 \delta^{n-1}+\dots+A_n$, où $A_i \in \K(z)$ et $\delta=z\gf{\mathrm{d}}{\mathrm{d}z}$.

Selon le théorème d'Ostrowski, les valuations non archimédiennes sur $\K$ sont équivalentes à $| \cdot |_{\mathfrak{p}}$ pour $\mathfrak{p}$ idéal premier de $\Oal_\K$. On peut donc définir la valuation de Gauss correspondante sur $\K(z)$, notée $| \cdot |_{\mathfrak{p}, \mathrm{Gauss}}$. On note
$$\mathcal{S}=\left\{ \mathfrak{p} \in \mathrm{Spec}(\Oal_\K) / \exists j \in \{ 1, \dots, n \} : | A_j |_{\mathfrak{p}, \mathrm{Gauss}} >1 \right\rbrace.$$
L'ensemble $\mathcal{S}$ est fini car pour $\alpha \in \K$, les premiers $\mathfrak{p}$ tels que $|\alpha|_{\mathfrak{p}} \neq 1$ sont les premiers intervenant dans la décomposition en facteurs premiers de l'idéal fractionnaire $(\alpha)$, qui sont en nombre fini.  

Si $\mathfrak{p} \not\in \mathcal{S}$, on note $\overline{\K_{\mathfrak{p}}}=\Oal_\K/\mathfrak{p}$ et $L_{\mathfrak{p}} \in \overline{\K_{\mathfrak{p}}}(z)\left[\gf{\mathrm{d}}{\mathrm{d}z}\right]$ la réduction de $L$ modulo $\mathfrak{p}$ :

$$ L_{\mathfrak{p}}=\delta^n+ \overline{A_1} \delta^{n-1} + \dots + \overline{A_n}. $$

\begin{defi}
Soit $T \subset \mathrm{Spec}(\Oal_\K)$. On dit que $T$ admet $d \geqslant 0$ pour \emph{densité de Dirichlet} lorsque $$ \gf{-1}{\log(s-1)} \sum_{\mathfrak{p} \in T} \gf{1}{N(\mathfrak{p})^s} \longrightarrow d$$ quand $s$ tend vers $1$ pour $s$ réel, $s >1$. \end{defi}

Notons en particulier qu'un ensemble fini est de densité de Dirichlet nulle. En annexe, d'autres notions de densité (naturelle et polaire) et leur lien avec la densité de Dirichlet sont présentées.

\begin{defi}
On dit qu'une propriété est valable pour \emph{presque tout} premier $\mathfrak{p}$ d'un corps de nombres (relativement à la densité de Dirichlet) si elle vaut pour un ensemble de premiers de densité de Dirichlet $1$. Ceci est vrai en particulier, si elle vaut pour tous les premiers sauf un nombre fini.
\end{defi}

\begin{defi} \label{defglobalementnilpotent}
L'opérateur $L$ est dit \emph{globalement nilpotent} lorsque pour presque tout $\mathfrak{p}$, $L_{\mathfrak{p}}$ est nilpotent.
\end{defi}

\begin{prop} \label{gaussserieformelle}
Considérons une valeur absolue $| \cdot |$ non archimédienne sur $\K$. Alors si $\xi \in \K(z)$ n'a pas de pôle dans $D_{|\cdot |}(0,1) \setminus \{ 0 \}$, où $| \cdot |$ désigne une extension fixée de la valeur absolue considérée à un corps de décomposition de $\xi$, et $\xi=\sum\limits_{j=-m}^{\infty} c_j z^j$ dans $\K((z))$, alors $|\xi|_{\mathrm{Gauss}}=\sup\limits_{j \geqslant -m} |c_j|$.
\end{prop}

\begin{dem}
\textbf{Étape 1} : Montrons d'abord qu'on peut se ramener au cas où $\xi=\gf{f(z)}{g(z)}$, avec $g(z)=1+a_1 z+\dots + a_{\ell} z^\ell$ et pour tout $j$, $|a_j| \leqslant 1$.

On a $\xi=P(z)/Q(z)$, avec $Q$ sans pôle dans $D_{| \cdot |}(0,1) \setminus \left\{0 \right\}$ par hypothèse. Donc dans un corps de décomposition de $Q$, on a $$Q(z)=rz^q (z-b_1) \dots (z-b_{\ell}),$$ où $q \in \N$ et $|b_j| \geqslant 1$. On peut donc écrire $\xi=\gf{f(z)}{g(z)}$, avec $$f(z)=\gf{(-1)^\ell P(z)}{r b_1 \dots b_{\ell} z^q} \quad \mathrm{et} \quad g(z)= \gf{(-1)^\ell (z-b_1) \dots (z-b_{\ell})}{b_1 \dots b_{\ell}} = 1 + a_1 z+ \dots + a_{\ell} z^\ell.$$

Si $j \in \{ 1, \dots, \ell \}$, $a_j$ s'exprime comme somme symétrique de $b_1, \dots, b_{\ell}$ :

$$ a_j=\gf{(-1)^j}{b_1 \dots b_{\ell}} \sum_{1 \leqslant i_1 < \dots < i_{\ell-j} \leqslant \ell} b_{i_1} \dots b_{i_{\ell-j}}=(-1)^j \sum_{i_1 < \dots < i_{\ell-j}} \prod_{k \not\in \left\{ i_1, \dots, i_{\ell-j} \right\}} \gf{1}{b_k}. $$

Il s'ensuit que les $a_j$ sont bien de valeur absolue inférieure ou égale à 1 car $| \cdot |$ est non-archimédienne. 

De plus, on peut supposer que $q=0$, car la multiplication par $z^{-q}$ est un simple décalage d'indice dans l'écriture en série formelle.

\bigskip

\textbf{Étape 2} : Supposons que $\xi=\gf{f(z)}{g(z)}$, avec $g(z)=1+a_1 z+\dots + a_{\ell} z^\ell$ et pour tout $j$, $|a_j| \leqslant 1$, de sorte que $|f|_{\mathrm{Gauss}}=|\xi|_{\mathrm{Gauss}}$. Écrivons $\xi=\sum\limits_{j=0}^{\infty} c_j z^j$. Alors en développant en série formelle, on a $$\xi=f(z) \sum\limits_{s=0}^{\infty} (-1)^s (a_1 z+ \dots + a_{\ell} z^\ell)^s=f(z) \sum\limits_{s=0}^{\infty} (-1)^s \sum\limits_{i_1+\dots+i_{\ell}=s} \binom{s}{i_1, \dots, i_{\ell}} a_1^{i_1} \dots a_{\ell}^{i_{\ell}} z^{i_1+2 i_2+\dots \ell i_{\ell}}$$

De plus, comme $| \cdot |$ est non archimédienne, par hypothèse sur les $a_j$, on a $$\forall (i_1, \dots, i_{\ell})\in \N^\ell, \;\; \left\vert \dbinom{s}{i_1, \dots, i_{\ell}} a_1^{i_1} \dots a_{\ell}^{i_{\ell}} \right\vert \leqslant 1,$$ d'où $\gf{1}{g(z)}=\sum\limits_{k=0}^{\infty} g_k z^k$, avec $|g_k| \leqslant 1$. En effectuant le produit avec $f(z)$ et en identifiant terme à terme, il vient que $\forall j \in \N, |c_j| \leqslant |f|_{\mathrm{Gauss}}$. En particulier, $\sup\limits_{j \in \N} |c_j| < \infty$ et $\sup\limits_{j \in \N} |c_j| \leqslant |f|_{\mathrm{Gauss}}$.

\medskip

Par ailleurs, \begin{align*}f(z) &=g(z) \sum\limits_{j=0}^{\infty} c_j z^j=\sum\limits_{j=0}^{\infty} c_j z^j+\sum\limits_{j=0}^{\infty} a_1 c_j z^{j+1} + \dots + \sum\limits_{j=0}^{\infty} a_{\ell} c_j z^{j+\ell} \\
&= \sum\limits_{j=0} (c_j+a_1 c_{j-1}+ \dots + a_{\ell} c_{j-\ell}) z^j\end{align*} en posant $c_j=0$ pour $j <0$. Donc $$|f|_{\mathrm{Gauss}} \leqslant \sup\limits_{j \in \N} (\max(|c_j|, |a_i c_{j-i}|, 1 \leqslant i \leqslant \ell)) \leqslant \sup\limits_{j \in \N} |c_j|.$$ Comme $|f|_{\mathrm{Gauss}} = |\xi|_{\mathrm{Gauss}}$, on a le résultat voulu.\end{dem}

On utilisera la conséquence suivante du théorème de Chebotarev (voir l'annexe pour la démonstration de cette proposition).

\begin{prop}\label{conseqchebokatz}
Soit $\K$ un corps de nombres et $\alpha \in \K$. Si pour tout premier $\mathfrak{p}$ de $\K$ dans un ensemble de densité $1$, $|\alpha|_{\mathfrak{p}} \leqslant 1$ et $\alpha \mod \mathfrak{p} \in \mathbb{F}_p$, où $(p)=\mathfrak{p} \cap \Z$, alors $\alpha \in \Q$. 
\end{prop}

Le théorème suivant est le résultat principal de cette sous-section (cf \cite[Theorem 8.1 p. 223]{Katzarticle})

\begin{Th}[Katz] \label{katz}
\begin{enumerateth}
\item Si $L_{\mathfrak{p}}$ est nilpotent pour un nombre infini de premiers $\mathfrak{p}$ de $\Oal_\K$, alors tout point de $\mathbb{P}^1(\K)$ est une singularité régulière.
\item Si $L$ est globalement nilpotent, alors le point \textbf{a)} est vrai et de plus les exposants de $L$ en tout point sont rationnels.
\end{enumerateth}
\end{Th}

\begin{dem}

\textbf{a)} Montrons l'assertion pour le point $0$ (par changement de variable, cela vaudra pour tous les points). 

Il existe un nombre fini de premiers $\mathfrak{p}$ tels que au moins un $A_j \mod \mathfrak{p}$ a un pôle dans $D_{|\cdot|_{\mathfrak{p}}}(0,1) \setminus~\{0\}$. En effet, étendons la valeur absolue $| \cdot |_{\mathfrak{p}}$ en une valeur absolue non archimédienne $| \cdot |_{\mathfrak{P}}$ sur l'extension finie $\mathbb{L}$ de $\K$ engendrée par les pôles dans $\Qbar$ de $A_1, \dots, A_n$. Écrivons pour tout $j$, $$A_j=U_j \prod\limits_{i=1}^{\ell_j} (z-\alpha_{i,j})^{-1}, \quad \text{avec} \;\; \alpha_{i,j} \in \mathbb{L},\;\; U_j \in \Qbar[z].$$ Alors il n'y a qu'un nombre fini de premiers $\mathfrak{p}$ de $\K$ tels qu'il existe $j$ vérifiant, pour tout $i \in \{ 1, \dots, \ell_j \}$, $|\alpha_{i,j}|_{\mathfrak{P}} <~1$.

Par conséquent, selon la proposition \ref{gaussserieformelle}, quitte à ajouter un nombre fini de premiers à l'ensemble $\mathcal{S}$, on a pour $\mathfrak{p} \not\in \mathcal{S}$,

$$ \forall j \in \{ 1, \dots, 	n \}, \quad |A_j|_{\mathfrak{p}, \mathrm{Gauss}}=\sup\limits_{1 \leqslant \ell \leqslant m \atop \ell \in \N} (|A_{j,k}|_{\mathfrak{p}, \mathrm{Gauss}}, |B_{j,\ell}|_{\mathfrak{p}, \mathrm{Gauss}}) \leqslant 1,$$ où $A_j=\sum\limits_{k=1}^{m} \gf{A_{j,k}}{z^k}+\sum\limits_{\ell=0}^{\infty} B_{j,\ell} z^\ell$. Si $\mathfrak{p} \not\in \mathcal{S}$ et $L_{\mathfrak{p}}$ est nilpotent alors par le théorème \ref{katzhonda}, pour tout $j \in \{ 1, \dots, n \}$, $A_j \mod \mathfrak{p}$ n'a pas de pôle en $0$. Comme les opérations de réduction modulo $\mathfrak{p}$  dans $\K(z)$ et dans $\K((z))$ coïncident, on a pour tout $k \in \{ 1, \dots, m \}$, $A_{j,k} \mod \mathfrak{p}=0$. Puisque cela vaut pour une infinité de premiers $\mathfrak{p}$, on a pour tout $k$, $A_{j,k}=0$, donc $A_j \in \K\llbracket z\rrbracket$, de sorte que $0$ est une singularité régulière pour $L$. 

\bigskip

\textbf{b)} Supposons que $L$ est globalement nilpotent. Selon le premier point, $0$ est une singularité régulière de $L$. Soit $\phi_0(x) \in \K[x]$ le polynôme indiciel de $L$ en 0.

Si pour tout $j \in \{ 1, \dots n \}, A_j = \sum\limits_{k=0}^{\infty} A_{j,k} z^k$, alors la définition \ref{defexposants} donne $$\phi_0(x)=x^n+\sum\limits_{j=1}^{n} A_{j,0} x^{n-j}.$$ Cette expression montre que si $\mathfrak{p} \not\in \mathcal{S}$, le polynôme indiciel $\phi_{0, \mathfrak{p}}$ de $L_{\mathfrak{p}}$ est la réduction modulo $\mathfrak{p}$ de $\phi_0$. 

Soit $\alpha \in \K$ une racine de $\phi_0$. Quitte à ajouter un nombre fini de premiers à $\mathcal{S}$, on peut supposer que $|\alpha |_{\mathfrak{p}} \leqslant 1$. De plus, $0=\phi_0(\alpha) \mod \mathfrak{p}=\phi_{0, \mathfrak{p}}(\alpha \mod \mathfrak{p})$, de sorte que, par le théorème \ref{katzhonda}, $\alpha \mod \mathfrak{p} \in \mathbb{F}_p$, où $(p)=\mathfrak{p} \cap \Z$. Ainsi, puisque cela vaut pour un ensemble de premiers de densité 1, par la proposition \ref{conseqchebokatz}, on a $\alpha \in \Q$. Ceci conclut la preuve de \textbf{b)}.\end{dem}

\begin{rqus}
\begin{enumerateth}[label=\textbf{\textit{\roman*)}}]
\item En réalité, si $L \in \K(z)\left[\gf{\mathrm{d}}{\mathrm{d}z}\right]$ est globalement nilpotent, alors les exposants de $L$ en tout point de $\Qbar$ (et non pas seulement de $\K$) sont rationnels.

En effet, soit $\alpha \in \Qbar$. Notons $\K'=\K(\alpha)$, et prenons $\mathfrak{p}' \in \Spec(\Oal_{\K'})$. Si $\mathfrak{p}=\mathfrak{p}' \cap \K$, alors comme $L$ est à coefficients dans $ \K(z)$, on a $L_{\mathfrak{p'}}=L_{\mathfrak{p}}$.

De plus, si $T$ est l'ensemble de densité nulle (prenons par exemple la densité de Dirichlet, voir définition \ref{densitedirichlet}) des premiers $\mathfrak{p} \in \Spec(\Oal_\K)$ tels que $L_{\mathfrak{p}}$ n'est pas nilpotent, et $T'$ est l'ensemble des premiers $\mathfrak{p}' \in \Spec(\Oal_{\K'})$ au dessus des éléments de $T$, alors

\begin{align*}\forall s >1, \quad \gf{-1}{\log(s-1)}\sum\limits_{\mathfrak{p}' \in T'} \gf{1}{N_{\K'/\Q}(\mathfrak{p}')^s} &\leqslant \gf{-1}{\log(s-1)}\sum\limits_{\mathfrak{p} \in T} \sum\limits_{\mathfrak{p}' \mid \mathfrak{p}} \gf{1}{N_{\K/\Q}(\mathfrak{p})^{s}} \\ &\leqslant \gf{-[\K':\K]}{\log(s-1)} \sum\limits_{\mathfrak{p} \in T} \gf{1}{N_{\K/\Q}(\mathfrak{p})^{[\K':\K] s}}\end{align*} car il y a au plus $[\K':\K]$ premiers au-dessus de $\mathfrak{p} \in T$. Donc comme le terme de droite tend vers $0$ quand $s$ tend vers $1$ par valeurs supérieures, il en va de même du terme de gauche, si bien que $T'$ est de densité de Dirichlet nulle. Or, $T'$ contient tous les premiers $\mathfrak{p'}$ de $\K'$ tels que $L_{\mathfrak{p'}}$ n'est pas nilpotent, donc $L \in \K'(z)[\mathrm{d}/\mathrm{d}z]$ est globalement nilpotent. Selon le théorème de Katz, $\alpha \in \K'$ est une singularité régulière de $L$ en laquelle ses exposants sont rationnels.

    \item Comme le remarque Dwork \cite[p. 100]{Dwork}, le théorème de Katz peut encore être raffiné de la manière suivante : soit $L$ un opérateur différentiel et $\mathcal{S}$ un ensemble de premiers de $\K$ tel que
\begin{itemize}
    \item L'ensemble de premiers de $\Z$ $$\mathcal{S}_{\Z} := \{ \mathfrak{p} \cap \Z, \mathfrak{p} \in \mathcal{S} \}$$ est de densité de Dirichlet  strictement supérieure à $\gf{1}{2}$.
    \item L'opérateur $L_{\mathfrak{p}}$ est nilpotent pour tout premier $\mathfrak{p}$ dans $\mathcal{S}$.
\end{itemize} Alors le théorème de Katz est vérifié pour un tel $L$. Ceci est une conséquence du corollaire \ref{corochebo1} de l'annexe.

\end{enumerateth}
\end{rqus}

\newpage

\section{$G$-fonctions et condition de Galochkin} \label{sec:GfonsChud}

\subsection{Présentation et exemples} \label{subsec:exempleGfons}

Le but de cette section est d'introduire la définition des $G$-fonctions de Siegel \cite{Siegelarticle}. Elles constituent une classe de fonctions spéciales incluant de nombreux exemples.

\begin{defi} \label{maisonetdenominateur}
Soit $\alpha \in \Qbar$. La \emph{maison} de $\alpha$ est le maximum des modules des conjugués (au sens de Galois) de $\alpha$, on la note $\house{\alpha}$.
	
Le \emph{dénominateur} de $\alpha$ est $\mathrm{den}(\alpha)$, plus petit entier $d \in \N^*$ tel que $d \alpha$ est un entier algébrique. 

Le dénominateur de $\alpha_1, \dots, \alpha_n \in \Qbar$ est $\mathrm{den}(\alpha_1, \dots, \alpha_n)=\mathrm{ppcm}(\mathrm{den}(\alpha_1), \dots, \mathrm{den}(\alpha_n))$.
\end{defi}

\begin{defi}\label{defgfonction}
Une $G$-fonction est une série $f(z)=\sum\limits_{n=0}^{\infty} a_n z^n \in \Qbar\llbracket z\rrbracket $ telle que

\begin{enumerateth}[label=\textbf{\alph*)}]
\item $f$ est solution d'une équation différentielle linéaire à coefficients dans $\Qbar(z)$.
\item Il existe une constante $C_1 >0$ telle que $\forall n \in \N, \house{a_n} \leqslant C_1^{n+1}$.
\item  Il existe une constante $C_2 >0$ telle que $\forall n \in \N, \mathrm{den}(a_0, \dots,a_n) \leqslant C_2^{n+1}$.
\end{enumerateth}
\end{defi}

\begin{rqu}
En réalité, la condition de croissance géométrique sur les coefficients de $f$ implique que $f$ est une fonction holomorphe autour de $0$.
\end{rqu}

Cette définition concerne une large classe de fonctions \og{}intéressantes \fg{} pour lesquelles on recherche des résultats de transcendance.

\begin{prop} \label{exemplesgfonctions}
Les fonctions suivantes sont des $G$-fonctions :

\begin{enumerateth}
\item La fonction logarithme complexe $-\log(1-z)=\sum\limits_{n=0}^{\infty} \gf{z^n}{n}$, et plus généralement, les fonctions polylogarithmes $\mathrm{Li}_s(z)=\sum\limits_{n=0}^{\infty} \gf{z^n}{n^s}$ pour $s \in \N^*$.

\item Toute fonction $f \in \Qbar\llbracket z\rrbracket $ algébrique sur $\Qbar(z)$ et holomorphe en $0$.

\item Les fonctions hypergéométriques : si  $\boldsymbol{\alpha}=(\alpha_1, \dots, \alpha_n) \in \Q^n$ et $\boldsymbol{\beta}=(\beta_1, \dots, \beta_n) \in (\Q \setminus \Z_{-})^n$, $$_n F_{n-1} (\boldsymbol{\alpha} ; \boldsymbol{\beta} ;z)=\sum\limits_{m=0}^{\infty} \gf{(\alpha_1)_m \dots (\alpha_n)_m}{(\beta_1)_m \dots (\beta_{n-1})_m m!} z^m,$$  où pour $x \in \C$ et $m \geqslant 1$, $(x)_m=x(x+1)(x+2) \dots (x+m-1)$, et $(x)_0=1$.
\end{enumerateth}
\end{prop}

On aura besoin pour la preuve du point \textbf{b)} du théorème suivant, dont une preuve peut être trouvée dans \cite[pp. 28-30]{Cassels}.

\begin{Th}[Eisenstein] \label{eisenstein}
Soit $y(z) \in \C\{ z \}$ algébrique sur $\Qbar(z)$. Alors il existe $c \in \N^*$ tel que $y(cz) \in \Oal\llbracket z\rrbracket $. 
\end{Th}
 
Pour la preuve du point \textbf{c)}, un point clef est le résultat suivant prouvé dans \cite[p. 57]{Siegel}.

\begin{prop} \label{prop:siegelcroissancepochhamer}
Si $u, v \in \Q$ et $v \not\in \Z_{-}^*$, alors le dénominateur commun des $\gf{(u)_k}{(v)_k}$, $k \in \{ 0, \dots, n \}$ a une croissance géométrique avec $n$. \end{prop}

\begin{dem}
\begin{itemize}
\item On a $((1-z)\mathrm{Li}'_1(z))'=0$ et si $s \geqslant 2, z\mathrm{Li}'_s(z)=\mathrm{Li}_{s-1}(z)$, donc on obtient par récurrence une équation différentielle à coefficients dans $\Qbar(z)$ satisfaite par $\mathrm{Li}_s$. 

La condition sur la taille des coefficients de $\mathrm{Li}_s$ étant trivialement satisfaite, il reste à vérifier celle sur le dénominateur. Soit, pour $n \in \N^*$, $d_n=\mathrm{ppcm}(1, 2, \dots, n)$. Alors $$\mathrm{den}\left(1, \gf{1}{2^s}, \dots, \gf{1}{n^s} \right)=\mathrm{ppcm}(1, 2^s, \dots, n^s),$$ qui divise $d_n^s$. Mais si $p$ est un premier inférieur à $n$, et $k \in \N$, alors $$p^k \leqslant n \ssi k \log p \leqslant \log n \ssi k \leqslant \left\lfloor \frac{\log n}{\log p} \right\rfloor,$$ donc $d_n \leqslant \prod\limits_{p \in \mathcal{P} \atop p \leqslant n} p^{\left\lfloor \frac{\log n}{\log p} \right\rfloor}$, d'où $$\log D_n \leqslant \sum\limits_{p \leqslant n} \left\lfloor \frac{\log n}{\log p} \right\rfloor \log p \leqslant \sum\limits_{p \leqslant n} \log n=(\log n) \pi(n),$$ où $\pi$ est la fonction de comptage des nombres premiers. Mais selon le théorème des nombres premiers, $\pi(n) \sim n \log n$, donc on peut trouver $c >0$ tel que $\log d_n \leqslant cn$, de sorte que $d_n \leqslant e^{cn}$, ce qui est la condition voulue.

\bigskip
\item Soit $f(z) \in \Qbar\llbracket z\rrbracket $ algébrique sur $\Qbar(z)$ et holomorphe en $0$. On admet que $f$ est solution d'une équation différentielle à coefficients dans $\Qbar(z)$.

Si $f(z)=\sum\limits_{n=0}^{\infty} a_n z^n$, alors comme $f(z)$ est holomorphe en $0$, par la formule de Hadamard sur le rayon de convergence, on peut trouver $C_1>0$ tel que $\forall n \in \N, |a_n| \leqslant C_1^{n+1}$. De plus, $f$ vérifie $p_d(z) f^{(d)}(z)+ \dots + p_0(z)=0$, avec $p_i(z) \in \Qbar[z]$, de sorte qu'en identifiant de part et d'autre de l'égalité les coefficients du développement en série entière, on obtient, pour tout $n \in \N$, $a_n \in \K$, où $\K$ est le corps de nombres engendré sur $\Q$ par les coefficients des $p_i(z)$ et $a_0, \dots, a_d$. 

Quitte à remplacer $\K$ par sa clôture galoisienne, on peut supposer que $\K/\Q$ est normale. Soit $\sigma \in \Gal(\K/\Q)$. En appliquant le morphisme d'anneaux (injectif car $\sigma$ l'est)
$$\fonction{\varphi_{\sigma}}{\K\llbracket z\rrbracket }{\K\llbracket z\rrbracket }{\sum\limits_{n=0}^{\infty} b_n z^n}{\sum\limits_{n=0}^{\infty} \sigma(b_n) z^n} $$ à l'équation algébrique satisfaite par $f$ sur $\Qbar(z)$, on obtient que $\varphi_{\sigma}(f)$ est algébrique sur $\Qbar(z)$. Par suite, on peut trouver une constante $C_{\sigma} >0$ telle que $\forall n \in \N, |\sigma(a_n)| \leqslant C_{\sigma}^{n+1}$. Cela nous fournit une constante $C_2 >0$ telle que, pour tout $n \in \N$, $\house{a_n} \leqslant C_2^{n+1}$.

\medskip

La condition de croissance géométrique sur les dénominateurs découle du théorème \ref{eisenstein}. En effet, fixons $c$ tel que $f(cz) \in \Oal[z]$, alors si $n \in \N$, $c^n a_0, \dots, c^n a_n \in \Oal$, donc $\mathrm{den}(a_0, \dots, a_n) \leqslant c^n$.

 \bigskip
\item L'existence d'une équation différentielle à coefficients dans $\Qbar(z)$ satisfaite par une série hypergéométrique est expliquée plus loin, voir l'équation \eqref{eqdiffhypergeom} après la proposition \ref{galochkinhypergeom}. De plus, proposition \ref{prop:siegelcroissancepochhamer} nous assure que la condition sur les dénominateurs est satisfaite.

On vérifie pour finir la condition sur la taille des coefficients. Si $\alpha=\gf{p}{q}, \beta=\gf{u}{v} \in \Q$, on a pour tout $n \in \N$, $(\alpha)_n=\gf{1}{q^n} p(p+q) \dots (p+(n-1)q)$ et de même $(\beta)_n=\gf{1}{v^n} u(u+v) \dots (u+(n-1)v)$, de sorte qu'on peut supposer $p$ et $v$ positifs puisque $p+kq \geqslant 0$ pour $k$ assez grand. 

De plus, si $m \in \N$ est tel que $p \leqslant mq$, 

$$p(p+q) \dots (p+(n-1)q) \leqslant n! \times \binom{m+n}{m} q^{n-1} \leqslant n! 2^{m+n} q^n.$$

Puisque $u \geqslant 0$, on a par ailleurs $u(u+v) \dots (u+(n-1)v) \geqslant v \times (2v) \times \dots \times ((n-1)v) \geqslant v^{n-1} (n-1)!$, ce qui donne l'existence d'une constante $c >0$ telle que $$ \forall n \in \N^*, \;\; \gf{(\alpha)_n}{(\beta)_n} \leqslant c 2^n.$$ \end{itemize} 
\noindent C'est bien la borne géométrique désirée, ce qui conclut la preuve du point \textbf{c)}.\end{dem}

\begin{rqu}
Le point \textbf{b)} de la proposition \ref{exemplesgfonctions} exhibe une classe importante de $G$-fonctions prenant une infinité de valeurs algébriques à des points algébriques.

On peut également citer les deux exemples suivants (cf \cite[p. 19]{Beukers}) :

$$ _2 F_1 \left(\gf{1}{12}, \gf{5}{12} ; \gf{1}{2} ; \gf{1323}{1331} \right)=\gf{3}{4} (11)^{\frac{1}{4}}. $$

$$ _2 F_1 \left(\gf{1}{12}, \gf{7}{12} ; \gf{2}{3} ; \gf{64000}{64009} \right)=\gf{2}{3} (253)^{\frac{1}{6}}. $$
\end{rqu}

\bigskip

Tous les exemples \og{}célèbres\fg{} de $G$-fonctions sont en réalité à coefficients rationnels, en particuliers ceux qui sont issus de la physique. A ce sujet, on pourra consulter \cite{Andrews}.

\begin{Ex}
Les fonctions suivantes ne sont pas des $G$-fonctions :
\begin{itemize}
\item $f(z)=\sum\limits_{n=0}^{\infty} n! z^n$, car $(n!)_{n \in \N}$ ne peut être bornée par une suite géométrique.
\item les $E$-fonctions qui ne sont pas polynomiales. On rappelle qu'une $E$-fonction est une série $f(z)=\sum\limits_{n=0}^{\infty} \gf{b_n}{n!} z^n$ vérifiant la condition \textbf{a)} de la définition \ref{defgfonction} et où $(b_n)_{n \in \N}$ vérifie les conditions \textbf{b)} et \textbf{c)} de la définition \ref{defgfonction}. 

En effet, soit $\K$ un corps de nombres et $f(z)=\sum\limits_{n=0}^{\infty} \gf{b_n}{n!} z^n=\sum\limits_{n=0}^{\infty} a_n z^n \in \K\llbracket z\rrbracket $ qui est à la fois une $G$-fonction et une $E$-fonction. L'existence d'un tel corps $\K$ est assurée par celle d'une équation différentielle à coefficients dans $\Qbar(z)$ satisfaite par $f$. On peut donc trouver $C>0$ tel que $\forall n \in \N, \house{b_n} \leqslant C^{n+1}$ et $d_n \leqslant C^{n+1}$, où $d_n a_n \in \Oal_\K, d_n \in \N$. Donc en passant à la norme, on peut trouver $\tilde{C} >0$ telle que

$$\forall n \in \N, \quad \gf{\tilde{C}^{2n+2}}{n!} \geqslant \gf{\tilde{C}^{n+1} d_n}{n!} \geqslant \left\vert N_{\K/\Q} \left(\gf{d_n b_n}{n!} \right) \right\vert.$$

La suite d'entiers $\left(N_{\K/\Q} \left(\gf{d_n b_n}{n!} \right)\right)_{n \in \N}$ tend donc vers $0$, donc est nulle à partir d'un certain rang $n_0$, si bien que $\forall n \geqslant n_0, b_n=0$. Ainsi, $f(z) \in \K[z]$. Le résultat voulu s'ensuit par contraposée. 
\end{itemize}
\end{Ex}

\medskip

\bigskip

Dans le cas hypergéométrique, la caractérisation suivante a été prouvée par Galochkin dans \cite[p. 8]{Galochkin81}.

\begin{prop}[Galochkin] \label{galochkinhypergeom}
Soit $F(z)=_n F_{n-1} (\boldsymbol{\alpha}, \boldsymbol{\beta} ;z)$, avec $\boldsymbol{\alpha} \in (\C \setminus \Z_{-})^n$, $\boldsymbol{\beta} \in (\C \setminus \Z_{-})^{n-1}$ tels que $\forall i,j, \alpha_i \neq \beta_j$. Alors $F$ est une $G$-fonction si et seulement les deux conditions suivantes sont vérifiées :

\begin{itemizeth}
\item $\boldsymbol{\alpha} \in \Qbar^n$ et  $\boldsymbol{\beta} \in \Qbar^{n-1}$ ;
\item Les $\alpha_i$ et $\beta_j$ qui ne sont pas rationnels peuvent être regroupées en paires $(\alpha_i, \beta_j)$ telles que $\alpha_i - \beta_j \in \N$.
\end{itemizeth}
\end{prop}

En particulier, dès que $\boldsymbol{\alpha} \in \Q^n$ et $\boldsymbol{\beta} \in \Q^{n-1}$, on retrouve le fait que $_n F_{n-1} (\boldsymbol{\alpha}, \boldsymbol{\beta} ;z)$ est une $G$-fonction. Rappelons (cf \cite[p. 3]{Beukershypergeom}) que l'équation différentielle d'ordre $n$ satisfaite par $F(z)=_n~F_{n-1} (\boldsymbol{\alpha}, \boldsymbol{\beta} ;z)$ est \begin{equation} \label{eqdiffhypergeom} \left[\delta(\delta+\beta_1-1)\dots(\delta+\beta_{n-1}-1)-z(\delta+\alpha_1)\dots (\delta+\alpha_n) \right](y(z))=0,\end{equation} avec $\delta=z\gf{\mathrm{d}}{\mathrm{d}z}$. Ses singularités sont $0$, $1$ et $\infty$ et les exposants 

\begin{itemize}
\item en $0$ sont $0$, $1-\beta_1$, \dots, $1-\beta_{n-1}$ ;
\item en $1$ sont $0$, $1$, \dots, $n-2$, $-\alpha_n+\sum\limits_{i=1}^{n-1} \beta_i- \alpha_i$ ;

\item en $\infty$ sont $\alpha_1$, \dots, $\alpha_n$.
\end{itemize}

En substance, l'exemple suivant illustre le seul cas générique de fonction hypergéométrique de paramètres non rationnels qui est une $G$-fonction (comme somme de $G$-fonctions) : 

$$_2 F_1\left[\sqrt{2}+1, \gf{1}{2} ; \sqrt{2} ; z \right]=\sum\limits_{n=0}^{\infty} \gf{\left(\frac{1}{2}\right)_n}{n!} z^n+\gf{1}{\sqrt{2}} \sum\limits_{n=0}^{\infty} n \gf{\left(\frac{1}{2}\right)_n}{n!} z^n$$

car $\gf{(\sqrt{2}+1)_n}{(\sqrt{2})_n}=1+\gf{n}{\sqrt{2}}$. 

On peut voir \og{} à la main \fg{} que l'équation différentielle minimale sur $\Qbar(z)$ satisfaite par cette fonction est fuchsienne à exposants rationnels. En effet, si $a \in \Qbar$ et $a > 1$, l'équation différentielle minimale sur $\Qbar(z)$ vérifiée par $$f_0(z) :=_2 F_1\left[a+1, \gf{1}{2} ; a ; z \right]=\gf{1}{\sqrt{1-z}}+\gf{z}{2a (1-z)^{3/2}}$$ est $$2(1-z)((1-2a)z+2a)f'_0(z)=((1-2a)z+2a+2)f_0(z).$$ Ses points singuliers sont :

\begin{itemize}
\item $1$ avec pour exposant $-\gf{3}{2}$ ;
\item $\infty$ avec pour exposant $\gf{1}{2}$ :
\item $\gf{2a}{2a-1} \in ]1,2[$ qui est une singularité apparente puisque $(f_0)$ constitue une base de solutions holomorphes au voisinage de ce point.
\end{itemize}

\medskip

 Dans certains cas, on peut prouver en étudiant les exposants de l'équation \eqref{eqdiffhypergeom} qu'une fonction hypergéométrique n'est pas une $G$-fonction sans utiliser la proposition \ref{galochkinhypergeom}. L'exemple suivant illustre cette technique et l'importance de la condition $\alpha_i - \beta_j \in \N$ dans la proposition \ref{galochkinhypergeom}.

\begin{Ex}

Montrons que $$f_1(z):=\sum\limits_{n=0}^{\infty} \gf{z^n}{n+\sqrt{2}}=\gf{1}{\sqrt{2}} {}_2 F_1\left[\sqrt{2}, 1 ; \sqrt{2}+1 ; z \right] $$ n'est pas une $G$-fonction à l'aide du théorème d'André-Chudnovsky-Katz (théorème \ref{katzchudandre}). 

Supposons que l'équation différentielle minimale satisfaite par $f_1$ sur $\Qbar(z)$ soit d'ordre 1, donc de la forme 
\begin{equation}\label{EDordre1ex} y'(z)=p(z)y(z), p(z) \in \Qbar(z). \end{equation} Selon le théorème \ref{katzchudandre}, cette équation est fuchsienne, donc, selon \cite[pp. 145--148]{Hille}, la solution générale de \eqref{EDordre1ex} est $y(z)=\exp\left(\int_{[z_0,z]} p(s) \mathrm{d}s\right)$, et $p(z)$ n'a que des pôles simples ; de plus, on peut écrire 

$$ p(z)=\sum_{j=1}^m \gf{r_j}{z-a_j}, $$ où $a_j \in \Qbar$ et $r_j$ est l'exposant de \eqref{EDordre1ex} en $a_j$. Toujours selon le théorème \ref{katzchudandre}, on sait que $\forall j \in \{1, \dots, m \}, r_j \in \Q$. Par conséquent, en intégrant, on a  $$f_1(z)=C_0 \prod\limits_{j=1}^m (z-a_j)^{r_j},$$ avec $C_0 \in \C^*$ constante, si bien qu'il existe $q \in \N^*$ tel que $f_1(z)^q \in \Qbar(z)$ (on a bien $C_0 \in \Qbar$ car $f_1(z) \in \Qbar\llbracket z\rrbracket $ et les $a_j$ sont algébriques) : $f_1$ est algébrique sur $\Qbar(z)$.

Selon le théorème d'Eisenstein (théorème \ref{eisenstein}), il existe donc $C \in \N^*$ tel que pour tout $n \in \N$, $\gf{C^n}{n+\sqrt{2}} \in \Oal_{\Q(\sqrt{2})}=\Z[\sqrt{2}]$. Donc en passant à la norme dans $\Q(\sqrt{2})$, on obtient $$ \forall n \in \N, \gf{C^{2n}}{n^2-2} \in \Z,$$ c'est-à-dire que $n^2-2$ divise $C^{2n}$ pour tout $n$. Ainsi, tout facteur premier de $n^2-2$ est également un facteur premier de $C$. On va montrer que ce n'est pas possible.

En effet, si $p$ est un premier de $\Z$, il existe $n \in \N$ tel que $p$ divise $n^2-2$ si et seulement si $2$ est un carré dans $\mathbb{F}_p$, ce qui est équivalent, selon la loi de réciprocité quadratique faible, à $p \equiv \pm 1 [8]$. Le théorème de Dirichlet nous assurant qu'il existe une infinité de tels premiers, fixons en un qui ne fait pas partie des diviseurs de $C$ et considérons $n \in \N$ tel que $p$ divise $n^2-2$. On obtient ainsi une contradiction.

Donc l'équation différentielle minimale de $f_1$ est d'ordre 2, et (cf \eqref{eqdiffhypergeom} ci-dessus) cette équation est $$\left[\delta(\delta+\sqrt{2})-z(\delta+\sqrt{2})(\delta+1)\right](f_1(z))=0. $$ Ses exposants en $0$ sont $0$ et $-\sqrt{2} \not\in \Q$. Donc par la contraposée du théorème \ref{katzchudandre}, $f_1$ n'est pas une $G$-fonction.

\end{Ex}

\bigskip

Pour finir, citons ce résultat de structure sur les $G$-fonctions, déjà utilisé dans l'exemple de $_2 F_1\left[\sqrt{2}+1, \gf{1}{2} ; \sqrt{2} ; z \right]$ ci-dessus.

\begin{prop}
Les $G$-fonctions forment un sous-anneau de $\Qbar\llbracket z\rrbracket $ stable par dérivation et primitivation (avec constante d'intégration algébrique) et même une $\Qbar$-algèbre différentielle.
\end{prop}

\subsection{Le théorème des Chudnovsky}\label{subsec:thchudnovsky}

Soit $\mathbf{f}={}^t\left(
f_1(z), \dots, f_n(z) \right)\in \Qbar\llbracket z\rrbracket ^n$ vérifiant $\mathbf{f}'=G\mathbf{f}$, avec $G \in \mathcal{M}_n(\Qbar(z))$. Soit $G_s \in \mathcal{M}_n(\Qbar(z))$ la matrice telle que $\mathbf{f}^{(s)}=G_s \mathbf{f}$. Un raisonnement par récurrence montre que les $G_s$, $s \in \N$, sont liées par la relation 

\begin{equation}\label{recurrenceGs}
G_{s+1}=G_s G+G'_s,
\end{equation} où $G'_s$ désigne la matrice $G_s$ dérivée coefficient par coefficient. On prend $T(z) \in \Qbar[z]$ le plus petit dénominateur commun de tous les coefficients de la matrice $G(z)$. 
\begin{prop}
Pour tout $s\in \N^*$, $T^s G_s \in \mathcal{M}_n(\Qbar[z])$.
\end{prop}

\begin{dem}
On procède par récurrence sur $s$ :
\begin{itemize}
\item C'est la définition de $T$ pour $s=1$.

\item Soit $s \in \N^*$, supposons que $T^s G_s \in \mathcal{M}_n(\Qbar[z])$. Alors en utilisant l'équation \eqref{recurrenceGs},$$T^{s+1} G_{s+1}=(T^s G_s)(T G)+T^{s+1} G'_s.$$\end{itemize}
\noindent Or, par hypothèse de récurrence, $(T^s G_s)'$ est une matrice à coefficients polynomiaux et $$(T^s G_s)'=sT' T^{s-1} G_s+T^s G'_s,$$ donc en multipliant cette équation par $T$, on a $T^{s+1} G'_s \in \mathcal{M}_n(\Qbar[z])$, de sorte que $T^{s+1} G_{s+1} \in \mathcal{M}_n(\Qbar[z])$, ce qui conclut la preuve.\end{dem}

\subsubsection{Condition de Galochkin}

\begin{defi} \label{galochkin}
On note, pour $s \in \N$, $q_s$ le plus petit dénominateur supérieur ou égal à $1$ de tous les coefficients des coefficients des matrices $T(z)^m \gf{G_m(z)}{m !}$, quand $m \in \{ 1, \dots, s \}$. On dit que le système $y'=Gy$ \textit{vérifie la condition de Galochkin} si $$\exists C >0 : \; \forall s \in \N, \quad q_s \leqslant C^{s+1}.$$
\end{defi}

\begin{defi}
Soit $L \in \Qbar(z)\left[\gf{\mathrm{d}}{\mathrm{d}z}\right]$. On dit que $L$ est un \emph{$G$-opérateur} si la matrice compagnon de $L$ vérifie la condition de Galochkin.
\end{defi}

La terminologie de \emph{$G$-opérateur} est justifiée par la proposition suivante :

\begin{prop} \label{prop:baseGfonctionspointordinaire}
Soit $L \in \Qbar(z)\left[\mathrm{d}/\mathrm{d}z\right]$ un $G$-opérateur non nul d'ordre $\mu$. Soit $\alpha \in \Qbar$ un point ordinaire de $L$, alors il existe une base de solutions de l'équation $L(y(z))=0$ autour de $\alpha$ de la forme $(f_1(z-\alpha), \dots, f_{\mu}(z-\alpha))$, où les $f_i(u)$ sont des $G$-fonctions.
\end{prop}

\begin{dem}
Notons $G(z)=A_L(z)$ la matrice compagnon de $L$. Comme $\alpha$ est un point ordinaire, on sait qu'il existe au voisinage de $\alpha$ une base de solutions $(f_1(z-\alpha), \dots, f_{\mu}(z-\alpha))$ de l'équation $L(y(z))=0$, où les $f_i$ sont holomorphes au voisinage de $0$. On sait aussi que la matrice wronskienne de cette base $Y(z) \in \mathcal{M}_n(\Qbar\llbracket z\rrbracket)$ a un rayon de convergence non nul et est telle que $Y(\alpha) \in \GL_n(\Qbar)$ et $Y'(z)=G(z)Y(z)$, de sorte que $Y^{(s)}(z)=G_s(z) Y(z)$ pour tout entier $s$.

D'où \begin{equation}\label{eq:YetGs} Y(z)=\sum\limits_{n=0}^{\infty} Y^{(n)}(\alpha)(z-\alpha)^n=\left(\sum\limits_{n=0}^{\infty}  \gf{G_n(\alpha)}{n!}(z-\alpha)^n\right)Y(\alpha).\end{equation} Puisque $\alpha$ est un point ordinaire, $G(z)$ n'a pas de pôle en $\alpha$ et la condition de Galochkin implique qu'il existe une suite d'entiers positifs $(q_n)_{n \in \N}$ telle que $$\forall n \in \N,\; \forall k \leqslant n, \;\; q_n \gf{G_k(\alpha)}{k!} \in \Oal_{\Qbar} \quad \mathrm{et} \quad \exists C_1>0, \;\; \forall n \in \N, q_n \leqslant C_1^{n+1}.$$ Ainsi, selon \eqref{eq:YetGs}, on a $\forall n \in \N, \mathrm{den}(Y(\alpha)) q_n Y^{(n)}(\alpha) \in \mathcal{O}_{\Qbar}$. Ainsi, les $f_i(u)$ vérifient la condition \textbf{c)} de la définition \ref{defgfonction}. 

Par ailleurs, soit $\K$ un corps de nombres galoisien contenant $\alpha$ tel que $L \in \K(z)\left[\gf{\mathrm{d}}{\mathrm{d}z}\right]$. Soit $\sigma \in \Gal(\K/\Q)$. Si $L=\sum\limits_{k=0}^{\mu} a_k(z) \left(\gf{\mathrm{d}}{\mathrm{d}z}\right)^k$, $a_k(z) \in \K[z]$, on définit $L^{\sigma}=\sum\limits_{k=0}^{\mu} a_k^{\sigma}(z) \left(\gf{\mathrm{d}}{\mathrm{d}z}\right)^k$ en étendant $\sigma$ à $\Qbar\llbracket z\rrbracket$ par l'action coefficient par coefficient. 

\bigskip

Alors pour tout $i \in \{1, \dots, \mu\}$, $L^{\sigma}(f_i^{\sigma}(z-\sigma(\alpha)))=0$. De plus, comme $a_{\mu}(\alpha) \neq 0$, on a $a_{\mu}^{\sigma}(\sigma(\alpha))=\sigma(a_{\mu}(\alpha)) \neq 0$, de sorte que $\sigma(\alpha)$ est un point ordinaire de $L^{\sigma}$. Ainsi, $f_i^{\sigma}$ est analytique au voisinage de $0$. Ceci valant pour tout $\sigma$, on en déduit que, si $f_i(u)=\sum\limits_{i=0} b_{i,n} u^n$, il existe une constante $C>0$ telle que $\forall n \in \N, \house{b_{i,n}} \leqslant C^{n+1}$, ce qui prouve que les $f_i(z)$ vérifient la condition \textbf{b)} de la définition \ref{defgfonction}\end{dem}

Le théorème suivant est le résultat principal de cette section. 

\begin{Th}[Chudnovsky, 1984] \label{chudnovsky}
Soit $\mathbf{f}={}^t\left(
f_1(z), \dots, f_n(z) \right)\in \Qbar\llbracket z\rrbracket ^n$ vérifiant $\mathbf{f}'=G\mathbf{f}$, avec $G \in \mathcal{M}_n(\Qbar(z))$. Si, pour tout $i \in \{ 1, \dots, n \}$, $f_i(z)$ est une $G$-fonction et $(f_1(z), \dots, f_n(z))$ est une famille libre sur $\Qbar(z)$, alors $G$ vérifie la condition de Galochkin.
\end{Th}

\begin{rqu}
Mentionnons une conséquence importante de ce théorème : soit $f(z) \in \Qbar\llbracket z\rrbracket $ une $G$-fonction et $$L=\left(\gf{\mathrm{d}}{\mathrm{d}z}\right)^n+a_{1}(z) \left(\gf{\mathrm{d}}{\mathrm{d}z}\right)^{n-1}+\dots+a_n(z) \not\equiv 0$$ un opérateur différentiel de $\Qbar(z)\left[\gf{\mathrm{d}}{\mathrm{d}z} \right]$ d'ordre minimal $n$ tel que $L(f(z))=0$. Alors $L$ est un $G$-opérateur.

En effet, introduisons le vecteur $$\mathbf{f}=\begin{pmatrix}
f \\ 
f' \\ 
\vdots \\ 
f^{(n-1)}
\end{pmatrix} \quad \textup{et} \quad A_L :=\begin{pmatrix}
0 & 1 &  & (0) \\ 
 & \ddots & \ddots &  \\ 
(0) &  & 0 & 1 \\ 
-a_n & \dots &  & -a_1
\end{pmatrix}$$ la matrice compagnon associée à $L$, de sorte que $\mathbf{f}'=A_L \mathbf{f}$. La condition de minimalité de $n$ impose que $(f, \dots, f^{(n-1)})$ est libre sur $\Qbar(z)$. Le théorème des Chudnovsky assure donc la condition de Galochkin est vérifiée pour les $T^m \gf{A_{L,m}}{m!}$ si $T(z)$ est un dénominateur commun des $a_i(z)$. En d'autres termes, $L$ est un $G$-opérateur.
\end{rqu}

\bigskip

Le théorème des Chudnovsky permet d'obtenir des résultats d'irrationalité sur les valeurs de $G$-fonctions, comme le théorème \ref{thmgalochkin}, dû à Galochkin, cité en introduction.

\medskip

Un autre exemple d'application du théorème d'André-Chudnovsky-Katz (théorème \ref{katzchudandre} en introduction), qui découle, entre autres, du théorème des Chudnovsky, est la description des $G$-opérateurs d'ordre 1. 
\begin{prop} \label{prop:Gopordre1}
Les $G$-opérateurs d'ordre $1$ sont les $$L=\gf{\mathrm{d}}{\mathrm{d}z}-\sum\limits_{j=1}^m \gf{r_j}{z-a_j}, \quad r_j \in \Q,\;\; a_j \in \Qbar,\;\; m \in \N.$$ Les solutions de $L(y(z))=0$ sont des $G$-fonctions algébriques sur $\Qbar(z)$.
\end{prop}

\begin{dem}
\begin{itemizeth}
\item Soit une équation différentielle \begin{equation}\label{EDordre1ex2} y'(z)=p(z)y(z), p(z) \in \Qbar(z). \end{equation} vérifiant la condition de Galochkin. Selon le théorème \ref{katzchudandre}, cette équation est fuchsienne, donc $p(z)$ n'a que des pôles simples ; de plus, on peut écrire 
$$ p(z)=\sum_{j=1}^m \gf{r_j}{z-a_j}, $$ où $a_j \in \Qbar$ et $r_j$ est l'exposant de \eqref{EDordre1ex2} en $a_j$. Toujours selon le théorème \ref{katzchudandre}, on sait que $\forall j \in \{1, \dots, m \}, r_j \in \Q$.

\item Soit $L=\gf{\mathrm{d}}{\mathrm{d}z}-p(z)$ de la forme de l'énoncé. Selon \cite[ pp. 145--148]{Hille}, les solutions de l'équation $L(y(z))=0$ sont les $$y(z)=C_0\exp\left(\int_{[z_0,z]} p(s) \mathrm{d}s\right)=\prod\limits_{j=1}^m (z-a_j)^{r_j},$$ où $C_0 \in \Qbar^*$ si bien qu'il existe $q \in \N^*$ tel que $y(z)^q \in \Qbar(z)$  : $y$ est ainsi algébrique sur $\Qbar(z)$. 
\end{itemizeth}
Donc selon la proposition \ref{exemplesgfonctions}, $y$ est une $G$-fonction, de sorte que $L$ est un $G$-opérateur par le théorème \ref{chudnovsky}.\end{dem}

\begin{rqu}
En revanche, la méthode proposée dans la démonstration ci-dessus ne semble pas permettre de déterminer les $G$-opérateurs d'ordre 2.
\end{rqu}
\subsubsection{Démonstration du théorème \ref{chudnovsky}}

Le résultat suivant, énoncé par Siegel dans \cite{Siegelarticle} (mais prouvé antérieurement par Thue dans \cite{Thue}), est un outil classique en approximation diophantienne qui nous servira au cours de la démonstration. Une preuve peut être trouvée dans \cite[p. 37]{Siegel}.

\begin{lem}[Lemme de Siegel] \label{siegellemme}
Soit $\K$ un corps de nombres. Considérons un système de $m$ équations linéaires

\begin{equation}\label{systemesiegel}
 \sum\limits_{j=1}^{n} a_{ij} x_j=0, \forall 1 \leqslant i \leqslant m,
\end{equation} où $\forall i,j, a_{ij} \in \Oal_\K$. On note $A=\max\limits_{i,j} \house{a_{ij}}$. Alors si $n>m$, \eqref{systemesiegel} a une solution non nulle $(x_j)_{1 \leqslant j \leqslant n} \in \Oal_\K^n$ vérifiant$$ \max\limits_{1 \leqslant j \leqslant n} \house{x_j} \leqslant c_1 (c_1 n A)^{\frac{m}{n-m}},$$ où $c_1 >0$ est une constante dépendant uniquement de $\K$. 
\end{lem}

Le lemme \ref{siegellemme} servira de manière essentielle dans l'étape 7 de la démonstration pour montrer l'existence de certains approximants de Padé.

Passons à présent à la preuve du théorème \ref{chudnovsky}. Nous compléterons l'esquisse donnée par Beukers dans \cite[pp. 21--22]{Beukers}. 
Remarquons tout d'abord que si $\K$ est un corps de nombres contenant les coefficients des coefficients de $G(z)$ et les $f_i(0)$, $1 \leqslant i \leqslant \mu$, alors $G \in \mathcal{M}_{n}(\K(z))$ et $\mathbf{f} \in \K\llbracket z\rrbracket^\mu$. En effet, cela découle de l'équation $\mathbf{f}'=G \mathbf{f}$ en écrivant $G$ comme un élément de $\mathcal{M}_{n}\big(\K((z))\big)$ et identifiant les coefficients du développement en série de Laurent de part et d'autre.

\bigskip

\textit{Notations et hypothèses} :

On a $T(z) \in \K[z]$, mais quitte à multiplier par un entier adapté, on peut supposer que $T(z) \in \Oal_\K[z]$ et $T(z) G(z) \in \mathcal{M}_{n}(\Oal_\K [z])$.

On note $D=\gf{\mathrm{d}}{\mathrm{d}z}$ la dérivation usuelle sur $\K\llbracket z\rrbracket$.

Si $A \in \K\llbracket z\rrbracket$ et $\ell \in \N$, on note $A = O\left(z^{\ell} \right)$ s'il existe $B \in \K\llbracket z\rrbracket$ tel que $A=z^{\ell} B$.

On note $\delta=[\K:\Q]$ le degré du corps de nombres $\K$.

\medskip

\textbf{Étape 1} : Soient $N, M \in \N$. On introduit des approximants de Padé $(Q, \mathbf{P})$ de type II de paramètres $(N,M)$ associés à $\mathbf{f}$ dont on laisse les paramètres libres pour l'instant, c'est-à-dire des polynômes $Q, \mathbf{P}_1, \dots, \mathbf{P}_{n} \in \K[z]$ tels que $\deg(Q) \leqslant N$, $\max\limits_{1 \leqslant i \leqslant n} \deg(P_i) \leqslant N$ et si $\mathbf{P}=(P_1, \dots, P_n)$, alors $$ Q \mathbf{f} - \mathbf{P}=O\left(z^{N+M}\right).$$ On ne discutera des conditions d'existence de tels approximants de Padé que dans l'étape~7.

On a pour tout $m < N+M$, $\gf{T^m}{m!} (D-G)^m \mathbf{P} \in \K[z]$, ce qui est immédiat en utilisant la formule de Leibniz. De plus, on va montrer par récurrence sur $m$ que \begin{equation}\label{eq:derivpade}
\forall m \in \N^*, \quad \gf{T^m}{m!} Q^{(m)} \mathbf{f} - \gf{T^m}{m!} (D-G)^m \mathbf{P}=O\left(z^{N+M-m}\right).
\end{equation}

Pour $m=1$, c'est la définition. 

Supposons la relation vraie au rang $m$. Alors en dérivant \eqref{eq:derivpade} et en multipliant par $T$, on~a 
\begin{multline*}
(mT' T^m Q^{(m)}+T^{m+1} Q^{(m+1}) \mathbf{f}+T^{m+1} Q^{(m)} \mathbf{f}' \\ - mT' T^m (D-G)^m \mathbf{P} - T^{m+1} D(D-G)^m \mathbf{P} =O\left(z^{N+M-(m+1)}\right).
\end{multline*}
Or, $\mathbf{f}'=G\mathbf{f}$ , donc en multipliant \eqref{eq:derivpade} par la matrice polynomiale $TG$ et en retranchant à l'équation précédente, on obtient
\begin{multline*}
(mT' T^m Q^{(m)}+T^{m+1} Q^{(m+1}) \mathbf{f}- mT' T^m (D-G)^m \mathbf{P}\\ - T^{m+1} (D-G)^{m+1} \mathbf{P}=O\left(z^{N+M-(m+1)}\right). 
\end{multline*}
Finalement, comme $T^m Q^{(m)}\mathbf{f}- T^m (D-G)^m \mathbf{P}=O\left(z^{N+M-m}\right)$, on a le résultat voulu.

\medskip

\textbf{Étape 2} : Montrons que
\begin{equation} \label{eq:similileibniz}
\forall s \in \N^*, \quad \gf{G_s}{s!} \mathbf{P}=\sum_{j=0}^s \gf{(-1)^j}{(s-j)! j!} D^{s-j}(D-G)^j \mathbf{P}.
\end{equation}
On procède par récurrence :
\begin{itemize}[label=\textbullet]
\item pour $s=1$, $G_s=G$ et $G \mathbf{P}=D(\mathbf{P})-(D-G)\mathbf{P}$.

\item Soit $s \in \N^*$, supposons la formule \eqref{eq:similileibniz} vraie pour $s$. 

Alors $G_{s+1}=G_s G+G'_s$, donc $\gf{G_{s+1}}{(s+1)!}=\gf{1}{s+1} \left( \gf{G_s}{s!} G+\gf{G'_s}{s!} \right)$. Donc en appliquant l'hypothèse de récurrence au vecteur $G \mathbf{P}$, on a $$ \gf{G_{s+1}}{(s+1)!} \mathbf{P}=\gf{1}{s+1} \left(\sum_{j=0}^s \gf{(-1)^j}{(s-j)! j!} D^{s-j}(D-G)^j (G\mathbf{P})+\gf{G'_s}{s!} \mathbf{P} \right).$$ Or, $$\gf{G'_s}{s!} \mathbf{P}=\left(\gf{G_s}{s!} \mathbf{P}\right)' - \gf{G_s}{s!} \mathbf{P}'=\sum_{j=0}^s \gf{(-1)^j}{(s-j)! j!} D^{s+1-j}(D-G)^j \mathbf{P}-\sum_{j=0}^s \gf{(-1)^j}{(s-j)! j!} D^{s-j}(D-G)^j D \mathbf{P}.$$ Donc, en remarquant que $\forall j \in \{ 0, \dots, s \}, (D-G)^j D=(D-G)^{j+1}+(D-G)^j G $, on a 
\begin{align*}
\gf{G_{s+1}}{(s+1)!} \mathbf{P} &= \gf{1}{s+1} \left(\sum_{j=0}^s \gf{(-1)^j}{(s-j)! j!} D^{s+1-j} (D-G)^j \mathbf{P} - \sum_{j=0}^s \gf{(-1)^j}{(s-j)! j!} D^{s-j}(D-G)^{j+1} \mathbf{P} \right) \\
&= \gf{1}{s+1} \left(\sum_{j=0}^s \gf{(-1)^j}{(s-j)! j!} D^{s+1-j} (D-G)^j \mathbf{P} + \sum_{k=1}^{s+1} \gf{(-1)^k}{(s-k+1)! (k-1)!} D^{s+1-k}(D-G)^{k} \mathbf{P} \right) \\
&= \gf{1}{s+1} \left(\sum_{j=1}^s \gf{(-1)^j (s+1-j)+j}{(s+1-j)! j!} D^{s+1-j} (D-G)^j \mathbf{P} + \right. \\ 
& \qquad \left. \sum_{k=1}^{s+1} \gf{(-1)^k}{(s-k+1)! (k-1)!} D^{s+1-k}(D-G)^{k} \mathbf{P} + \gf{D^{s+1}}{s!} \mathbf{P}+\gf{(-1)^{s+1}}{s!}(D-G)^{s+1} \mathbf{P} \right) \\
&= \sum\limits_{j=0}^{s+1} \gf{(-1)^j}{(s+1-j)! j!} D^{s+1-j} (D-G)^j \mathbf{P}. 
\end{align*}
Cela conclut la récurrence.
\end{itemize}

\bigskip

\textbf{Étape 3} : Utilisation du lemme de Shidlovskii.

On note pour $h \in \N$, $\mathbf{P}_h=\gf{1}{h!} (D-G)^h \mathbf{P}$ et $R_{(h)} \in \mathcal{M}_{n}(\K(z))$ la matrice dont la $j$\up{ème} colonne est $\binom{h+j-1}{j-1} \mathbf{P}_{h+j-1}$. Alors la formule \eqref{eq:similileibniz} implique immédiatement que

$$ \gf{G_s}{s!} R_{(0)}=\sum\limits_{j=0}^s \gf{(-1)^j}{(s-j)!} D^{s-j} R_{(j)}.$$ Selon le lemme de Shidlovskii (théorème \ref{Shidlovskii}) qui sera démontré dans la partie \ref{subsec:shidlovskii},, la matrice $R_{(0)}$ est inversible pourvu que $M$ soit assez grand ce qui sera réalisé quand on spécifiera $M$ et $N$ dans l'étape 7. Donc \begin{equation} \label{eq:conclusionetape3}
T^s \gf{G_s}{s!}=\sum\limits_{j=0}^s (-1)^j T^{s+n-1} \gf{D^{s-j} R_{(j)}}{(s-j)!} (T^{n-1} R_{(0)})^{-1}.
\end{equation}

\medskip

\textbf{Étape 4} :
Soit $d$ le dénominateur commun des coefficients d'ordres inférieurs à $N+M$ du développement en série entière de $\mathbf{f}$. On suppose trouvés des approximants de Padé $Q, \mathbf{P}_1, \dots, \mathbf{P}_{n} \in \K[X]$ de $d\mathbf{f}$, c'est-à-dire des polynômes tels que $\deg(Q) \leqslant N$, $\max\limits_{1 \leqslant i \leqslant n} \deg(P_i) \leqslant N$ et $$ Q(d\mathbf{f}) - \mathbf{P}=O\left(z^{N+M}\right),$$ avec $\mathbf{P}=(P_1, \dots, P_{n})$. On fait l'hypothèse supplémentaire que $Q$ est un polynôme à coefficients entiers algébriques. On peut alors appliquer les résultats des trois étapes précédentes, dont on conservera les notations, à $Q$ et $\mathbf{P}$ puisque $d\mathbf{f}$ est encore solution du système différentiel $y'=Gy$.

Selon \eqref{eq:derivpade}, on a
\begin{equation} \label{eq:derivpadedenom}
 \forall m \leqslant N+M, \quad  \gf{T^m}{m!} Q^{(m)} (d\mathbf{f}) - T^m \mathbf{P}_m=O\left(z^{N+M-m}\right).
\end{equation}
On remarque que $\gf{Q^{(m)}}{m!} \in \Oal_\K[z]$ puisque $Q$ est à coefficients dans $\Oal_\K$. On en déduit que si $N+M-m > \max\limits_{1 \leqslant i \leqslant n} \deg(\mathbf{P}_{i,m})$, les coefficients de $T^m \mathbf{P}_m$ sont des éléments de $\Oal_\K[z]$.

\medskip

Montrons par récurrence sur $m$ que \begin{equation} \label{eq:bornedegrePm}
\max\limits_{1 \leqslant i \leqslant n} \deg(T^m \mathbf{P}_{i,m}) \leqslant N+tm.
\end{equation}

\begin{itemize}[label=\textbullet]
\item Pour $m=0$, il s'agit simplement du fait que les composantes de $\mathbf{P}$ sont de degrés inférieurs à $N$.

\item Pour $m=1$, $T \mathbf{P}'-T G \mathbf{P}$ a des composantes de degrés inférieurs à $N+t$, car les coefficients de $T(z)G$ sont de degrés bornés par $t$.

\item Soit $m \in \N$, supposons le résultat vrai au rang $m$. Alors\begin{align*}
(D-G)(T^m (D-G)^m \mathbf{P}) &= m T' T^{m-1} (D-G)^m \mathbf{P}+T^m D(D-G)^m \mathbf{P}-T^m G(D-G)^m \mathbf{P} \\
&= m T' T^{m-1} (D-G)^m \mathbf{P}+T^m (D-G)^{m+1} \mathbf{P}.
\end{align*} Donc $$ T^{m+1} (D-G)^{m+1} \mathbf{P}=(D-G)(T^m (D-G)^m \mathbf{P})- m T' T^m (D-G)^m \mathbf{P}. $$
Or, en utilisant à la fois l'hypothèse de récurrence et le cas $m=1$, on voit que $T (D-G)(T^m (D-G)^m \mathbf{P})$ a ses composantes de degrés bornés par $N+mt+t=N+(m+1)t$ ; par ailleurs, l'hypothèse de récurrence nous assure que $m T' T^m (D-G)^m \mathbf{P}$ a des composantes de degrés inférieurs à $t+N+mt=N+(m+1)t$. On en déduit le résultat souhaité~\eqref{eq:bornedegrePm}.
\end{itemize}

On déduit, à condition que $N+M-m \geqslant N+tm$, c'est-à-dire 
\begin{equation} \label{eq:conditionentieral}
m \leqslant \gf{M}{t+1},
\end{equation}
que $T^m \mathbf{P}_m \in \Oal_\K[z]^n$. Ceci implique immédiatement que si $j$ est un entier naturel tel que $j+n-1 \leqslant \gf{M}{t+1}$, alors $T^{j+n-1} R_{(j)}$ est une matrice à coefficients dans $\Oal_\K [z]$. En particulier, $T^{n-1} R_{(0)} \in \mathcal{M}_{n}(\Oal_\K[z])$.

\medskip

Montrons à présent que \begin{equation} \label{eq:derivrjentier}
\forall s \in \N, \;\; s+n-1 \leqslant \gf{M}{t+1}, \;\; \forall j \in \{ 0, \dots, s \}, \quad \gf{T^{s+n-1}}{(s-j)!} D^{s-j} R_{(j)} \in \mathcal{M}_{n}(\Oal_\K [z]).
\end{equation}
Rappelons la formule de Leibniz généralisée : si $\ell \in \N^*$ et $f_1, \dots, f_{\ell}$ sont $l$ fonctions dérivables $k$ fois, alors $$ (f_1 \dots f_{\ell})^{(k)}=\sum_{i_1+\dots+i_{\ell}=k} \binom{k}{i_1, \dots, i_{\ell}} \prod_{1 \leqslant t \leqslant k} f_t^{(i_t)}.$$ En particulier, considérons un élément quelconque $W \in \Oal_\K[z]$. Alors $$\gf{(W^\ell)^{(k)}}{k!}=\sum\limits_{i_1+\dots+i_{\ell}=k} \prod\limits_{1 \leqslant t \leqslant \ell} \gf{W^{(i_t)}}{i_t!}.$$ Si $k <\ell$, pour tout $(i_1, \dots, i_{\ell})$ intervenant dans la somme, chaque terme $\prod\limits_{1 \leqslant t \leqslant \ell} \gf{W^{(i_t)}}{i_t!}$ contient au moins $\ell-k$ indices de dérivation nuls. Ainsi, comme pour tout entier $s$,$\gf{W^{(s)}}{s!} \in \Oal_\K[z]$, on a $\gf{(W^\ell)^{(k)}}{k!} \in W^{\ell-k} \Oal_\K[z]$.

\medskip

Déduisons de cela par récurrence le résultat \eqref{eq:derivrjentier}. Pour $s=0$, c'est évident.

Soit $s \in \N^*$, supposons \eqref{eq:derivrjentier} vrai pour $s' \in \{ 0, \dots, s-1 \}$. Par la formule de Leibniz, on a, pour $j \in \{0, \dots, s \}$, \begin{align*}
\gf{D^{s-j}\big(T^{s+n-1} R_{(j)}\big)}{(s-j)!} &= \sum_{k=0}^{s-j} \binom{s-j}{k} \times \gf{1}{(s-j)!} \big(T^{s+n-1}\big)^{(k)} D^{s-j-k} R_{(j)} \\
&= T^{s+n-1} \gf{D^{s-j} R_{(j)}}{(s-j)!}+ \sum_{k=1}^{s-j} \gf{\big(T^{s+n-1}\big)^{(k)}}{k!} \gf{D^{s-j-k} R_{(j)}}{(s-j-k)!} \\
&= T^{s+n-1} \gf{D^{s-j} R_{(j)}}{(s-j)!}+ \sum_{k=1}^{s-j} U_k T^{s+n-1-k} \gf{D^{s-k-j} R_{(j)}}{(s-k-j)!},
\end{align*} avec $U_k \in \Oal_\K[z]$, en utilisant la remarque précédente.

Or, d'une part, $$\gf{D^{s-j}\big(T^{s+n-1} R_{(j)}\big)}{(s-j)!} \in \mathcal{M}_{n}(\Oal_\K[z]),$$ puisque $T^{s+n-1}=T^{s-j} T^{j+n-1} R_{(j)} \in \mathcal{M}_{n}(\Oal_\K[z])$, et d'autre part, par hypothèse de récurrence, $$\forall k \in \{ 1, \dots, s-j \}, \quad T^{s-k+n-1} \gf{D^{s-k-j} R_{(j)}}{(s-k-j)!} \in \mathcal{M}_{n}(\Oal_\K[z]).$$ Par conséquent, $$T^{s+n-1} \gf{D^{s-j} R_{(j)}}{(s-j)!} \in \mathcal{M}_{n}(\Oal_\K[z]),$$ ce qu'il fallait démontrer.

\bigskip

\textbf{Étape 5} : Lien entre $q_s$ et taille des coefficients de $\det(T^{n-1} R_{(0)})$.

Selon \eqref{eq:conclusionetape3}, on a pour tout $s \in \N$, \begin{equation} \label{eq:formuleetape5} T^s \gf{G_s}{s!}=\gf{1}{V} \sum\limits_{j=0}^s (-1)^j T^{s+n-1} \gf{D^{s-j} R_{(j)}}{(s-j)!}  \left(\mathrm{com}\left(T^{n-1} R_{(0)}\right)\right)^{^T},\end{equation} où $V=\det(T^{n-1} R_{(0)}) \in \Oal_\K[z]$, et selon \eqref{eq:derivrjentier}, tous les termes de la somme sont à coefficients entiers algébriques à condition que $s+n-1 \leqslant \gf{M}{t+1}$.

Prenons pour $s \in \N$, $q_s$ le dénominateur des coefficients des coefficients des matrices $TG, T^2 \gf{G_2}{2!}, \dots, T^s \gf{G_s}{s!}$, comme dans la définition \ref{galochkin}. Pour estimer $q_s$ sous la condition $s+n-1 \leqslant \gf{M}{t+1}$, il suffit donc d'obtenir une estimation de la maison des coefficients de $V_s$. C'est ce que permet de faire ce lemme :

\begin{lem} \label{lem:etape5}
Soient $U \in \Oal_\K[z]$, $V \in \Oal_\K[z]$, $W \in \K[z]$ tels que $U=VW$. Notons $V=\sum\limits_{i=0}^\ell v_i z^i$, alors pour tout $k$ tel que $v_k \neq 0$, on a $N_{\K/\Q}(v_k) W \in \Oal_\K[z]$.
\end{lem}

\begin{dem}
Introduisons la valuation de Gauss associée à un premier $\mathfrak{p}$ de $\Oal_\K$ : $$v_{\mathfrak{p}}\left(\sum\limits_{i=0}^{q} a_i z^i \right) :=\min\limits_{0 \leqslant i \leqslant q}(v_{\mathfrak{p}}(a_i)),$$ où $v_{\mathfrak{p}}$ est la valuation $\mathfrak{p}$-adique associée à l'idéal premier $\mathfrak{p}$ de $\Oal_\K$. En utilisant les propriétés de valuation, on a $v_{\mathfrak{p}}(U)=v_{\mathfrak{p}}(V)+v_{\mathfrak{p}}(W)$, donc $v_{\mathfrak{p}}(W) \geqslant - v_{\mathfrak{p}}(V)$ car, comme $U$ est à coefficients entiers algébriques, $v_{\mathfrak{p}}(U) \geqslant 0$.

Notons $S$ l'ensemble fini des premiers divisant tous les coefficients de $V$. Alors si $W=\sum\limits_{i=0}^{d} w_i z^i$, pour tout $i \in \{ 0, \dots, d \}$ et $\mathfrak{p} \in S$, on a $v_{\mathfrak{p}}(w_i) + v_{\mathfrak{p}}(V) \geqslant v_{\mathfrak{p}}(V)+v_{\mathfrak{p}}(W) \geqslant 0$, donc $ \prod\limits_{\mathfrak{p} \in S} \mathfrak{p}^{v_{\mathfrak{p}}(V)} (w_i) \subset \Oal_\K$. En particulier, si $v_{k} \neq 0$, comme $v_{k} \in \prod\limits_{\mathfrak{p} \in S} \mathfrak{p}^{v_{\mathfrak{p}}(V)}=\mathrm{pgcd}((v_0), \dots, (v_{\ell}))$, on a $$\forall i \in \{ 0, \dots, d \}, \quad (v_k) (w_i) \subset \Oal_\K.$$ D'où comme $N_{\K/\Q}(v_k) \in (v_k) \cap \Z$, on a $N_{\K/\Q}(v_k) W \in \Oal_\K[z]$.\end{dem}

Pour $W=\sum\limits_{i=0}^\ell w_i z^i \in \K[z]$, on définit $\sigma(W)=\max\limits_{0 \leqslant i \leqslant \ell} \house{w_i}$, la \textit{maison} de $W$. Selon le lemme \ref{lem:etape5} appliqué à la formule \eqref{eq:formuleetape5} avec $$U=\sum\limits_{j=0}^s (-1)^j T^{s+n-1} \gf{D^{s-j} R_{(j)}}{(s-j)!}  \left(\mathrm{com}\left(T^{n-1} R_{(0)}\right)\right)^{^T},$$ $V=\det(T^{n-1} R_{(0)})$ et $W=T^s \gf{G_s}{s!}$, on a donc ici \begin{equation} \label{eq:conclusionetape5}
\forall s \in \N \;\; \text{tel que} \;\; s+n-1 \leqslant \gf{M}{t+1}, \;  \; q_s \leqslant \sigma(V)^{\delta}.
\end{equation}

\bigskip

\textbf{Étape 6} : Majoration de la taille des coefficients de $\det(T^{n-1} R_{(0)}) \in \Oal_\K [z]$ en fonction de la maison de $Q$. 

Par commodité, on s'intéresse à $$\widetilde{V}=\det(\mathbf{P}, T \mathbf{P}_1, \dots, T^{n-1} \mathbf{P}_{n-1})=T^{\frac{n(n-1)}{2}} \det(R_{(0)})= T^{-\frac{n(n-1) }{2}} V.$$ Le lemme suivant nous assure que ce changement n'introduit qu'une constante multiplicative dépendant seulement de $G$ dans la majoration recherchée.

\begin{lem} \label{lem:etape6}
Soient $A, B \in \K[z]$ et $C=AB$, alors $\sigma(C) \leqslant (\deg(A)+\deg(B)+1) \sigma(A) \sigma(B)$.
\end{lem}

\begin{dem}
On écrit $A=\sum\limits_{i=0}^{p} a_i z^i$ et $B=\sum\limits_{i=0}^q b_i z^i$, de sorte que $C=\sum\limits_{j=0}^{p+q} c_j z^j$, avec, pour tout $j \in \{ 0, \dots, p+q \}$, $c_j=\sum\limits_{i=0}^{j} a_i b_{j-i}$.

Soit $\tau : \K \hookrightarrow \Qbar$ un plongement. Alors $$ |\tau(c_j)| \leqslant \sum_{i=0}^{j} |\tau(a_i)| |\tau(b_{j-i})| \leqslant (j+1) \sigma(A) \sigma(B) \leqslant (p+q+1) \sigma(A) \sigma(B),$$ si bien qu'en prenant le maximum sur $j$ et sur $\tau$, on obtient l'inégalité voulue.\end{dem}

Soit $m \in \{0, \dots, n-1\}$. Si $Q=\sum\limits_{i=0}^N q_i z^i$, alors $$\gf{Q^{(m)}}{m!}=\sum\limits_{i=0}^{N-m} \gf{(i+1) \dots (i+m)}{m!} q_{m+i} z^m=\sum\limits_{i=0}^{N-m} \binom{m+i}{i} q_{m+i} z^m.$$
De la majoration $\binom{m+i}{i} \leqslant 2^{m+i} \leqslant 2^{N}$, il s'ensuit que $\sigma\left(\gf{Q^{(m)}}{m!}\right) \leqslant 2^N \sigma(Q)$. Selon le lemme \ref{lem:etape6} appliqué à $A=T^m$ et $B=Q^{(m)}/m!$,

\begin{multline*}\sigma\left(T^m \gf{Q^{(m)}}{m!} \right) \leqslant 2^N \sigma(Q) \sigma(T^m)(mt+N-m+1) \\ \leqslant 2^N \sigma(Q) \sigma(T^m)((n-1)(t-1)+N+1) \leqslant c_1^N \sigma(Q) \sigma(T^m), 
\end{multline*} avec $c_1$ constante dépendant seulement de $G$ pour $N$ suffisamment grand. 

En appliquant $m$ fois le le lemme \ref{lem:etape6}, on obtient $$\sigma(T^m) \leqslant c_2^N \sigma(T) \sigma(T^{m-1}) \leqslant \dots \leqslant c_3^N \sigma(T)^m,$$ où $c_2$ et $c_3$ sont des constantes. Dans ce qui suit, les $c_i$ désigneront des constantes.

Ainsi, pour $N$ suffisamment grand, $$\sigma\left(T^m \gf{Q^{(m)}}{m!} \right) \leqslant c_4^N \sigma(Q) \sigma(T)^m \leqslant c_5^N \sigma(Q).$$

Soit $\theta_{N+M}$ le maximum des maisons des $N+M$ premiers coefficients du développement en série entière des $f_i$, et $d_{N+M}$ leur dénominateur commun. En répétant le raisonnement de la preuve du lemme \ref{lem:etape6}, on voit que la  maison de la partie polynomiale tronquée à l'ordre $N+tm$ de $T^m \gf{Q^{(m)}}{m!} (d_{N+M} f_i)$ est majorée par $$c_5^N \sigma(Q) (d_{N+M} \theta_{N+M})(N+tm+1) \leqslant c_6^N \sigma(Q) (d_{N+M} \theta_{N+M}).$$ Or, si $N+M-(n-1) \geqslant N+t(n-1)$, selon \eqref{eq:derivpadedenom}, cette partie polynomiale est $T^m \mathbf{P}_m$, donc, avec une extension de la notation $\sigma$ aux vecteurs colonnes, $$\forall m \in \{ 0, \dots, n-1 \}, \quad \sigma(T^m \mathbf{P}_m) \leqslant c_7^N \sigma(Q) d_{N+M} \theta_{N+M}.$$ On a $\widetilde{V}=\sum\limits_{\tau \in \mathfrak{S}_{n}} \varepsilon(\tau) \prod\limits_{j=0}^{n-1} T^j \mathbf{P}_{\tau(j),j}$. Pour $\tau \in \mathfrak{S}_{n}$, en appliquant le lemme \ref{lem:etape6} à $A=\mathbf{P}_{\tau(0),0}$ et $B=\prod\limits_{j=1}^{n-1} T^j \mathbf{P}_{\tau(j),j}$, on a 

$$\sigma\left(\prod\limits_{j=0}^{n-1} T^j \mathbf{P}_{\tau(j),j}\right) \leqslant \sigma(\mathbf{P}_{\tau(0),0}) \sigma\left(\prod\limits_{j=1}^{n-1} T^j \mathbf{P}_{\tau(j),j}\right) (n(N+t(n-1)+1)$$ en utilisant \eqref{eq:bornedegrePm}. En itérant le procédé, on obtient une constante $c_8$ telle que

$$\sigma\left(\prod\limits_{j=0}^{n-1} T^j \mathbf{P}_{\tau(j),j}\right) \leqslant c_7^{n N} \sigma(Q)^n d_{N+M}^n \theta_{N+M}^n c_8^N \leqslant c_9^N \sigma(Q)^n d_{N+M}^n \theta_{N+M}^n.$$ Donc par inégalité triangulaire $$\sigma(\widetilde{V}) \leqslant c_{10}^N \sigma(Q)^n d_{N+M}^n \theta_{N+M}^n,$$ si bien que, comme $V=T^{n(n-1)/2} \widetilde{V}$, $\sigma(V) \leqslant c_{11}^N \sigma(Q)^n d_{N+M}^n \theta_{N+M}^n$. D'où selon \eqref{eq:conclusionetape5},

\begin{equation} \label{eq:conclusionetape6}
\forall s \in \N, \;\; s+n-1 \leqslant \gf{M}{t+1}, \quad q_s \leqslant \sigma(V)^{\delta} \leqslant c_{12}^N \sigma(Q)^{n \delta} d_{N+M}^{n \delta} \theta_{N+M}^{n \delta}.
\end{equation}

\bigskip

\textbf{Étape 7} : Conclusion à l'aide du lemme de Siegel.

Soit $s \in \N^*$, supposé suffisamment grand. Soit $N=2n(t+1)(s+n-1)$ et $M=\gf{N}{2n} \in \N^*$. Alors l'équation de Padé $Q (d_{N+M} \mathbf{f})-\mathbf{P}=O\left(z^{N+M} \right)$ se traduit par un système linéaire de $\gf{nN}{2n}=\gf{N}{2}$ équations à $N+1$ inconnues (les coefficients de $Q$). Selon le lemme \ref{siegellemme}, il existe une solution $Q \in \Oal_\K[z]$ telle que $$\sigma(Q) \leqslant c_{10} (c_{10} (N+1) \theta_{N+M})^{\frac{N/2}{N+1-N/2}} \leqslant c_{11}^N \theta_{N+M}$$ car $\gf{N}{2} \leqslant N+1-\gf{N}{2}$. Puisque les composantes de $\mathbf{f}$ sont des $G$-fonctions, on peut trouver des constantes $c_{12}$ et $c_{13}$ telles que $\theta_{N+M} \leqslant c_{12}^{N+M}$ et $d_{N+M} \leqslant c_{13}^{N+M}$.

On a, par définition de $M$, $\gf{M}{t+1} \geqslant s+n-1$, donc selon \eqref{eq:conclusionetape6}, $$q_s \leqslant c_9^N (c_{11}^N \theta_{N+M})^{n \delta} d_{N+M}^{\delta n} \theta_{N+M}^{\delta n} \leqslant c_{14}^N$$ de sorte qu'on peut trouver une constante $c_{15} >0$ dépendant seulement de $G$ telle que $\forall s \in \N^*, q_s \leqslant c_{15}^{s}$. Ceci conclut la démonstration du théorème \ref{chudnovsky}

\subsection{Lemme de Shidlovskii pour les approximants de Padé de type II} \label{subsec:shidlovskii}

Dans l'étape 3 de la démonstration du théorème des Chudnovsky, nous avons utilisé le lemme de Shidlovskii pour les approximants de Padé de type II (théorème \ref{Shidlovskii}), et l'objet de cette partie est de le démontrer.

Soit $\K$ un corps de nombres. Soit $\mathbf{f}={}^t (f_1(z), \dots, f_n(z)) \in \K\llbracket z\rrbracket ^n$ vérifiant un système d'équations $\mathbf{y}'=G \mathbf{y}$, $G \in \mathcal{M}_n(\K(z))$. L'hypothèse que les $f_i$ sont des $G$-fonctions n'est ici pas nécessaire. Supposons de plus que $(f_1(z), \dots, f_n(z))$ constitue une famille libre sur $\K(z)$.

On note $T(z) \in \K[z]$ un dénominateur commun des coefficients de $G$ et on définit $t$ le maximum du degré de $T$ et des degrés des coefficients de $T G$. 

Soient $N, M \in \N$. On se donne $Q, P_1, \dots, P_n$ des approximants de Padé de type II\footnote{Un approximant de Padé de type I et de paramètre $N$ de $\mathbf{f}$ est une fonction $\sum\limits_{i=1}^n P_i(z) f_i(z), \deg(P_i) \leqslant N$, d'ordre en $0$ suffisamment grand par rapport à $N$. On dispose pour le type I d'un résultat analogue au théorème \ref{Shidlovskii}, voir \cite[p.100]{Shidlovskii} ou \cite[p. 112]{Baker}} de $\mathbf{f}$ : $\deg(Q) \leqslant N$, $\max\limits_{1 \leqslant i \leqslant n} \deg (P_i) \leqslant N$ et

$$ Q \mathbf{f} - \mathbf{P}=O\left(z^{N+M}\right), $$ où $\mathbf{P}={}^t (P_1, \dots, P_n)$. Comme on l'a vu dans l'étape 1 de la sous-section \ref{subsec:thchudnovsky}, on a alors pour tout $m \in \N$,

$$ \gf{T^m}{m!} Q^{(m)} \mathbf{f} - \gf{T^m}{m!} (D-G)^m \mathbf{P}=O\left(z^{N+M-m}\right),$$ de sorte qu'en posant $$Q_m=\gf{T^m}{m!} Q^{(m)} \quad  \textup{et} \quad \mathbf{P}_m=\gf{T^m}{m!} (D-G)^m \mathbf{P}={}^t (\mathbf{P}_{1,m}, \dots, 
\mathbf{P}_{n,m}),$$ $(Q_m, \mathbf{P}_{1,m}, \dots, \mathbf{P}_{n,m})$ est un système d'approximants de Padé de type II  de $\mathbf{f}$.

On pose par convention $\mathbf{P}_0=\mathbf{P}$. On s'intéresse à la matrice $R_{(0)} \in \mathcal{M}_n(\K[z])$ de $j$-ième colonne $\mathbf{P}_{j-1}$, et à son déterminant $\Delta(z)=\det(R_{(0)}) \in \K[z]$. Le théorème suivant, qui est un point crucial de la démonstration du théorème \ref{chudnovsky}, a été prouvé dans \cite[pp. 42--43]{Chudnovsky}, avec une correction issue de \cite[pp. 115--119]{Andre}.

\begin{Th}[Shidlovskii pour le type II] \label{Shidlovskii}
Le déterminant $\Delta(z)$  n'est pas identiquement nul si $M$ est assez grand.
\end{Th}

Le reste de cette section est consacré à la démonstration de ce théorème.

\subsubsection{Lemmes techniques} 

On utilise dans la preuve du théorème \ref{Shidlovskii} deux lemmes techniques, les lemmes \ref{lemmebaker} et \ref{ordannulCL} ci-dessous. L'objet de cette partie est de les démontrer.

On commence par introduire la notion de degré total d'une fraction rationnelle :
\begin{defi}
Le \emph{degré total} de $F=\gf{P}{Q} \in \K(z)$, où $P,Q \in \K[z]$ sont premiers entre eux est~$\deg_t(F)=\deg(P)+ \deg(Q)$.
\end{defi}

\begin{rqu}
Si $F=\gf{P_1}{Q_1}$, avec $P_1, Q_1 \in \K[z]$ quelconques, alors on voit que $$\deg_t(F) \leqslant \deg(P_1)+ \deg(Q_1).$$ De plus, si $F_1, F_2 \in \K(z)$, $\deg_t (F_1 F_2) \leqslant \deg_t(F_1) + \deg_t(F_2)$.
\end{rqu}

\begin{lem}[\cite{Baker}, pp. 112--113] \label{lemmebaker}
Soient $G \in \mathcal{M}_n(\K(z))$ et $\mathbf{P} \in \K(z)^n$, on définit $\widetilde{\mathbf{P}}_m=(D-G)^m \mathbf{P}$. Notons $\ell$ le rang de la famille $(\mathbf{P}, \widetilde{\mathbf{P}}_1, \dots, \widetilde{\mathbf{P}}_{n-1})$. Alors $\widetilde{F}=(\mathbf{P}, \widetilde{\mathbf{P}}_1, \dots, \widetilde{\mathbf{P}}_{\ell-1}) \in \mathcal{M}_{n,\ell}(\K(z))$ est de rang $\ell$. Quitte à renuméroter les composantes de $\mathbf{P}$, on peut supposer que la matrice $\widetilde{R} \in \mathcal{M}_{\ell}(\K(z))$ formée des $\ell$ premières lignes de $\widetilde{F}$ est inversible. On note $\widetilde{S} \in \mathcal{M}_{n-\ell,\ell}(\K(z))$ la matrice formée des $n-\ell$ dernières lignes de $\widetilde{F}$.

Alors les coefficients de $\widetilde{S}\widetilde{R}^{-1}$ sont des fractions rationnelles de degrés totaux bornés par une constante $c_0$ ne dépendant que de $G$. 
\end{lem}
 
Le lemme suivant sera utilisé dans la preuve du lemme \ref{lemmebaker}.
\begin{lem}[\cite{Shidlovskii}, p. 86] \label{Shidlovskiifracanalytique}
Soit $\mathbb{L}$ une extension de $\K(z)$ et soient $\varphi_1, \dots, \varphi_s$, $\psi_1, \dots, \psi_m$ des éléments de $\mathbb{L}$ tels que les $\psi_i$ ne sont pas tous nuls. Alors il existe $N_0 \in \N$ tel que pour tous $(\alpha_1, \dots, \alpha_s) \in \C^s$ et $(\beta_1, \dots, \beta_m) \in \C^m$  tels que $\sum\limits_{i=1}^m \beta_i \psi_i \neq 0$,

$$ \omega=\left(\sum\limits_{i=1}^s \alpha_i \varphi_i\right)\left(\sum\limits_{i=1}^m \beta_i \psi_i\right)^{-1} \in \C(z) \Rightarrow \deg_t(\omega) \leqslant N_0.$$
\end{lem}

\begin{dem}
Soient $(\alpha_1, \dots, \alpha_s)$ et $(\beta_1, \dots, \beta_m)$ comme dans l'énoncé du théorème. Soit $r$ le rang sur $\C(z)$ de $\left\{\varphi_1, \dots, \varphi_s, \psi_1, \dots, \psi_m \right\}$. On choisit $r$ fonctions extraites de cette famille $(g_1, \dots, g_r)$ formant une famille libre sur $\C(z)$. 

Notons pour tout $i \in \{ 1, \dots, s \}, \varphi_i=\sum\limits_{j=1}^r A_{ij} g_j, A_{ij} \in \C(z)$, et pour  $i \in \{ 1, \dots, m \} $, $\psi_i=\sum\limits_{j=1}^r B_{ij} g_j$, $B_{ij} \in \C(z)$. Notons $N$ le maximum des degrés totaux des $A_{ij}$ et des $B_{ij}$. C'est un nombre qui dépend uniquement des $\varphi_i$ et des $\psi_i$. On a donc 

$$ \omega=\gf{\sum\limits_{i=1}^s \sum\limits_{j=1}^r \alpha_i A_{ij} g_j}{\sum\limits_{i=1}^m \sum\limits_{j=1}^r \beta_i B_{ij} g_j}= \gf{\sum\limits_{j=1}^r C_j g_j}{\sum\limits_{j=1}^r D_j g_j}, $$ où $C_j=\sum\limits_{i=1}^s \alpha_i A_{ij}$, $D_j=\sum\limits_{i=1}^m \beta_i B_{ij}$.

Par hypothèse, $\sum\limits_{j=1}^r D_j g_j$ n'est pas nulle donc comme $(g_1, \dots, g_r)$ est libre sur $\C(z)$, il existe $j_0 \in \{ 1, \dots, r \}$ tel que $D_{j_0} \neq 0$. 

De plus, $\omega \in \C(z)$ et $\sum\limits_{j=1}^r (D_j \omega - C_j) g_j=0$, si bien que $\forall j \in \{ 1, \dots, r \}, D_j \omega - C_j=0$. En prenant $j=j_0$, on obtient $\omega=\gf{C_{j_0}}{D_{j_0}}$. Il est clair que $\deg_t(C_{j_0}) \leqslant 2sN$ et $\deg_t(D_{j_0}) \leqslant 2mN$, donc $\deg_t(\omega) \leqslant 2(s+m)N$ est borné par une constante ne dépendant que des $\varphi_i$ et des $\psi_i$.\end{dem}

\begin{dem}[du lemme \ref{lemmebaker}]

Commençons par justifier que $\widetilde{F}$ est de rang $\ell$ en montrant que $(\mathbf{P}, \mathbf{P}_1, \dots, \mathbf{P}_{\ell-1})$ est libre sur $\K(z)$. En effet, supposons le contraire et fixons $k \in \{ 1, \dots, \ell-1 \}$ tel que $\mathbf{P}_k \in \Vect_{\K(z)}(\mathbf{P}, \mathbf{P}_1, \dots, \mathbf{P}_{k-1})$. Alors on montre par récurrence, en utilisant la formule $\widetilde{\mathbf{P}}_m=(D-G)^m \mathbf{P}$ que $\widetilde{\mathbf{P}}_{k+s} \in \Vect_{\K(z)}(\mathbf{P}, \widetilde{\mathbf{P}}_1, \dots, \widetilde{\mathbf{P}}_{k-1})$ pour tout $s \in \N$, de sorte que le rang de $(\widetilde{\mathbf{P}}, \widetilde{\mathbf{P}}_1, \dots, \widetilde{\mathbf{P}}_{n-1})$ est strictement inférieur à $\ell$, ce qui est absurde.

\medskip

On considère à présent le système, dont l'inconnue est un vecteur ligne $y$, \begin{equation} \label{systbaker} y'=y(-^{t} G).\end{equation}
On peut trouver une matrice $Q \in \mathcal{M}_{\ell}(\K(z))$ telle que si $y=(y_1, \dots, y_n)$ est solution de \eqref{systbaker}, alors le vecteur ligne $Y=yF$ vérifie $Y'=YQ$. En effet, notons $G=(g_{ij})_{1 \leqslant i,j \leqslant n}$ et $Y=(Y_1, \dots, Y_{\ell})$. Alors on a, pour $1 \leqslant j \leqslant \ell-1$, $$ Y'_{j} = \left(\sum\limits_{i=1}^n y'_i \widetilde{\mathbf{P}}_{i,j-1} \right)'=\sum\limits_{i=1}^n y'_i \widetilde{\mathbf{P}}_{i,j-1} + y_i \widetilde{\mathbf{P}}'_{i,j-1} = \sum\limits_{i=1}^n \sum\limits_{k=1}^n (-g_{ki}) y_k \widetilde{\mathbf{P}}_{i,j-1}+y_i \widetilde{\mathbf{P}}'_{i,j-1}.$$ Or, $\widetilde{\mathbf{P}}_{j}=\widetilde{\mathbf{P}}'_{j-1}-G \widetilde{\mathbf{P}}_{j-1}$, donc $$\forall i \in \{ 1, \dots, n \}, \widetilde{\mathbf{P}}'_{i,j-1}=\widetilde{\mathbf{P}}_{i,j}+\sum\limits_{k=1}^n g_{ik} \widetilde{\mathbf{P}}_{k,j-1},$$ de sorte que $$Y'_{j}=-\sum\limits_{i=1}^n \sum\limits_{k=1}^n g_{ki} y_k \widetilde{\mathbf{P}}_{i,j-1}+ \sum\limits_{i=1}^n \sum\limits_{k=1}^n g_{ik} \widetilde{\mathbf{P}}_{k,j-1} y_i+\sum\limits_{i=1}^n y_i \widetilde{\mathbf{P}}_{i,j}=\sum\limits_{i=1}^n y_i \widetilde{\mathbf{P}}_{i,j}=Y_{j+1}.$$ Le même calcul nous donne $Y'_{\ell}=\sum\limits_{i=1}^n y_i \widetilde{\mathbf{P}}_{i,\ell}$, donc, comme $\widetilde{\mathbf{P}}_{\ell} \in \Vect_{\K(z)}(\mathbf{P}, \widetilde{\mathbf{P}}_1, \dots, \widetilde{\mathbf{P}}_{\ell-1})$, $Y'_{\ell}$ est une combinaison linéaire à coefficients dans $\K(z)$ (ne dépendant pas de $y$) des $Y_j, 1 \leqslant~j \leqslant~\ell$.

\medskip

Soit $\mathbb{L}$ une extension de Picard-Vessiot de $\K(z)$ telle que le système \eqref{systbaker} a une base de solutions $(w_1, \dots, w_n) \in \mathbb{L}^n$ sur $\K(z)$. Par les propriétés du wronskien, c'est également une base sur $\K$ de ce système (voir \cite[pp. 7--9]{Singer}). On introduit la matrice wronskienne associée $$W=\begin{pmatrix}
w_1 \\ 
 \vdots \\ 
w_n
\end{pmatrix} \in \mathcal{M}_n(\mathbb{L}).$$
Donc pour tout $i \in \{ 1, \dots, n \}$, la ligne $i$ $Y_i=w_i F$ de $WF$ vérifie $Y'_i=Y_i Q$. Or, le système $Y'=YQ$ n'a que $\ell$ solutions linéairement indépendantes sur $\K$. Par suite, en effectuant des opérations à coefficients dans $\K$ sur les lignes et les colonnes de $WQ$, on obtient une matrice $M \in \mathcal{M}_{n-\ell,n}(\K)$ de rang $n-\ell$  telle que $MWQ=0$. 

On note $U$ la matrice formée des $\ell$ premières colonnes de $MW$ et $V$ la matrice composée des $n-\ell$ dernières colonnes de $MW$. En écrivant le produit de matrices par blocs, on constate que $U\widetilde{R}+V\widetilde{S}=0$, de sorte que $U=-V\widetilde{S}\widetilde{R}^{-1}$. Donc $M=(-V\widetilde{S}\widetilde{R}^{-1} \quad V)$. 

Si $u={}^t (u_1, \dots, u_n) \in \mathcal{M}_{n,1}(\K(z))$, alors $$Mu=V \left( (-\widetilde{S}\widetilde{R}^{-1}) \begin{pmatrix}
u_1 \\ 
\vdots \\ 
u_{\ell}
\end{pmatrix}+ \begin{pmatrix}
u_{\ell+1} \\ 
\vdots \\ 
u_n
\end{pmatrix} \right) \in \Img(V).$$ Donc le rang sur $\K(z)$ de $M$ est inférieur au rang de $V$. Or, $M$ est de rang $n-\ell$ sur $\K$ donc également sur $\K(z)$ ; par suite, $V$ est de rang $n-\ell$ et est donc inversible. De plus, $\widetilde{S}\widetilde{R}^{-1}=-V^{-1}U$.

\medskip

Mais si $W=(w_{ij})_{1 \leqslant i,j \leqslant n}$, alors on voit avec la formule $V^{-1}=\gf{1}{\det V} (\mathrm{com}(V))^{T}$ (formule de la comatrice) que les coefficients de $V^{-1}U$ sont des fractions rationnelles en les $(w_{i,j})$ de degré total borné par une constante qui ne dépend que de $n$. De plus, comme $\widetilde{S}\widetilde{R}^{-1} \in \mathcal{M}_{n-\ell,\ell}(\K(z))$, ces coefficients sont également des éléments de $\K(z)$.

Finalement, le lemme \ref{Shidlovskiifracanalytique} nous fournit une borne supérieure $N_0 >0$ ne dépendant que de $W$ (donc \emph{in fine} de $G$) sur les degrés totaux des coefficients de $V^{-1}U$, donc sur ceux des coefficients de $\widetilde{S}\widetilde{R}^{-1}$.\end{dem}

\begin{lem}[\cite{Shidlovskii}, p. 85] \label{ordannulCL}
Soient $f_1, \dots, f_n \in \C\llbracket z \rrbracket$ tous non nuls et $R_1, \dots, R_n \in \C(z)$ de degrés totaux bornés par $c >0$. Alors $R \cdot f=\sum\limits_{i=1}^n R_i f_i$ a un zéro d'ordre au plus $d >0$, où $d$ est une constante dépendant uniquement des $f_i$ et de $c$. 
\end{lem}
\begin{dem}
\textbf{Étape 1} : réduction au cas où $\forall 1 \leqslant i \leqslant n, R_i \in \C[z]$.

Écrivons $R_i$ sous forme réduite $\gf{P_i}{Q_i}$, où $\deg(P_i) + \deg(Q_i) \leqslant c$. Alors $$R \cdot f=\gf{1}{Q_1 \dots Q_n} \sum\limits_{i=1}^n \left(\prod\limits_{j \neq i} Q_j \right) P_i f_i$$ et pour tout $i$, $\deg \left(\prod\limits_{j \neq i} Q_j \right) P_i \leqslant (n-1)c+c=nc$. De plus,

$$\ord_{z=0}(R \cdot f)=\ord \left(\sum\limits_{i=1}^n \left(\prod\limits_{j \neq i} Q_j \right) P_i f_i \right)- \ord(Q_1 \dots Q_n) \leqslant \ord \left(\sum\limits_{i=1}^n \left(\prod\limits_{j \neq i} Q_j \right) f_i\right).$$ On peut donc se ramener au cas où $\forall i \in \{ 1, \dots, n \}, R_i \in \C[z]$.

\medskip

\textbf{Étape 2} : réduction au cas où $\forall 1 \leqslant i \leqslant n, R_i \in \C$. On écrit donc pour tout $ i \in \{ 1, \dots, n \}, R_i=\sum\limits_{j=0}^c r_{j,i} z^j$, avec $r_{j,i} \in \C$. Donc $R \cdot f=\sum\limits_{i=1}^n \sum\limits_{j=0}^c r_{j,i} z^j f_i=\sum\limits_{1 \leqslant n \atop 0 \leqslant j \leqslant c} r_{j,i} g_{i,j}$, où $g_{i,j}(z)=z^j f_i(z) \in \C \{ z \}$. On remarque que les $g_{i,j}$ sont tous non nuls car les $f_i$ le sont.

\medskip

\textbf{Étape 3} : traitons le cas $\forall 1 \leqslant i \leqslant n, R_i \in \C$.

Introduisons $\mathcal{E}$ l'ensemble des $R \in \C^n$ tel que $R \cdot f$ n'est pas identiquement nul. Remarquons que si $g_1, \dots, g_{\ell} \in \C \{ z \}$ ont des ordres d'annulation en $0$ deux à deux distincts, alors $(g_1, \dots, g_{\ell})$ est libre sur $\C$. 

En effet, supposons que $\ord(g_1)=n_1 < \ord(g_2)=n_2 < \dots < \ord(g_{\ell})=n_{\ell}$. Supposons sans perte de généralité que $\alpha_1 \neq 0$. Alors, comme $\ord$ est une valuation sur $\C\llbracket z \rrbracket$, $$g_1=-\sum\limits_{i=2}^s \gf{\alpha_i}{\alpha_1} g_i $$ est d'ordre au moins $\min\limits_{i \geqslant 2, \alpha_i \neq 0} \ord(g_i) \geqslant n_2$ ce qui est absurde. Donc $\alpha_1=0$.

En itérant le procédé, on obtient $\alpha_1=\dots=\alpha_{s}=0$, donc $(g_1, \dots, g_{s})$ est libre sur $\C$.

\medskip

Notons $\nu$ le maximum des $k \in \N$ tel que $\mathcal{E}$ contient $k$ éléments d'ordre deux à deux distincts. Selon ce qui précède, $\nu \leqslant \dim_{\C} \C^n=n$. Soient $R_1 \cdot f, \dots, R_{\nu} \cdot f \in \mathcal{E}$ d'ordres deux à deux distincts.

Si $R \in \mathcal{E}$, $(R \cdot f, R_1 \cdot f, \dots, R_{(\nu)} \cdot f)$ est une famille à $\nu+1$ éléments, donc elle contient au moins deux éléments d'ordre égal. D'où $\ord(R \cdot f) \in \left\{ \ord(R_{(i)} \cdot f), 1 \leqslant i \leqslant \nu \right\}$.\end{dem}

\subsubsection{Démonstration du théorème \ref{Shidlovskii}}

On peut à présent passer à la preuve du théorème \ref{Shidlovskii}. On procède par l'absurde. 

\medskip

Supposons que $\Delta(z)$ est identiquement nul. Notons $\ell$ le rang sur $\K(z)$ de $(\mathbf{P}, \mathbf{P}_1, \dots, \mathbf{P}_{n-1})$ et $F=(\mathbf{P}, \mathbf{P}_1, \dots, \mathbf{P}_{\ell-1}) \in \mathcal{M}_{n,\ell}(\K[z])$. Alors en reprenant les notations du lemme \ref{lemmebaker}, on a $\mathbf{P}_m=T^m \widetilde{\mathbf{P}}_m$, de sorte que si $D=\mathrm{Diag}(1,T, \dots, T^{\ell-1})$, on a $F=\widetilde{F}D$, donc comme $\widetilde{F}$ est de rang $\ell$, il en va de même de $F$. Quitte à renuméroter les composantes de $\mathbf{f}$, on peut supposer que la matrice $R \in \mathcal{M}_{\ell}(\K(z))$ formée des $\ell$ premières lignes de $F$ est inversible. On note $S \in \mathcal{M}_{n-\ell,\ell}(\K(z))$ la matrice formée des $n-\ell$ dernières lignes de~$F$.

Le lemme \ref{lemmebaker} ci-dessous nous assure que les degrés totaux des coefficients de $SR^{-1}$ sont bornés par une constante $c_0 >0$ qui ne dépend que de $G$. En effet, toujours avec les notations du lemme \ref{lemmebaker}, on a $R=\widetilde{R}D$ et $S=\widetilde{S}D$, donc $SR^{-1}=\widetilde{S} \widetilde{R}^{-1}$.

Considérons 
$$H=\begin{pmatrix}
f_n(z) & 0 & \dots & 0 & 0 & \dots & -f_1(z) \\ 
f_2(z) & -f_1(z) & \dots & 0 &  &  &  \\ 
\vdots &  & \ddots & \vdots &  & (0) &  \\ 
f_{\ell}(z) & 0 & \dots & -f_1(z) & 0 & \dots & 0
\end{pmatrix}. $$
On note $H_0 \in \mathcal{M}_{\ell}(\K(z))$ la matrice composée des $\ell$ premières colonnes de $H$ et $H_1$ la matrice de $\mathcal{M}_{\ell,n-\ell}(\K(z))$ formée des $n-\ell$ dernières colonnes de $H$. Notons $U=HF$. On a 
$$\forall (i,j) \in \{ 1, \dots, \ell \}^2, \quad U_{ij}=  \begin{cases} \mathbf{P}_{1,j-1} f_n(z)-\mathbf{P}_{n,j-1} f_1(z) & \text{si} \;\; i=1 \\
\mathbf{P}_{1,j-1} f_i(z)-\mathbf{P}_{i,j-1} f_1(z) & \text{si} \;\; i \geqslant 2. \end{cases}$$
Un calcul de matrices par blocs nous montre que $U=H_0 R+H_1 S$, donc $UR^{-1}=H_0+H_1 SR^{-1}$.

Notons $B=SR^{-1}=(B_{i,j})_{i,j} \in \mathcal{M}_{n-\ell,\ell}(\K(z))$. La matrice $$H_1 B=(-f_1(z) B_{n-\ell,1}, \dots, -f_1(z) B_{n-\ell,\ell})=(-b_1 f_1(z), \dots, -b_{\ell} f_{1}(z)),$$ avec $b_i=B_{n-\ell,i} \in \K(z)$, est de degré total borné par $c_0$. Donc $$ UR^{-1}=\begin{pmatrix}
f_n(z)-b_1 f_1(z) & -b_2 f_1(z) & \dots & -b_{\ell} f_1(z) \\ 
f_2(z) & -f_1(z) & \dots & 0 \\ 
\vdots & \vdots & \ddots & \vdots \\ 
f_{\ell}(z) & 0 & \dots & -f_1(z)
\end{pmatrix}.$$ Par suite, en effectuant les opérations sur les lignes $L_1 \leftarrow L_1 - b_i L_i$ pour $i \in \{ 2, \dots, \ell \}$, on obtient
\begin{align}
\det(UR^{-1})&= \left\vert\begin{array}{cccc}
f_n(z)-\sum\limits_{i=1}^{\ell} b_1 f_i(z) & 0 & \dots & 0 \\ 
f_2(z) & -f_1(z) & \dots & 0 \\ 
\vdots & \vdots & \ddots & \vdots \\ 
f_{\ell}(z) & 0 & \dots & -f_1(z)
\end{array}\right\vert \notag \\
&= \left(f_n(z)-\sum\limits_{i=1}^{\ell} b_i f_i(z)\right) \times \det(-f_1(z) I_{\ell-1}) \label{calculdetur}
\end{align} en développant selon la première colonne. Donc $$\det(UR^{-1})=\left(f_n(z)-\sum\limits_{i=1}^{\ell} b_i f_i(z)\right) \times (-f_1(z))^{\ell-1} \neq 0,$$ car $(f_1(z), \dots, f_n(z))$ est supposée libre sur $\K(z)$.

\medskip

Pour aboutir à une contradiction, on va minorer puis majorer l'ordre d'annulation en $0$ de $\det(UR^{-1})$. D'abord, si $j \in \{ 1, \dots, \ell \}$ et $i \in \{ 2, \dots, n \}$, on a $$\begin{cases} Q_{j-1} f_1(z) - \mathbf{P}_{1,j-1} &=O\left(z^{N+M-(j-1)}\right)=O\left( z^{N+M-\ell}\right) \\
Q_{j-1} f_i(z) - \mathbf{P}_{i,j-1} &=O\left(z^{N+M-\ell}\right). \end{cases}$$

Donc en multipliant la première ligne par $f_i(z)$, la deuxième par $-f_1(z)$ et en les additionnant, on obtient que $$\mathbf{P}_{i,j-1} f_1(z) - \mathbf{P}_{1,j-1} f_i(z)=O\left(z^{N+M-\ell}\right).$$ Or, $\det(U)$ est une somme de produits de $n$ termes de cette forme, donc $$\det(U)=O\left(z^{\ell(N+M-\ell)}\right).$$

Mais comme on l'a montré dans l'étape 4 de la sous-section \ref{subsec:thchudnovsky}, si $j \in \{ 1, \dots, \ell \}$, les composantes de $\mathbf{P}_{j-1}$ sont de degrés inférieurs à $N+t(j-1) \leqslant N+t(\ell-1)$. Donc les coefficients de la matrice $R$ sont de degrés inférieurs à $N+t(\ell-1)$. 
Ainsi, $\det R \in \K[z]$ a un degré inférieur à $\ell(N+t(\ell-1))$, de sorte que $\det(UR^{-1})=\det U (\det R)^{-1}$ a un ordre d'annulation en $0$ au moins $$\ell(N+M-\ell)-\ell(N+t(\ell-1))=\ell[M-(t+1)(\ell-1)+1] \geqslant M-(n-1)(t+1)+1$$ car $\ell \geqslant 1$. D'autre part, comme $b_1, \dots, b_{\ell} \in \K(z)$ sont de degrés totaux bornés par $c_0$, le lemme~\ref{ordannulCL} ci-dessous nous assure, puisque les $f_i(z)$ sont tous non nuls, que l'ordre d'annulation en 0 de $f_n(z)-\sum\limits_{i=1}^\ell b_1 f_i(z)$, donc de $\det(UR^{-1})$ selon \eqref{calculdetur}, est borné par une constante $c_1 >0$ dépendant uniquement du système $\mathbf{y}'=G \mathbf{y}$ et de $\mathbf{f}$.

Par conséquent, si $M-(n-1)(t+1)+1 > c_1$, on aboutit à une contradiction. Donc $\Delta(z) \neq 0$, ce qui conclut la preuve du théorème \ref{Shidlovskii}.
\newpage

\section{Lien entre la nilpotence et la condition de Galochkin} \label{sec:nilpotencegalochkin}

Le but de cette partie est de définir une notion de nilpotence pour les systèmes différentiels, généralisant celle définie pour les opérateurs différentiels dans la sous-section \ref{subsec:thkatz}, que nous caractériserons à l'aide de rayons de convergence $p$-adiques. Dans un second temps, nous démontrerons que tout système différentiel satisfaisant la condition de Galochkin est nilpotent.

La référence principale pour cette section est \cite{Dwork}.

\subsection{Systèmes différentiels nilpotents}

\subsubsection{Nilpotence locale}\label{subsubsec:nilpotencelocale}

Dans cette partie,  nous allons donner une généralisation de la notion de nilpotence définie dans la sous-section \ref{subsec:nilpotentkatzhonda} et étudier ses propriétés.

Soit $H$ un corps de caractéristique $p>0$. On note $\mathcal{F}=H(z)$ et $D=\gf{\mathrm{d}}{\mathrm{d}z}$. On considère $G \in \mathcal{M}_n(H(z))$ et le système différentiel associé $y'=Gy$. 

On rappelle qu'on peut définir une suite de matrices $(G_s)_{s \in \N}$ par $$
G_{s+1}=G_s G+G'_s,$$ où $G'_s$ désigne la matrice $G_s$ dérivée coefficient par coefficient, de manière à ce que pour tout $y$ vérifiant $y'=Gy$, on ait $y^{(s)}=G_s y$.

\begin{defi} \label{systdiffnilpotent}
On dit que $y'=Gy$ \emph{est un système différentiel nilpotent} si $G_p$ est une matrice nilpotente.
\end{defi}

\begin{lem} \label{relationGp}
Pour tout $n \in \N$, $G_{pn}=(G_p)^n$.
\end{lem}

\begin{dem}
Soit $y$ un vecteur colonne tel que $y'=Gy$. Par définition, $y^{(p)}=G_p y$. Par la formule de Leibniz, on a donc

$$ y^{(2p)}=(y^{(p)})^{(p)}=\sum_{k=0}^p \binom{p}{k} G_p^{(k)} y^{(p-k)}=G_p y^{(p)}+G_p^{(p)} y $$ car si $1 \leqslant k \leqslant p-1, \binom{p}{k}$ est divisible par $p$. De plus comme $H$ est de caractéristique $p$ et $G \in \mathcal{M}_n(H(z))$, $(G_p)^{(p)}=0$, donc $y^{(2p)}=G_p y^{(p)}=G_p^2 y$. Ainsi, $G_{2p}=G_p^2$. On obtient le résultat annoncé par récurrence sur $n$ en écrivant $y^{(np)}=\left(y^{(n-1)p}\right)^{(p)}$.\end{dem}

\begin{prop} 
Soit $\Omega$ une extension de $H(z)$ sur laquelle $y'=Gy$ a une base de solutions (on dit que ce système est trivial sur $H(z)$) et $D^{p^s}$ est trivial, avec $pn \geqslant p^s$. On note $\Omega_s$ le corps $\mathrm{Sol}(D^{p^s}, \Omega)$. Alors $y'=Gy$ est nilpotent si et seulement si $y'=Gy$ est trivial sur $\Omega_s$.
\end{prop}

\begin{dem}
\begin{itemize}
\item Supposons que $y'=Gy$ est un système nilpotent. Soit $U \in \mathrm{GL}_n(\Omega)$ une matrice fondamentale de solutions de $y'=Gy$. Alors $G_p$ est nilpotent et classiquement, l'indice de nilpotence de $G_p$ est inférieur à $n$, donc $U^{(pn)}=G_{pn} U=G_p^n U=0$, donc comme $p^s \geqslant pn$, on a $U \in \GL_n(\Omega_s)$.

\item Si $y'=Gy$ est trivial sur $\Omega_s$, soit $U \in \GL_n(\Omega_s)$ une matrice fondamentale de solutions.\end{itemize} \noindent Alors $G_{p^s} U=U^{(p^s)}=0$, donc $0=G_{p^s}=G_p^{p^{s-1}}$, de sorte que $G_p$ est nilpotent.\end{dem}

\begin{rqu}
La notion de nilpotence ainsi définie généralise bien celle donnée dans la partie \ref{subsec:nilpotentkatzhonda}. Selon la proposition \ref{caracterisationnilpotenceomegas} et le résultat ci-dessus, un opérateur différentiel $L=D^n + a_1 D^{n-1}+ \dots + a_n \in H(z)[D]$ est nilpotent si et seulement si $y'=A_L y$ est un système différentiel nilpotent, où $$A_L=\begin{pmatrix}
0 & 1 &  & (0) \\ 
 & \ddots & \ddots &  \\ 
(0) &  & 0 & 1 \\ 
-a_n & \dots &  & -a_1
\end{pmatrix}$$ est la matrice compagnon associée à $L$.
\end{rqu}

\subsubsection{Nilpotence globale} \label{subsubsec:nilpotenceglobale}

Soit $\K$ un corps de nombres, on considère le système différentiel $y'=Gy$, $G \in \mathcal{M}_n(\K(z))$. Le but de cette partie est de définir une notion de nilpotence globale pour un tel système différentiel selon le procédé employé dans la sous-section \ref{subsec:thkatz}.
\medskip

Soit $\mathfrak{p}$ un premier de $\K$ définissant une valuation de Gauss sur $\K(z)$. Le lemme suivant est une conséquence de la proposition \ref{effet||0deriv} ci-dessous.

\begin{lem} \label{derivationgauss}
Si $f \in \K(z)$ est telle que $|f|_{\mathfrak{p},\mathrm{Gauss}} \leqslant 1$, alors $|f'|_{\mathfrak{p},\mathrm{Gauss}} \leqslant 1$.
\end{lem}

\begin{defi}
La \emph{norme de Gauss} d'une matrice $H \in \mathcal{M}_n(\K(z))$ est définie par $|H|_{\mathrm{Gauss}}=\sup\limits_{1 \leqslant i,j \leqslant n} |h_{ij}|_{\mathrm{Gauss}}$.
\end{defi}
\begin{rqu}
On vérifie que si $H_1, H_2 \in \mathcal{M}_n(\K(z))$, $|H_1 H_2|_{\mathrm{Gauss}} \leqslant |H_1|_{\mathrm{Gauss}} |H_2|_{\mathrm{Gauss}}$.
\end{rqu}

On déduit de la relation $G_{s+1}=G_s G+G'_s$ et du lemme \ref{derivationgauss} la proposition suivante

\begin{prop} \label{lienGsG}
\begin{itemizeth}
\item Si $|G|_{\mathfrak{p}, \mathrm{Gauss}} \leqslant 1$, alors $\forall s \in \N, |G_{s+1}|_{\mathfrak{p},\mathrm{Gauss}} \leqslant |G_s|_{\mathfrak{p},\mathrm{Gauss}}$. 
\item Pour tout $s \in \N$, $|G_s|_{\mathfrak{p},\mathrm{Gauss}} \leqslant \sup(1, |G|_{\mathfrak{p},\mathrm{Gauss}})^s$.
\end{itemizeth}
\end{prop}

Si $|G|_{\mathfrak{p},\mathrm{Gauss}} \leqslant 1$ (ce qui est vérifié pour presque tout premier $\mathfrak{p}$), alors selon la remarque suivant la proposition \ref{valuationgauss}, on peut réduire $G$ modulo $\mathfrak{p}$ coefficient par coefficient, si bien qu'on obtient une matrice $G \mod \mathfrak{p} \in \mathcal{M}_n(\overline{\K_{\mathfrak{p}}}(z))$ définissant un système différentiel $y'=(G \mod \mathfrak{p})y$. De plus, la proposition \ref{lienGsG} montre que pour tout $s \in \N$, $(G \mod \mathfrak{p})_s=G_s \mod \mathfrak{p}$.

\begin{defi} \label{defisystglobalementnilpotent}
On dit que $y'=Gy$ est un \emph{système globalement nilpotent} si pour tout $\mathfrak{p} \in \Spec(\Oal_\K)$ dans un ensemble de densité de Dirichlet $1$, $y'=(G \mod \mathfrak{p}) y$ est un système différentiel nilpotent.
\end{defi}

La définition et les propriétés de la densité de Dirichlet dans $\Spec(\Oal_{\K})$ sont présentées en annexe.

\begin{rqu}
Selon la remarque clôturant la partie \ref{subsubsec:nilpotencelocale} et la définition \ref{defglobalementnilpotent}, si $L \in \K(z)\left[\mathrm{d}/\mathrm{d}z\right]$, alors $L$ est globalement nilpotent si et seulement si le système différentiel $y'=A_L y$ est globalement nilpotent, où $A_L$ est la matrice compagnon associée à $L$.
\end{rqu}

\subsection{Rayon de convergence $p$-adique} \label{subsec:rcvpadique}

Soit $\K$ un corps de caractéristique $0$ muni d'une valeur absolue non archimédienne $| \cdot |$. On note $A$ l'anneau de valuation de $\K$ muni de son unique idéal maximal $\mathfrak{p}$ et $\overline{\K}=A/\mathfrak{p}$ le corps résiduel de $\K$, qu'on suppose de caractéristique $p >0$. Par exemple, on prend $\K$ un corps de nombres, $\mathfrak{p}$ un premier de $\K$, et $| \cdot |=| \cdot |_{\mathfrak{p}}$.

On s'intéresse au système différentiel $y'=Gy$, avec $G \in \mathcal{M}_n(\K(z))$. Soit $t$ une variable libre indépendante de $z$ sur $\K$ et $\Omega$ une extension complète et algébriquement close de $\K$ telle que $t \in \Omega$. Concrètement, $\Omega$ est le complété de la clôture algébrique du complété de $\K(t)$ muni de la valuation de Gauss (voir \cite[p. 24]{Dwork} pour une preuve de l'existence de $\Omega$). On note $B$ l'anneau de valuation de $\Omega$ et $\overline{\Omega}$ son corps résiduel.

Le but de cette section est de démontrer le théorème \ref{rcvnilpotence} donnant une caractérisation de la nilpotence du système réduit $y'=(G \mod \mathfrak{p}) y$ à l'aide d'une inégalité sur le rayon de convergence $p$-adique d'une matrice de solutions au point \emph{générique} $t$ du système $y'=Gy$.

\bigskip

\begin{rqu}
La réduction $\overline{t} \in \overline{\Omega}$ est transcendante sur $\overline{\K}$. En particulier, $\overline{\Omega}$ est infini.

En effet, soit $P \in A[x]$ tel que $\overline{P}(\overline{t})=0=\overline{P(t)}$. Alors $|P(t)|_{\mathrm{Gauss}} <1$, de sorte que, par définition de la valuation de Gauss, $P \in \mathfrak{p}[x]$, donc $\overline{P}=0$.

\end{rqu}

\medskip

En définissant pour $s \in \N$, $G_s$ comme ci-dessus, une matrice fondamentale de solutions $\mathcal{U} \in \GL_n(\Omega((z)))$ au point $t$ de $y'=Gy$ -- c'est-à-dire que $\mathcal{U}'(t)=G(t) \mathcal{U}(t)$, avec la dérivation $\left(\sum\limits_{s=0}^{\infty} A_s z^s\right)'=\sum\limits_{s=1}^{\infty} sA_s z^{s-1}$ sur $\Omega\llbracket z\rrbracket$ -- vérifiant $\mathcal{U}(t)=I_n$ est donnée par 

$$ \mathcal{U}(z)=\sum_{s=0}^{\infty} \gf{G_s(t)}{s!} (z-t)^s. $$

\begin{defi}
Le \emph{rayon de convergence générique} du système $y'=Gy$ est le rayon de convergence dans $\Omega\llbracket z\rrbracket$ de la série définissant $\mathcal{U}$.
\end{defi}

Si $|G|_{\mathrm{Gauss}} \leqslant 1$, alors comme dans la sous-section précédente, on peut définir un système réduit modulo $\mathfrak{p}$, $y'=\overline{G}y$, en caractéristique $p$.

\begin{Th} \label{rcvnilpotence}
Supposons que $|G|_{\mathrm{Gauss}} \leqslant 1$. Soit $R$ le rayon de convergence générique du système $y'=Gy$ et $\pi$ un élément de $\Omega$ algébriquement clos tel que $\pi^{p-1}=-p$. Alors le système $y'=\overline{G}y$ est nilpotent si et seulement si $R > |\pi|$.
\end{Th}

Notons que $|\pi|=|p|^{\frac{1}{p-1}}=p^{-\frac{1}{p-1}}$. Pour la preuve, on utilisera les deux lemmes suivants. Le premier a été démontré par Kummer dans \cite[p. 115]{Kummer}  :

\begin{lem} \label{ordrefactorielle}
Soit $p$ est un nombre premier et $n=a_0+a_1 p+ \dots +a_{\ell} p^{\ell}\in \N$, avec $0 \leqslant a_i < p-1$ et $a_{\ell} \neq 0$ un entier écrit en base $p$. Alors si $S_n=a_0+\dots+a_{\ell}$, on a $$v_p(n!)=\gf{n-S_n}{p-1}=\sum\limits_{k=0}^{\ell}\left\lfloor \gf{n}{p^k} \right\rfloor .$$
\end{lem}

\begin{dem}
Pour tout $j \in \{0, \dots, \ell \}$, le nombre de multiples de $p^j$ dans $\{0, \dots, n \}$ est $\left\lfloor \gf{n}{p^j} \right\rfloor$. Donc, en comptant $j$ fois les multiples de $p^j$, il y a $$\left\lfloor \gf{n}{p} \right\rfloor + \left\lfloor \gf{n}{p^2} \right\rfloor+ \dots + \left\lfloor \gf{n}{p^\ell} \right\rfloor$$ multiples de $p$ \og{}avec multiplicité\fg{} dans $\{0, \dots, n\}$. Or, pour tout $j$, $$\gf{n}{p^j}=\gf{a_0}{p^j}+\gf{a_1}{p^{j-1}}+\dots+\gf{a_{j-1}}{p}+a_j+\dots+a_{\ell} p^{\ell-j}$$ et $$0 \leqslant \gf{a_0}{p^j}+\gf{a_1}{p^{j-1}}+\dots+\gf{a_{j-1}}{p} \leqslant \gf{p-1}{p} \gf{1-\frac{1}{p^{j}}}{1-\frac{1}{p}} =1- \gf{1}{p^j} \leqslant 1,$$ de sorte que $\left\lfloor \gf{n}{p^j} \right\rfloor = a_j+a_{j+1} p+ \dots + a_{\ell} p^{\ell-j}$. Donc 

\begin{align*}
v_p(n!) &= \sum_{k=0}^n v_p(k)=\sum_{j=1}^\ell \left\lfloor \gf{n}{p^j} \right\rfloor = \sum_{j=1}^\ell a_j+a_{j+1} p+ \dots + a_{\ell} p^{\ell-j} \\
&= a_1+a_2(1+p)+a_3(1+p+p^2)+\dots + a_{\ell} (1+p+\dots+p^{\ell-1}) \\
&= a_1+a_2 \gf{p^2-1}{p-1} + \dots + a_{\ell} \gf{p^\ell-1}{p-1} =\gf{1}{p-1}(pa_1+\dots + p^\ell a_{\ell}+a_0 - S_n)= \gf{n-S_n}{p-1}.
\end{align*}\end{dem}

On munit $\Omega((z))$ de la dérivation $D=\gf{\mathrm{d}}{\mathrm{d}z}$.

\begin{lem} \label{valuationderiveepieme}
Si $f(z)=\gf{g(z)}{h(z)} \in \Omega(z)$, alors si $v$ est la valuation associée à $| \cdot |_{\mathrm{Gauss}}$, $v(D^p f) \geqslant 1 + v(f)$.
\end{lem}

\begin{dem}
Selon la remarque suivant la proposition \ref{effet||0deriv}, on a $|D^p f|_{\mathrm{Gauss}} \leqslant |p!| |f|_{\mathrm{Gauss}}$, donc $v(D^p f) \geqslant v_p(p!)+v(f) \geqslant 1+v(f)$.
\end{dem}

\begin{dem}[du théorème \ref{rcvnilpotence}]
\begin{itemize}
En vertu du lemme \ref{relationGp} , il faut montrer que $\overline{G}_{pn}=0$ si et seulement si $R > |\pi|$.

\item Supposons que $R > |\pi|$. Soit $R > R_1 > |\pi|$. Alors par définition du rayon de convergence, on a $$\left\vert \gf{G_s}{s!} \right\vert_{\mathrm{Gauss}} R_1^s \xrightarrow[s \rightarrow +\infty]{} 0.$$ Comme $\overline{\K}$ est de caractéristique $p$, par le théorème d'Ostrowski, $| \cdot |_{|\Q}=| \cdot |_p$. Donc selon le lemme \ref{ordrefactorielle}, $|(p^j)!|= p^{-\frac{p^j-1}{p-1}}$.

Comme $R_1 > |\pi|$, on peut trouver $\nu >1$ tel que $R_1=\nu |\pi|=\nu p^{-\frac{1}{p-1}}$, d'où

$$ \gf{R_1^{p^j}}{|(p^j)!|} \geqslant p^{\frac{p^j-1}{p-1}} \nu^{p^j} \left(p^{-\frac{1}{p-1}}\right)^{p^j}= \nu^{p^j} p^{-\frac{1}{p-1}} \xrightarrow[j \rightarrow +\infty]{} + \infty, $$
de sorte que $|G_{p^j}|_{\mathrm{Gauss}} \xrightarrow[j \rightarrow + \infty]{} 0$. Mais comme $|G|_{\mathrm{Gauss}} \leqslant 1$, selon la proposition \ref{lienGsG}, la suite $(|G_s|_{\mathrm{Gauss}})_{s}$ est décroissante, d'où $|G_s|_{\mathrm{Gauss}} \xrightarrow[s \rightarrow + \infty]{} 0$. 

En particulier, il existe $s_0 \in \N$ tel que pour tout $s \geqslant s_0$, $|G_s|_{\mathrm{Gauss}} <1$. Donc si $p \mu > s_0$, la réduction modulo $\mathfrak{p}$ de $G_{p \mu}$ est nulle. Par suite, selon le lemme \ref{relationGp}, le système réduit $y'=\overline{G} y$ est nilpotent.

\item Supposons que $y'=\overline{G} y$ est nilpotent. Alors $\overline{G_{pn}}=\overline{G_p}^n=0$, de sorte que $|G_{pn}|_{\mathrm{Gauss}} <1$. On peut donc trouver $\sigma \in \N^*$ tel que $$v(G_{pn}) := \min\limits_{1 \leqslant i,j \leqslant n} \gf{-\log|(G_{pn})_{i,j}|_{\mathrm{Gauss}}}{\log p} \geqslant \gf{1}{\sigma}.$$ On veut montrer par récurrence sur $s \in \N$ que

\begin{equation} \label{ordregpns}
 \forall s \in \N, \quad v(G_{pns}) \geqslant \gf{s}{\sigma}.
 \end{equation}

Soit $s \in \N$, supposons \eqref{ordregpns} vrai pour $s$. Soit $\mathcal{U}$ la matrice fondamentale de solutions générique au point $t$ de $y'=Gy$. 

Alors il est immédiat que 

\begin{align*}
G_{pn(s+1)}(t) \mathcal{U}(t) &= \left(D^{pn}(G_{pns} \mathcal{U})\right)(t)\overset{\mathrm{Leibniz}}{=} \sum\limits_{i=0}^{pn} \binom{pn}{i} (D^i G_{pns})(t) (D^{pn-i} \mathcal{U})(t) \\
&= \sum\limits_{i=0}^{pn} \binom{pn}{i} (D^i G_{pns})(t) G_{pn-i}(t) \mathcal{U}(t).
\end{align*}

Donc comme $\mathcal{U}(t)=I_n$, en substituant $z$ à $t$ dans $\K(z,t)$, on a

$$G_{pn(s+1)}= \sum\limits_{j=0}^n \binom{pn}{pj} (D^{pj} G_{pns}) G_{pn-pj} +\sum\limits_{(p,i)=1} \binom{pn}{i} (D^i G_{pns}) G_{pn-i}.$$

Il s'agit pour conclure d'examiner l'ordre de chacun des termes de la somme : 

\begin{itemize}
\item Par hypothèse de récurrence, $v(G_{pns}) \geqslant \gf{s}{\sigma}$ et $v(G_{pn}) \geqslant \gf{1}{\sigma}$, donc $$v(G_{pns} G_{pn}) \geqslant v(G_{pns})+v(G_{pn}) \geqslant \gf{s+1}{\sigma}.$$

\item Pour $j \neq 0$, selon le lemme \ref{valuationderiveepieme} ci-dessous, on a $v(D^{pj} G_{pns}) \geqslant 1 + v(G_{pns})$, donc $$v((D^{pj} G_{pns}) G_{pn-pj}) \geqslant 1+\gf{s}{\sigma} \geqslant \gf{s+1}{\sigma}.$$

\item Si $i$ et $p$ sont premiers entre eux, alors $v_p\binom{pn}{i} \geqslant 1$, donc $$v\left( \binom{pn}{i} (D^i G_{pns}) G_{pn-i} \right) \geqslant 1+v(G_{pns}) \geqslant \gf{s+1}{\sigma}.$$
\end{itemize}

Ces trois arguments impliquent que $v(G_{pn(s+1)}) \geqslant \gf{s+1}{\sigma}$, ce qui montre \eqref{ordregpns} par récurrence. 

De plus, selon la proposition \ref{lienGsG}, $(v(G_{\ell}))_{\ell \in \N}$ est croissante, donc si $\ell \in \N$ est écrit sous la forme $\ell=pns+r, 0 \leqslant r < pn-1$, on a $$v(G_{\ell}) \geqslant v(G_{pns}) \geqslant \gf{s}{\sigma}=\left\lfloor \gf{s}{pn} \right\rfloor \gf{1}{\sigma}.$$

\bigskip

Comme $| \cdot |$ est non-archimédienne, si $x \in \Omega$, $\mathcal{U}(x)$ converge si et seulement si $$\left\vert \gf{G_s(t)}{s!} (x-t)^s \right\vert_{\mathrm{Gauss}} \xrightarrow[s \rightarrow + \infty]{} 0,$$ c'est-à-dire si $v(G_s(t))-v_p(s!)+s v(x-t) \xrightarrow[s \rightarrow + \infty]{} +\infty$. Mais selon le lemme \ref{ordrefactorielle}, $v_p(s!) \leqslant \gf{s}{p-1}$, donc 

$$ v(G_s(t))-v_p(s!)+s v(x-t) \geqslant \left\lfloor \gf{s}{pn} \right\rfloor \gf{1}{\sigma} + s v(x-t)- \gf{s}{p-1} \geqslant s\left(\gf{1}{\sigma pn}-\gf{1}{p-1}+v(x-t)\right). $$

Pour que $\mathcal{U}(x)$ converge, il suffit donc que $v(x-t)>\gf{1}{p-1}-\gf{1}{\sigma pn}=v(\pi p^{1/(\sigma pn)})$. Par suite, $R \geqslant |\pi| p^{1/(\sigma pn)} > |\pi|$, ce qu'on voulait démontrer.\end{itemize}\end{dem}

\begin{rqu}
Même si le système $y'=Gy$ n'est pas nilpotent, on peut montrer avec la proposition \ref{lienGsG} que $R \geqslant |\pi|$.
\end{rqu}

\subsection{Estimations sur le rayon de convergence}\label{subsec:estimationsrcv}

Cette sous-section s'appuie sur le chapitre IV de \cite{Dwork}. Le but est d'étudier le comportement de séries entières $p$-adiques sur le bord de leur disque ouvert de convergence à l'aide de la norme $\| \cdot \|_0$ qui sera introduite dans la proposition \ref{||0}, ce qui fournit d'une part la formule de Hadamard sur le rayon de convergence $p$-adique (corollaire \ref{Hadamard}), et d'autre part, à l'aide des propositions techniques 19 à 21, le théorème \ref{dworkrobba}, dû à Dwork et Robba, dont une conséquence cruciale est le théorème de Dwork-Robba pour les systèmes différentiels (théorème \ref{estimationrcvdwork}). L'intérêt de la norme $\| \cdot \|_0$ dans le cadre de la démonstration du théorème d'André-Bombieri est qu'elle coïncide pour les fractions rationnelles avec la valeur absolue de Gauss.

\bigskip

On se donne $(\Omega, | \cdot |)$ un corps valué de caractéristique $0$ non archimédien tel que $\overline{\Omega}$ est infini. En particulier, le corps $\Omega$ construit dans la sous-section précédente vérifie ces hypothèses selon la remarque initiale de la sous-section \ref{subsec:rcvpadique}.

On note $p$ la caractéristique de $\overline{\Omega}$, de sorte que $| \cdot |_{|\Q}=| \cdot |_{p}$, quitte à renormaliser.

\subsubsection{Formule de Hadamard}

\begin{prop} \label{similicauchy}
Soit $f(z)=\sum\limits_{s=m}^{\infty} A_s z^s \in \Omega((z))$. Soit $R \in G_{\Omega}=\{|\omega|, \omega \in \Omega \}$. Si $f(x)$ converge pour $|x|=R$, alors 

$$ \sup_{|x|=R} |f(x)|=\sup_{s \geqslant m} (|A_s| R^s).$$
\end{prop}

\begin{dem}
Notons $M=\sup\limits_{s \geqslant m} (|A_s| R^s)$. Comme $f(x)$ converge pour $|x|=R$, on a $$|A_s| R^s \xrightarrow[s \rightarrow + \infty]{} 0,$$ donc $M < \infty$.

D'une part, comme $| \cdot |$ est non archimédienne, on a pour tout $x$ tel que $|x|=R$, $|f(x)| \leqslant \sup\limits_{s \geqslant m} (|A_s| R^s) \leqslant M$.

Pour montrer la deuxième inégalité, supposons sans perte de généralité que $R=1$. Comme $(|A_s|)_s$ tend vers 0, il y a un nombre fini d'indices $s$ tels que $|A_s|=M$, de sorte qu'on peut écrire $f(z)=g(z)+h(z)$, avec $$g(z)=\sum\limits_{|A_s| <M} A_s z^s \in \Omega((z)) \quad \textup{et} \quad h(z)=\sum\limits_{|A_s| = M} A_s z^s \in \Omega\left[z, \gf{1}{z} \right].$$

Comme $|A_s| \xrightarrow[s \rightarrow + \infty]{} 0$, on peut trouver $M' < M$ tel que, pour tout $s$ vérifiant $|A_s| <M$, $|A_s| \leqslant M'$. D'où, $| \cdot |$ étant non archimédienne, pour $x$ tel que $|x|=1$, $|g(x)| \leqslant M'$.

Par ailleurs, $$\sup\limits_{|x|=1} |h(x)| = |h|_{\mathrm{Gauss}}=M.$$ En effet, $\sup\limits_{|x|=1} |h(x)| \leqslant M$ de manière évidente. D'autre part, soit $y \in \Omega$ tel que $|y|=M$ et $$h_1(z)=y^{-1} h(z) \in \Oal_\Omega[z],$$ où $\Oal_\Omega$ désigne l'anneau de valuation de $\Omega$, d'unique idéal premier $\mathfrak{p}_{\Omega}$. 

Soit $\overline{h_1}$ la réduction modulo $\mathfrak{p}_{\Omega}$ de $h_1$. Supposons $\overline{h_1}=0$. Alors pour tout $s$ tel que $|A_s|=M$, $y^{-1} A_s \in \mathfrak{p}_{\Omega}$, c'est-à-dire $|y^{-1} A_s| <1$, si bien que $|A_s| <M$, ce qui est absurde. Donc $\overline{h_1} \neq 0$. Comme $\overline{\Omega}$ est infini, il existe ainsi $t \in \Oal_{\Omega}$ tel que $|t|=1$ et $\overline{h_1}(\overline{t}) \neq 0$, \emph{i.e.} $|h_1(t)|=1$. D'où par définition de $h_1$, $|h(t)|=|y|=M$.

Finalement, $|f(t)|=\max(|g(t)|, |h(t)|)=M$, ce qui fournit la conclusion voulue.\end{dem}

\begin{coro} \label{Hadamard}
Le rayon de convergence de $f(z)=\sum\limits_{s=0}^{\infty} A_s z^s \in \Omega\llbracket z\rrbracket$ est donné par la formule d'Hadamard 

$$ R=\liminf_{s \rightarrow + \infty} |A_s|^{-1/s}.$$
\end{coro}

\begin{dem}
\begin{itemize}
\item Soit $R'$ tel que $f(x)$ converge pour $|x|=R'$. Alors selon la proposition \ref{similicauchy}, pour tout $s \in \N$, $|A_s| R^s \leqslant \sup\limits_{|x|=R} |f(x)|$. Donc en prenant la puissance $\gf{1}{s}$ de part et d'autre et en passant à la limite supérieure, on obtient 
$$ \left(\limsup\limits_{s \rightarrow + \infty} |A_s|^{1/s}\right) R' \leqslant 1, $$
de sorte que $R' \leqslant \liminf\limits_{s \rightarrow + \infty} |A_s|^{-1/s}$. Comme le rayon de convergence $R$ est la borne supérieure de l'ensemble des tels $R'$, on a $R \leqslant \liminf\limits_{s \rightarrow + \infty} |A_s|^{-1/s}$.

\item si $R' > \liminf\limits_{s \rightarrow + \infty} |A_s|^{-1/s}$, alors la contraposée du premier point nous donne qu'il existe $x$ de valeur absolue $R'$ tel que $f(x)$ ne converge pas pour $|x|=R'$, d'où $R' \geqslant R$. Ainsi, $R=\liminf\limits_{s \rightarrow + \infty} |A_s|^{-1/s}$.\end{itemize}\end{dem}

\begin{prop} \label{||0}
Soit $R \in G_{\Omega}$. Si $f(z)=\sum\limits_{s=0}^{\infty} A_s z^s \in \Omega\llbracket z\rrbracket$ est convergente pour $|x|=R$, on pose $|f|_0(R)=\sup\limits_{|x|=R} |f(x)|$. Alors l'ensemble des séries formelles convergeant pour $|x|=R$ forme un anneau intègre sur lequel $| \cdot |_0(R)$ est une valeur absolue. 
\end{prop}

\begin{dem}
Admis, voir \cite[p. 116]{Dwork} .\end{dem}

\begin{defi}
On définit $\mathcal{A}_0$ l'\emph{anneau des fonctions analytiques} sur $D(0,1)$, c'est-à-dire le sous-anneau de $\Omega\llbracket z\rrbracket$ constitué des séries convergeant pour $|x| <1$. On note $\mathcal{A}'_0=\mathrm{Frac}(\mathcal{A}_0)$.
\end{defi}

\begin{prop} \label{caracterisationnorme0}
Pour $f(z)=\sum\limits_{s=0}^{\infty} A_s z^s \in \mathcal{A}_0$, on définit $\| f \|_0=\limsup\limits_{R \rightarrow 1^{-}} |f|_0(R)$. Alors $\| f \|_0=\sup\limits_{s \in \N} |A_s|$.
\end{prop}

\begin{dem}
Selon la proposition \ref{similicauchy}, pour tout $R<1$, $|f|_0(R)=\sup\limits_{s \in \N} |A_s| R^s \leqslant \sup\limits_{s \in \N} |A_s|$, donc $\| f \|_0 \leqslant \sup\limits_{s \in \N} |A_s|$ par passage à la limite supérieure.

D'autre part, pour tout $s \in \N$ et $R <1$, $|A_s R^s| \leqslant |f|_0(R)$, donc $|A_s|=\limsup\limits_{R \rightarrow 1^{-}} |A_s| R^s \leqslant \| f \|_0$, d'où le résultat.\end{dem}

Le résultat suivant est une conséquence de la proposition \ref{||0}.

\begin{prop} \label{||0normesurB0}
La fonction $\| \cdot \|_0$ définit une valeur absolue non archimédienne sur l'anneau des \emph{fonctions analytiques bornées} $\mathcal{B}_0 :=\{ f \in \mathcal{A}_0 : \| f \|_0 < + \infty \}$ pouvant être étendue à $\mathcal{B}'_0=\mathrm{Frac}(\mathcal{B}_0)$.
\end{prop}

\begin{rqu}
Pour $f, g \in \mathcal{A}'_0$, on a $\| f g \|_0 \leqslant \sup(\| f \|_0, \| g \|_0)$, dès que l'une des quantités intervenant dans le terme de droite est finie.
\end{rqu}

\begin{prop}
Si $h(z) \in \Omega(z)$, alors $h \in \mathcal{B}'_0$ et $\| h \|_0 = |h|_{\mathrm{Gauss}}$.
\end{prop}

\begin{dem}
En effet, si $h=\gf{f}{g}, f, g \in \Omega[z]$, alors $\| h \|_0= \gf{\| f \|_0}{\| g \|_0}=\gf{|f|_0(1)}{|g|_0(1)}$ car $f(x)$ et $g(x)$ convergent pour $|x|=1$, donc selon la proposition \ref{caracterisationnorme0}, $\| h \|_0=\gf{|f|_{\mathrm{Gauss}}}{|g|_{\mathrm{Gauss}}}=|h|_{\mathrm{Gauss}}$.\end{dem}

On munit dans toute la suite de cette sous-section $\Omega((z))$ de la dérivation $D=\gf{\mathrm{d}}{\mathrm{d}z}$, avec $D_{|\Omega}=0$.

\begin{prop} \label{effet||0deriv}
Soit $\xi \in \mathcal{A}'_0$. Pour tout $n \in \N$, on a $\left\| \gf{D^n \xi}{s! \xi} \right\|_0 \leqslant 1$. 
\end{prop}

\begin{dem}
Soit $g \in \mathcal{A}_0$ tel que $g \xi \in \mathcal{A}_0$. Alors selon la formule de Leibniz,

$$ \gf{D^n \xi}{s! \xi}=\gf{D^n(g \xi)}{n! g \xi} - \sum\limits_{i=0}^{n-1} \gf{D^{n-i} g}{(n-i)! g} \gf{D^i \xi}{i! \xi}. $$

Ainsi, si l'on montre le résultat pour $\xi \in \mathcal{A}_0$, selon la remarque suivant la proposition \ref{||0normesurB0}, un raisonnement par récurrence sur $n$ prouvera le résultat pour $\xi \in \mathcal{A}'_0$.

\medskip

Soit donc $\xi=\sum\limits_{s=0}^{\infty} A_s z^s \in \mathcal{A}_0$ et $R<1$. Alors $\gf{D^n \xi}{n!}=\sum\limits_{s=n}^{\infty} \binom{s}{n} A_s z^{s-n}$, de sorte que, par la proposition \ref{similicauchy}, $$\left\vert \gf{D^n \xi}{n!} \right\vert_{0}(R) \leqslant \sup\limits_{s \in \N} (|A_s| R^{s-n})=\gf{1}{R^n} |\xi|_0(R).$$

En passant à la limite supérieure, on obtient $\left\| \gf{D^n \xi}{n! \xi} \right\|_0 \leqslant 1$.\end{dem}

\begin{rqu}
En particulier, la conjugaison des deux dernières propositions donne pour $h \in \Omega(z)$ et $s \in \N$, $\left\vert \gf{D^s h}{s!} \right\vert_{\mathrm{Gauss}} \leqslant |h|_{\mathrm{Gauss}}$.
\end{rqu}

\bigskip
\bigskip

\subsubsection{Théorème de Dwork-Robba}

\begin{defprop}\label{majorationaccolade}
Si $s,m \in \N$, on définit
$$\{s, m \}_p = \sup\limits_{1 \leqslant \lambda_1 < \dots < \lambda_{m} \leqslant s} \gf{1}{|\lambda_1 \dots \lambda_{m}|_{p}},$$ avec la convention que $\{s, 0\}_p=1$ pour tout $s \in \N$.

Soient $s \in \N^*$, $m \in \N$ et $p$ premier. Alors si $p \leqslant s$, $\{s, m \}_p  \leqslant s^{m}$ et si $p >s$, $\{s, m \}_p=1$.
\end{defprop}

\begin{dem}
Le cas $m=0$ est évident puisque $\{s, 0\}_p=1$ par convention. On prend donc $m \geqslant 1$. Supposons $p \leqslant m$ Soit $1 \leqslant \lambda_1 < \dots < \lambda_{m} \leqslant s$, on a $|\lambda_1 \dots \lambda_{m} |_p \geqslant \left(\gf{1}{s}\right)^{m}$, puisque $$p^{\sum_{i=1}^{m} v_p(\lambda_i)} \leqslant \prod\limits_{i=1}^{m} p^{v_p(\lambda_i)} \leqslant \prod\limits_{i=1}^{m} \lambda_i \leqslant s^{m}.$$

Le deuxième cas est clair.\end{dem}

Le théorème suivant est le résultat principal de cette sous-section.

\begin{Th}[Dwork-Robba] \label{dworkrobba}
Soient $u_1, \dots, u_n \in \mathcal{A}'_0$ de matrice wronskienne $W(u_1, \dots, u_n)$ inversible. Pour tout $s \in \N$, on définit $(G_{s,0}, \dots, G_{s,n-1}) \in \mathcal{A'}_{0}^n$ tel que $$ \gf{D^s}{s!}(u_1, \dots, u_n)=(G_{s,0}, \dots, G_{s,n-1}) W(u_1, \dots, u_n).$$ Alors $$ \forall s \in \N, \;\; \forall 0 \leqslant j \leqslant n-1, \quad \| G_{s,j} \|_{0} \leqslant \{s,n-1 \}_p. $$
\end{Th}

\begin{rqu}
Le vecteur ligne $(G_{s,0}, \dots, G_{s,n-1})$ existe bien pour tout $s \in \N$ car les lignes $L_i=(D^i u_1, \dots, D^i u_n)$ pour $0 \leqslant i \leqslant n$ forment une famille libre, de sorte que $(D^n u_1, \dots, D^n u_n) \in \mathrm{Vect}_{\mathcal{A}'_0}(L_1, \dots, L_{n-1})$ ce qui fournit $(G_{n,0}, \dots, G_{n,n-1})$. Une récurrence immédiate conclut à l'existence des $G_{s,i}$ pour tout $s \in \N$.
\end{rqu}

\begin{dem}[du théorème \ref{dworkrobba}]
Procédons par récurrence sur $n$. L'initialisation pour $n=0$ est une reformulation de la proposition \ref{effet||0deriv}.

Soit $n \in \N$, supposons le résultat vrai pour $n$. On écrit $(u_1, \dots, u_{n+1})=u(1, \tau_1, \dots, \tau_n)$ ($u_1 \neq 0$ car $w(u_1, \dots, u_{n+1}) \neq 0$), de sorte que si $s \in \N$,

\begin{equation} \label{dworkrobbaeq1}
 \gf{D^s}{s!}(u_1 \dots u_{n+1})=\sum\limits_{i=0}^s \gf{u^{(i)}}{i!} \gf{D^j}{j!} (1, \tau_1, \dots, \tau_n).
 \end{equation}

Fixons $s \in \N^*$ et introduisons $U= \begin{pmatrix}
u_1 & \dots & u_{n+1} \\ 
\vdots & \ddots & \vdots \\ 
\gf{u_1^{(s)}}{s!} & \dots & \gf{u_{n+1}^{(s)}}{s!}
\end{pmatrix} \in \mathcal{M}_{s+1,n+1}(\mathcal{A}'_0)$. Alors selon \eqref{dworkrobbaeq1}, $$U=u P \begin{pmatrix}
1 & \tau_1 & \dots & \tau_n \\ 
0 & \tau'_1 & \dots & \tau'_n \\ 
\vdots & \vdots & \ddots & \vdots \\ 
0 & \gf{\tau_1^{(s)}}{s!} & \dots & \gf{\tau_n^{(s)}}{s!}
\end{pmatrix}, \;\; \text{où} \;\; P=\gf{1}{u} \begin{pmatrix}
u &  &  & \\ 
u' & u & (0) &  \\ 
\vdots & \vdots & \ddots &  \\ 
\gf{u^{(s)}}{s!} & \gf{u^{(s-1)}}{(s-1)!} & \dots & u
\end{pmatrix} \in \mathcal{M}_{s+1}(\mathcal{A}'_0).$$ Or, par hypothèse de récurrence, 
$$ \begin{pmatrix}
\tau'_1 & \dots & \tau'_n \\ 
\vdots & \ddots & \vdots \\ 
\gf{(\tau'_1)^{(s-1)}}{(s-1)!} & \dots & \gf{(\tau'_n)^{(s-1)}}{(s-1)!}
\end{pmatrix} =H W(\tau'_1, \dots, \tau'_n),$$ où $H \in \mathcal{M}_{s,n}(\mathcal{A}'_0)$ a tous ses coefficients sur la ligne $i+1$ bornées en norme $\| \|_0$ par $\{ i, n-1 \}_p$.

Définissons la matrice diagonale $\Delta=\mathrm{Diag}\left(1, \gf{1}{2}, \dots, \gf{1}{s} \right)$. Alors 

\begin{equation} \label{equationdworkrobba2}
U=uP \begin{pmatrix}
1 & \tau \\ 
0 & \Delta H W(\tau'_1, \dots \tau'_n)
\end{pmatrix}, \quad \tau=(\tau_1, \dots, \tau_n).
\end{equation}

Or, $U=G W$, où $G \in \mathcal{M}_{s+1,n+1}(\mathcal{A}'_0)$ de $(i+1)$-ème ligne $(G_{i,0} \dots G_{i,n})$, et $W=W(u_1, \dots, u_{n+1})$ est une matrice wronskienne. Par la formule de Leibniz,

\begin{align*}
W &= \begin{pmatrix}
u & u \tau_1 & \dots & u \tau_n \\ 
u' & (u \tau_1)' & \dots & (u \tau_n)' \\ 
\vdots & \vdots & \ddots & \vdots \\ 
u^{(n)} & (u \tau_1)^{(n)} & \dots & (u \tau_n)^{(n)}
\end{pmatrix} \\
&= \underset{uP_n}{\underbrace{\begin{pmatrix}
u & 0 &  &  &  \\ 
u' & u &  & (0) &  \\ 
u'' & \binom{2}{1} u' & u &  &  \\ 
\dots & \dots & \vdots & \vdots &  \\ 
u^{(n)} & \binom{n}{n-1} u^{(n-1)} & \binom{n}{n-2} u^{(n-2)} & \dots & u
\end{pmatrix}}} \begin{pmatrix}
1 &  \tau_1 & \dots &  \tau_n \\ 
0 & \tau'_1 & \dots &  \tau'_n \\ 
\vdots & \vdots & \ddots & \vdots \\ 
0 & (\tau_1)^{(n)} & \dots & (\tau_n)^{(n)}
\end{pmatrix},  \\
&&
\end{align*} de sorte que \begin{equation} \label{equationdworkrobba3}
W=uP_n \begin{pmatrix}
1 & \tau \\ 
0 &  W(\tau'_1, \dots \tau'_n)
\end{pmatrix}.
\end{equation}

En comparant \eqref{equationdworkrobba2} et \eqref{equationdworkrobba3} et en se souvenant que $U=GW$, on a 

$$ G P_n \begin{pmatrix}
1 & \tau \\ 
0 & W(\tau'_1, \dots \tau'_n)
\end{pmatrix} = P \begin{pmatrix}
1 & \tau \\ 
0 & \Delta H W(\tau'_1, \dots \tau'_n)
\end{pmatrix}.$$ Comme $(\tau'_1, \dots, \tau'_n)$ est libre (le wronskien de cette famille est non nul selon \eqref{equationdworkrobba3}), on en déduit que 
$$ G P_n= P \begin{pmatrix}
1 & \tau \\ 
0 & \Delta H W(\tau'_1, \dots, \tau'_n)
\end{pmatrix} \begin{pmatrix}
1 & \tau \\ 
0 &  W(\tau'_1, \dots, \tau'_n)
\end{pmatrix}^{-1}=\begin{pmatrix}
1 & -\tau W(\tau'_1, \dots, \tau'_n)^{-1} \\ 
0 &  W(\tau'_1, \dots, \tau'_n)^{-1}
\end{pmatrix}=P \begin{pmatrix}
1 & 0 \\ 
0 & \Delta H 
\end{pmatrix}.$$

Donc $G=P \begin{pmatrix}
1 & 0 \\ 
0 & \Delta H 
\end{pmatrix} P_n^{-1}$. Or, de la proposition \ref{effet||0deriv} découle que $\| P \|_0 \leqslant 1$ et $\| P_n \|_0 \leqslant 1$, donc comme $\det P_n=1$, on a $\| P_n^{-1} \|_0 = \| (\mathrm{com}(P_n))^{T} \|_0 \leqslant 1$.

Soit $h \leqslant s$, on sait que la $h$-ième ligne de $H$ est bornée pour $\| \cdot \|_0$ par $\{ h-1, n-1 \}_p$, donc la $h$\up{ème} ligne de $\Delta H$ est bornée pour $\| \cdot \|_0$ par $\gf{1}{|h|_p} \{h-1, n-1 \}_p \leqslant \{ h, n \}_p$.

Donc la $(h+1)$\up{ème} ligne de $\begin{pmatrix}
1 & 0 \\ 
0 & \Delta H 
\end{pmatrix}$ est bornée pour $\| \cdot \|_0$ par $\{ h, n \}_p$. Comme $\| P_n^{-1} \|_0 \leqslant 1$, c'est aussi le cas de la $(h+1)$\up{ème} ligne de $\begin{pmatrix}
1 & 0 \\ 
0 & \Delta H 
\end{pmatrix} P_n^{-1}$. Or, $P$ est triangulaire inférieure, donc $P \begin{pmatrix}
1 & 0 \\ 
0 & \Delta H 
\end{pmatrix} P_n^{-1}$ a pour coefficient d'ordre $(h+1,j)$ une somme de produits impliquant les coefficients de $\begin{pmatrix}
1 & 0 \\ 
0 & \Delta H 
\end{pmatrix} P_n^{-1}$ situés sur la colonne $j$ de numéro de ligne inférieure ou égale à $h+1$ et les coefficients de $P$ situés sur sa $h+1$\up{ème} ligne. Comme $h \mapsto \{ h, n \}_p$ est croissante, ce coefficient est borné par $\sup\limits_{1 \leqslant i \leqslant h} \gf{1}{|i|_p} \{ h-1, n-1 \}_p \leqslant \{h, n \}_p$. On a ainsi obtenu $$ \forall h \in \N, \forall j \in \{0, \dots, n-1 \}, \| G_{h+1, j} \|_0 \leqslant \{h, n \}_p,$$ ce qu'on voulait démontrer.\end{dem}

On peut obtenir également une estimation du type de celle du théorème de Dwork-Robba dans le cadre de l'étude d'un système différentiel. 

\begin{prop}[Dwork-Robba pour les systèmes] \label{dworkrobbasystemes}
Soit $A \in \mathcal{M}_n(\mathcal{A}'_0)$, et pour tout $i \in \N$, $A_i$ la matrice telle que si $y'=Ay$, alors $\gf{D^i}{i!} y=A_i y$. Soient $s \in \N$ et $Y \in \GL_n(\mathcal{A}'_0)$ tels que $\gf{D^s Y}{s!}=A_s Y$, $A_s \in \mathcal{M}_n(\mathcal{A}'_0)$. Alors

$$ \| A_s \|_0 := \sup\limits_{1 \leqslant i, j \leqslant n} \| (A_s)_{i,j} \|_0 \leqslant \{ s,n-1 \}_p \sup_{0 \leqslant i \leqslant n-1} (\| i! A_i \|_0).$$
\end{prop}

\begin{dem}
Il suffit de montrer le résultat pour la première ligne $u=(u_1, \dots, u_n) \in \Omega((z))^n$ de $A_s$.

Notons $q=\dim \Vect_{\Omega}(u_1, \dots, u_n)$ et prenons une base $z_1, \dots, z_q$ de ce $\Omega$-espace vectoriel. En notant $z=(z_1, \dots, z_q)$, alors il existe $E \in \GL_n(\Omega)$ tel que $uE=(z \quad 0_{n-p})$. Considérons la matrice wronskienne de $z$ $$Z=\begin{pmatrix}
z \\ 
\vdots \\ 
D^{q-1} z
\end{pmatrix} \in \mathcal{M}_{q-1}(\mathcal{A}'_0) \quad \text{et} \;\; T=\begin{pmatrix}
u \\ 
\vdots \\ 
D^{q-1} u
\end{pmatrix} Y^{-1}.$$ Prenons $G_s=(G_{s,0}, \dots, G_{s,q-1})$ tel que $\gf{D^s z}{s!}=G_s Z$. Alors, comme $D_{|\Omega}=0$, $$\left( \gf{D^s u}{s!} \right) E= \gf{D^s(uE)}{s!}=\gf{D^s}{s!} (z \quad 0)=G_s (Z \quad 0_{n-q}).$$ D'autre part, $$TYE= \begin{pmatrix}
uE \\ 
\vdots \\ 
(D^{q-1}u) E
\end{pmatrix}=\begin{pmatrix}
uE \\ 
\vdots \\ 
D^{q-1}(uE)
\end{pmatrix}=\begin{pmatrix}
z \quad 0 \\ 
\vdots \\ 
D^{q-1}z \quad 0
\end{pmatrix}=(Z \quad 0).$$

Donc $\left(\gf{D^s}{s!} u\right) E=G_s TYE$, donc, comme $E$ est inversible, $\gf{D^s u}{s!}=G_s TY$.

\medskip

Notons pour tout $i \in \N$, $L_i$ la première ligne de $A_i$. 
Comme $A_s=\left(\gf{D^s Y}{s!}\right) Y^{-1}$, on a $$L_s=\gf{D^s u}{s!} Y^{-1}=G_s T, \; \; \text{d'où} \;\; \| L_s \|_0 \leqslant \| G_s \|_0 \| T \|_0.$$ D'autre part, si $0 \leqslant i \leqslant n-1$, la $(i+1)$\up{ème} ligne de $T$ est $(D^i u)Y^{-1}=(i!) L_i$, donc $$\| T \|_0 \leqslant \sup\limits_{0 \leqslant i \leqslant n-1} (\| i! A_i \|_0).$$ 

Le théorème \ref{dworkrobba} donne finalement $$ \| L_s \|_0  \leqslant \| G_s \|_0 \| T \|_0 \leqslant \{s, q-1 \}_p \| T \|_0 \leqslant \{ s, n-1\}_p \| T \|_0 \leqslant \{ s,n-1 \}_p \sup\limits_{0 \leqslant i \leqslant n-1} \| i! A_i \|_0.$$\end{dem}

On cherche à obtenir une version plus générale du théorème précédent appliquée à des matrices $A$ définies par une série entière de rayon de convergence quelconque autour d'un point $a$ de $\Omega$.

On définit pour $a \in \Omega$ et $\rho >0$, $\mathcal{A}_{a, \rho}$ l'anneau des séries formelles $\sum\limits_{s=0}^{\infty} A_s (z-a)^s$ de rayon de convergence autour de $a$ supérieur ou égal à $\rho$. On peut le munir d'une norme $\| \cdot \|_{a, \rho}$ de manière à ce que, si $\alpha \in \Omega$ est tel que $|\alpha|=\rho$, l'isomorphisme d'anneaux

$$ \fonction{\phi_{a, \rho}}{\mathcal{A}_{a, \rho}}{\mathcal{A}_{0}}{f(z)}{f\left(\alpha(z+a)\right)}$$ soit une isométrie. On définit de même que précédemment $\mathcal{A}'_{a, \rho}$, $\mathcal{B}_{a, \rho}$, $\mathcal{B}'_{a, \rho}$. On remarque que $D(f(\alpha(z+a))) = \alpha f(\rho(z+a))$. L'utilisation de l'isométrie $\phi_{a, \rho}$ permet d'obtenir la version généralisée du théorème \ref{dworkrobbasystemes} suivante :

\begin{Th} \label{estimationrcvdwork}
Soit $A \in \mathcal{M}_n(\mathcal{A}'_{a, \rho})$, et pour tout $i \in \N$, $A_i$ la matrice telle que si $y'=Ay$, alors $\gf{D^i}{i!} y=A_i y$. Soit $s \in \N$, $Y \in \GL_n(\mathcal{A}'_{a, \rho})$ tel que $\gf{D^s Y}{s!}=A_s Y$, $A_s \in \mathcal{M}_n(\mathcal{A}'_{a, \rho})$. Alors  

$$\| A_s \|_{a, \rho} \leqslant \rho^{-s} \{s, n-1 \}_p \sup\limits_{0 \leqslant i \leqslant n-1} \left( \rho^i \| i! A_i \|_{a, \rho} \right).$$
\end{Th}

\begin{dem}
Posons $\tilde{Y}=Y(\alpha(z+a))=\phi_{a,\rho}(Y) \in \mathrm{GL}_n(\mathcal{A}'_0)$ via l'isomorphisme d'anneaux $\phi_{a, \rho}$. On a pour tout $i \in \N$, $$\gf{D^i \tilde{Y}}{i!}=\alpha^i \left(\gf{D^i Y}{i!}\right)(\alpha(z+a))=\alpha^i A_i(\alpha(z+a)) \tilde{Y}=\tilde{A}_i \tilde{Y},$$ avec $\tilde{A}_i=\alpha^i A_i(\alpha(z+a))$, donc selon la proposition \ref{dworkrobbasystemes}, 

$$\| \alpha^s A_s(\alpha(z+a)) \|_0 \leqslant \{s,n-1 \}_p \sup\limits_{0 \leqslant i \leqslant n-1} \| i! \tilde{A}_i \|_0=\{s,n-1 \}_p \sup\limits_{0 \leqslant i \leqslant n-1} \| \rho^i i! A_i(\alpha(z+a)) \|_0$$ ce qui donne le résultat voulu car $|\alpha|=\rho$.\end{dem}

Dans la partie suivante, on va appliquer ces résultats aux $G$-opérateurs dans le théorème d'André-Bombieri.

\subsection{Le théorème d'André-Bombieri}

Le but de cette partie est de démontrer le théorème d'André-Bombieri (théorème \ref{bombieriandre}) qui établit que la condition de Galochkin pour un système différentiel $y'=Gy$ implique la nilpotence globale de ce système. Pour cela, nous allons introduire une condition sur les rayons de convergence génériques $p$-adiques de $y'=Gy$, la \emph{condition de Bombieri}, qui implique la nilpotence globale de $y'=Gy$. De plus, nous allons reformuler la condition de Galochkin à l'aide d'une quantité $p$-adique, la \emph{taille} de $G$.

\subsubsection{Rayon de convergence global}
On considère un corps de nombres $\K$ de degré $\delta$.

\begin{defi}
On note, pour tout $\mathfrak{p} \in \Spec(\Oal_\K)$, $R_{\mathfrak{p}}$ le rayon de convergence de la matrice solution générique $\mathcal{U}_{G, t_{\mathfrak{p}}}$ issue de la construction de la sous-section \ref{subsec:rcvpadique} associée au système $y'=Gy$ et à la valeur absolue $|\cdot |_{\mathfrak{p}}$ sur $\K$.

Le \emph{rayon de convergence générique inverse global} de $G$ est 

$$ \rho(G) :=\sum\limits_{\mathfrak{p} \in \Spec(\Oal_\K)} \log^{+}\left(\gf{1}{R_\mathfrak{p}} \right).$$

On dit que $G$ \emph{satisfait la condition de Bombieri} si $\rho(G) < \infty$.
\end{defi}

\begin{prop}
Si la matrice $G$ satisfait la condition de Bombieri, alors $y'=Gy$ est un système différentiel globalement nilpotent.
\end{prop}

\begin{dem}
Soit $\mathcal{S}$ l'ensemble des nombres premiers $p$ tels qu'il existe $\mathfrak{p} \in \Spec(\Oal_\K)$ au-dessus de $p$ tel que $R_{\mathfrak{p}} \leqslant |p|_{\mathfrak{p}}^{\frac{1}{p-1}}$.  Selon le théorème \ref{rcvnilpotence}, $\mathcal{S}$ est l'ensemble des $p$ tels que $y'=(G \mod \mathfrak{p}) y$ est non nilpotent pour au moins un premier $\mathfrak{p}$ au-dessus de $p$. Fixons pour chaque $p \in \mathcal{S}$ un tel premier $\mathfrak{p}(p)$. Alors 

$$\rho(G) \geqslant \sum\limits_{p \in \mathcal{S}} \log^{+} \left(\gf{1}{R_{\mathfrak{p}(p)}} \right) \geqslant \sum\limits_{p \in \mathcal{S}} \gf{1}{\delta} \gf{\log p}{p-1}.$$ En effet, $|p|_{\mathfrak{p}}=p^{-f_{\mathfrak{p}}/\delta}$, avec la normalisation choisie.

Par suite, si $\rho(G) < \infty$, on a $$\sum\limits_{p \in \mathcal{S}} \gf{1}{p} \leqslant \sum\limits_{p \in \mathcal{S}} \gf{\log p}{p-1} < \infty.$$ Ainsi (cf définition \ref{densitedirichlet} en annexe), $\mathcal{S}$ a une densité de Dirichlet nulle, c'est-à-dire que
$$\gf{-1}{\log(s-1)} \sum\limits_{p \in \mathcal{S}} p^{-s} \rightarrow 0$$ quand $s$ tend vers $1$, $s>1$. Or, si $\mathfrak{p} \in \Spec(\Oal_\K)$ est tel que $y'=(G \mod \mathfrak{p}) y$ est non nilpotent, alors $\mathfrak{p} \cap \Z \in \mathcal{S}$, donc l'ensemble $\mathcal{S}'$ des premiers vérifiant cette propriété est de densité de Dirichlet nulle. En effet, si $s>1$,
$$ \sum\limits_{\mathfrak{p} \in \mathcal{S}'} \gf{1}{N(\mathfrak{p})^s}=\sum\limits_{p \in \Spec(\Z)} \sum\limits_{\mathfrak{p} \in \mathcal{S} \atop \mathfrak{p} \mid p} \gf{1}{p^{f_{\mathfrak{p}}s}} \leqslant \sum\limits_{p \in \mathcal{S}} \gf{\mathrm{Card}\{ \mathfrak{p} \in \mathcal{S}' : \mathfrak{p} \mid p\}}{p^s} \leqslant \delta \sum\limits_{p \in \mathcal{S}} \gf{1}{p^s}$$ car le nombre de premiers $\mathfrak{p}$ au-dessus de $p$ est borné par $\delta$ selon la formule $\sum\limits_{\mathfrak{p} \mid p} d_{\mathfrak{p}}=\delta$. Donc en divisant par $-\log(s-1)$ de part et d'autre de l'inégalité et en passant à la limite, on obtient que $\mathcal{S}'$ a une densité de Dirichlet nulle. Selon la définition \ref{defisystglobalementnilpotent}, le système $y'=Gy$ est donc globalement nilpotent.\end{dem}

\subsubsection{Taille d'un système différentiel}

On considère un corps de nombres $\K$ de degré $\delta$.
Pour tout $\mathfrak{p} \in \Spec(\Oal_\K)$, on définit la valeur absolue $| \cdot |_{\mathfrak{p}}$ de manière à ce que si $p$ est le premier de $\Z$ au-dessous de $\mathfrak{p}$, $|p|_{\mathfrak{p}}=p^{-\frac{d_{\mathfrak{p}}}{\delta}}$, où $d_{\mathfrak{p}}=e_{\mathfrak{p}} f_{\mathfrak{p}}$, avec $e_{\mathfrak{p}}$ (resp. $f_{\mathfrak{p}}$) le degré d'inertie (resp. de ramification) de $\mathfrak{p}$ au-dessus de $p$. 

En remarquant que $v_{\mathfrak{p}}(p)=e_{\mathfrak{p}}$, on a donc par définition 

$$ \forall x \in \K, \quad |x|_{\mathfrak{p}}=\left(N_{\K/\Q}(\mathfrak{p})\right)^{-\delta v_{\mathfrak{p}}(x)}=p^{-\frac{f_{\mathfrak{p}}}{\delta} v_{\mathfrak{p}}(x)}. $$

Les résultats de la sous-section \ref{subsec:estimationsrcv} n'utilisant par cette normalisation, on les appliquera en tenant compte de cette adaptation technique rendue nécessaire par la démonstration du lemme \ref{denomhauteur} ci-dessous.

\begin{defi} \label{denomprime}
Si $\alpha_1, \dots, \alpha_n$ tels que $(\alpha_i)=\gf{\mathfrak{a}_i}{\mathfrak{b}_i}$, où $\mathfrak{a}_i$ et $\mathfrak{b}_i$ sont des idéaux de $\Oal_\K$ premiers entre eux, on définit $\den'(\alpha_1, \dots, \alpha_n)=N_{\K/\Q}(\mathfrak{b})$, où $\mathfrak{b}$ est le plus petit multiple commun au sens des anneaux de Dedekind des idéaux $\mathfrak{b}_1, \dots, \mathfrak{b}_n$.
\end{defi}

\begin{rqu}
Ce dénominateur n'est pas nécessairement le plus petit dénominateur commun des $\alpha_i$ au sens de la définition \ref{maisonetdenominateur}, comme le montre l'exemple de $\alpha=\gf{1+i}{2}$ dans $\K=\Q(i)$ qui vérifie $\mathrm{den}(\alpha)=2$ et $\mathrm{den}'(\alpha)=4$, puisque $(1+i)$ est un idéal premier de $\Z[i]=\Oal_\K$. Cependant on peut obtenir un encadrement $$ \den(\alpha_1, \dots, \alpha_n) \leqslant \den'(\alpha_1, \dots, \alpha_n) \leqslant (\den(\alpha_1, \dots, \alpha_n))^{\delta}. $$ En effet, notons $N=\den(\alpha_1, \dots, \alpha_n)$ et $M=\den'(\alpha_1, \dots, \alpha_n)$. Alors :

\begin{itemize}
\item On sait par les propriétés de la norme que $M \in \mathfrak{b}$ (c'est clair si $\mathfrak{b}$ est un idéal premier, sinon on utilise la multiplicativité de la norme et la décomposition en idéaux premiers dans l'anneau de Dedekind $\Oal_\K$), donc $(M)=\mathfrak{b} \mathfrak{b}'$ pour un certain idéal $\mathfrak{b}'$, de sorte que pour tout $i \in \{1, \dots, n \}$, $(M \alpha_i)= \mathfrak{b} \mathfrak{b}' \gf{\mathfrak{a}_i}{\mathfrak{b}_i} =\mathfrak{b}' \mathfrak{a}_i \gf{\mathfrak{b}}{\mathfrak{b}_i}$, donc, par définition de $\mathfrak{b}$, $M \alpha_i \in \Oal_\K$. Donc $N \leqslant M$.

\item D'autre part, pour tout $1 \leqslant i \leqslant n, (N) \gf{\mathfrak{a}_i}{\mathfrak{b}_i} \subset \Oal_\K$, donc $N \mathfrak{a}_i \subset \mathfrak{b}_i$ pour tout $i$, donc, comme $\mathfrak{a}_i$ et $\mathfrak{b}_i$ sont premiers entre eux, $(N) \subset \mathfrak{b}_i$, de sorte que $(N) \subset \mathfrak{b}$, le plus petit commun multiple des $\mathfrak{b}_i$. Donc $M=N_{\K/\Q}(\mathfrak{b})= |\Oal_\K/\mathfrak{b}| \leqslant |\Oal_\K/N \Oal_\K|=N^{\delta}$.

\end{itemize}
\end{rqu}

\begin{defi}
On définit les fonctions $$\fonction{\log^{+}}{]0,+\infty[}{\R_{+}}{x}{\log(\sup(1,x))} \quad \mathrm{et} \quad \fonction{\log^{-}}{]0,+\infty[}{\R_{-}}{x}{\log(\inf(1,x)).} $$
\end{defi}

On remarque immédiatement que $$\forall x, y >0,\;\; \log^{+}(xy) \leqslant \log^{+}(x)+\log^{+}(y) \quad \text{et} \;\; \log^{+}(x^y) \leqslant y \log^{+}(x).$$

\begin{lem} \label{denomhauteur}
Si $\alpha_1, \dots, \alpha_n \in \K$, alors 

$$\sum\limits_{\mathfrak{p} \in \Spec(\Oal_\K)}\sup\limits_{1 \leqslant i \leqslant n} \log^{+} |\alpha_i |_{\mathfrak{p}} = \gf{1}{\delta} \log\left( \den'(\alpha_1, \dots, \alpha_n)\right).$$
\end{lem}

\begin{dem}
En reprenant les notations de la définition \ref{denomprime}, on a pour tout $i \in \{1, \dots, n \}$, $\mathfrak{b}_i= \prod\limits_{v_{\mathfrak{p}}(\alpha_i) <0} \mathfrak{p}^{-v_{\mathfrak{p}}(\alpha_i)}$, de sorte que $\mathfrak{b}=\prod\limits_{\exists i : v_{\mathfrak{p}}(\alpha_i) <0} \mathfrak{p}^{-\min\limits_{i}(v_{\mathfrak{p}}(\alpha_i))}$.

Donc $$N_{\K/\Q}(\mathfrak{b})=\prod\limits_{p \in \Spec(\Z)} \prod\limits_{\mathfrak{p} \mid p \atop \exists i : v_{\mathfrak{p}}(\alpha_i) <0} p^{-f_{\mathfrak{p}} (\min\limits_{i}(v_{\mathfrak{p}}(\alpha_i))}$$ (le produit a bien un nombre de termes fini). D'où

$$\gf{1}{\delta} \log\left(\mathrm{den}'(\alpha_1, \dots, \alpha_n)\right)=\sum\limits_{p \in \mathrm{Spec}(\Z)} \sum\limits_{\mathfrak{p} \mid p \atop \exists i : v_{\mathfrak{p}}(\alpha_i) <0} \gf{-f_{\mathfrak{p}}}{\delta} (\min\limits_{i}(v_{\mathfrak{p}}(\alpha_i)) \log p. $$

Or, pour tout $i \in \{ 1, \dots, n\}$, $|\alpha_i|_{\mathfrak{p}}=p^{-\frac{\delta}{f_{\mathfrak{p}}} v_{\mathfrak{p}}(\alpha_i)}$, de sorte que $$\sup\limits_{1 \leqslant i \leqslant n} \log^{+} |\alpha_i|_{\mathfrak{p}}=\gf{-f_{\mathfrak{p}}}{\delta} \log p \min\limits_{1 \leqslant i \leqslant n \atop v_{\mathfrak{p}}(\alpha_i) <0}  v_{\mathfrak{p}}(\alpha_i).$$

En remarquant que cette dernière quantité est nulle si $\forall i, v_{\mathfrak{p}}(\alpha_i) \geqslant 0$, on obtient le résultat voulu en sommant sur $\mathfrak{p} \in \Spec(\Oal_\K)$.\end{dem}

\begin{defi}
La \emph{taille} de $G$ est $$\sigma(G)=\limsup\limits_{s \rightarrow + \infty} \gf{1}{s} \sum\limits_{\mathfrak{p} \in \Spec(\Oal_\K)} h(s, \mathfrak{p}),$$ où 

$$ \forall s \in \N, \quad h(s, \mathfrak{p})= \sup_{m \leqslant s} \log^{+} \left\vert \gf{G_m}{m!} \right\vert_{\mathfrak{p}, \mathrm{Gauss}}. $$
\end{defi}

La taille de $G$ est une quantité qui \og{} encode \fg{} la condition de Galochkin (définition \ref{galochkin}),  dans un sens précisé par la proposition suivante :

\begin{prop} \label{proplienq'ssigmaG}
Soit $T(z) \in \K[z]$ tel que $T(z)G(z) \in \mathcal{M}_n(\K[z])$, et soit, pour $s \in \N$, $q_s$ (resp. $q'_s$), le dénominateur (resp. le $\mathrm{den}'$) des coefficients des coefficients des matrices $TG, T^2 \gf{G_2}{2}, \dots, T^s \gf{G_s}{s!}$. Alors 
\begin{enumerateth}
\item On a, pour tout entier $s$,  $\gf{\log q_s}{s} \leqslant \gf{\log q'_s}{s} \leqslant \delta \gf{\log q_s}{s}$. 
\item Notons $$h^{-}(T)=\sum\limits_{\mathfrak{p} \in \Spec(\Oal_\K)} \log^{-} |T|_{\mathfrak{p}, \mathrm{Gauss}} \quad \mathrm{et} \quad h^{+}(T)=\sum\limits_{\mathfrak{p} \in \Spec(\Oal_\K)} \log^{+} |T|_{\mathfrak{p}, \mathrm{Gauss}}.$$ Alors $\delta \sigma(G)+h^{-}(T) \leqslant \limsup\limits_{s \rightarrow + \infty} \gf{\log q'_s}{s}  \leqslant \delta \sigma(G) + h^{+}(T)$.
\end{enumerateth}
\end{prop}

\begin{dem}
Le point \textbf{a)} est une conséquence directe de la remarque suivant la proposition \ref{denomprime}. 

\textbf{b)} Selon le lemme \ref{denomhauteur}, on a $$\log q'_s=\delta \sum\limits_{\mathfrak{p} \in \Spec(\Oal_\K)} \sup\limits_{m \leqslant s} \log^{+} \left\vert\gf{T^m G_m}{m!} \right\vert_{\mathfrak{p}, \mathrm{Gauss}}.$$
De plus, si $s\in \N^*$ et $m \leqslant s$, alors $$\log^{+} \left\vert T^m \gf{G_m}{m!} \right\vert_{\mathfrak{p}, \mathrm{Gauss}} =\log^{+} \left( |T|^m_{\mathfrak{p}, \mathrm{Gauss}} \left\vert\gf{G_m}{m!} \right\vert_{\mathfrak{p}, \mathrm{Gauss}}\right)\leqslant s \log^{+} |T|_{\mathfrak{p}, \mathrm{Gauss}} + \log^{+} \left\vert \gf{G_m}{m!} \right\vert_{\mathfrak{p}, \mathrm{Gauss}}.$$

Donc $$\sup_{m \leqslant s} \log^{+} \left\vert T^m \gf{G_m}{m!} \right\vert_{\mathfrak{p}, \mathrm{Gauss}} \leqslant  s \log^{+} |T|_{\mathfrak{p}, \mathrm{Gauss}}+h(s,\mathfrak{p}),$$ si bien que $$\sum\limits_{\mathfrak{p} \in \Spec(\Oal_\K)} \sup\limits_{m \leqslant s} \log^{+} \left\vert\gf{T^m G_m}{m!} \right\vert_{\mathfrak{p}, \mathrm{Gauss}} \leqslant s \sum\limits_{\mathfrak{p} \in \Spec(\Oal_\K)} \sup\limits_{m \leqslant s} \log^{+} |T|_{\mathfrak{p}, \mathrm{Gauss}}+h(s,\mathfrak{p})$$

Symétriquement, en remarquant que $$\gf{G_m}{m!}=\left(\gf{1}{T}\right)^m \gf{T^m G_m}{m} \;\;\; \text{et} \;\;\; \log^{+}\left\vert\gf{1}{T}\right\vert_{\mathfrak{p},\mathrm{Gauss}}=-\log^{-}|T|_{\mathfrak{p},\mathrm{Gauss}}\;,$$ on a $$h(s,\mathfrak{p}) \leqslant -s\sum\limits_{\mathfrak{p} \in \Spec(\Oal_\K)} \sup\limits_{m \leqslant s} \log^{-} |T|_{\mathfrak{p}, \mathrm{Gauss}}+\sum\limits_{\mathfrak{p} \in \Spec(\Oal_\K)} \sup\limits_{m \leqslant s} \log^{+} \left\vert\gf{T^m G_m}{m!}\right\vert_{\mathfrak{p},\mathrm{Gauss}}.$$

En divisant par $s$ et en passant à la limite supérieure, on obtient le résultat voulu.\end{dem}

\begin{coro}
Le système $y'=Gy$ satisfait la condition de Galochkin de la définition \ref{galochkin} si et seulement si $\sigma(G) < +\infty$.
\end{coro}

\begin{dem}
Selon le point \textbf{a)} de la proposition \ref{proplienq'ssigmaG}, on a $$\limsup\limits_{s \rightarrow + \infty} \gf{\log q'_s}{s} < \infty \ssi \limsup\limits_{s \rightarrow + \infty} \gf{\log q_s}{s} < \infty.$$ De plus, selon le point \textbf{b)}, on a $\limsup\limits_{s \rightarrow + \infty} \gf{\log q'_s}{s} < +\infty \ssi \sigma(G) < +\infty$. Par suite $$\sigma(G) < \infty \Longleftrightarrow \exists C >0 : \forall s \in \N, q_s \leqslant C^{s+1},$$ ce qui n'est autre que la condition de Galochkin.
 \end{dem}
 
\subsubsection{Théorème d'André-Bombieri}
Le théorème suivant donne le lien entre le rayon de convergence global et la taille, et donc entre la condition de Bombieri et la condition de Galochkin. 

\bigskip
\begin{Th}[André-Bombieri] \label{bombieriandre}
Si $\K$ est un corps de nombres et $G \in \mathcal{M}_n(\K(z))$, alors $$\rho(G) \leqslant \sigma(G) \leqslant \rho(G)+n-1.$$ En particulier, si la condition de Galochkin est vérifiée, alors $y'=Gy$ est un système globalement nilpotent.
\end{Th}

\begin{dem}

La preuve consiste pour l'essentiel à montrer en deux étapes que si $\mathfrak{p}$ est un premier de $\K$, alors $\log^{+} \left(\gf{1}{R_{\mathfrak{p}}} \right) = \lim\limits_{s \rightarrow + \infty} \gf{1}{s} h(s, \mathfrak{p})$, chacune des deux étapes fournissant une inégalité du théorème.

\begin{itemize}
\item \textbf{Étape 1} : soit $\mathfrak{p} \in \Spec(\Oal_\K)$, montrons que $$\log^{+} \left(\gf{1}{R_{\mathfrak{p}}} \right) = \limsup\limits_{s \rightarrow + \infty} \gf{1}{s} h(s, \mathfrak{p}).$$

Selon le corollaire \ref{Hadamard}, on a $$R_{\mathfrak{p}}=\liminf_{s \rightarrow + \infty} \left\vert \gf{G_s}{s!} \right\vert_{\mathfrak{p}, \mathrm{Gauss}}^{-1/s}.$$ Donc

\begin{align}
\log^{+} \left(\gf{1}{R_{\mathfrak{p}}} \right) &= \log^{+} \left( \limsup_{s \rightarrow + \infty} \left\vert \gf{G_s}{s!} \right\vert_{\mathfrak{p}, \mathrm{Gauss}}^{1/s} \right)&\overset{\log^{+} \, \mathrm{croissante}}{\leqslant} \limsup_{s \rightarrow + \infty} \gf{1}{s} \log^{+} \left( \left\vert \gf{G_s}{s!} \right\vert_{\mathfrak{p}, \mathrm{Gauss}}\right) \notag \\
& \leqslant  \limsup_{s \rightarrow + \infty} \gf{1}{s} h(s, \mathfrak{p}). \label{inegalite1etape1ba}
\end{align}

Par ailleurs, selon le corollaire \ref{estimationrcvdwork}, en tenant compte de la normalisation de $| \cdot |_{\mathfrak{p}}$ choisie en début de sous-section,
$$\forall m \in \N, \quad \left\vert \gf{G_m}{m!} \right\vert_{\mathfrak{p}, \mathrm{Gauss}} \leqslant \gf{1}{R_{\mathfrak{p}}^m} \{ m, n-1 \}_{\mathfrak{p}} \sup_{j \leqslant n-1} |G_j|_{\mathfrak{p}, \mathrm{Gauss}},$$

où $$\{ m, n-1 \}_{\mathfrak{p}}=\sup\limits_{1 \leqslant \lambda_1 < \dots < \lambda_{n-1} \leqslant m} \gf{1}{|\lambda_1 \dots \lambda_{n-1}|_{\mathfrak{p}}}=\{ m, n-1 \}_p^{d_{\mathfrak{p}}/\delta}.$$ Donc 

$$ \log^{+} \left\vert \gf{G_m}{m!} \right\vert_{\mathfrak{p}, \mathrm{Gauss}} \leqslant m \log^{+} \left( \gf{1}{R_{\mathfrak{p}}} \right)+\gf{d_{\mathfrak{p}}}{\delta} \log^{+} \{m,n-1 \}_p+C_{\mathfrak{p}},$$
où $C_{\mathfrak{p}}=\log^{+} \sup\limits_{j \leqslant n-1} |G_j|_{\mathfrak{p}, \mathrm{Gauss}}$ est indépendant de $m$. Selon la définition-proposition \ref{majorationaccolade}, si $p \leqslant~m$, on a $\{ m, n-1 \}_p \leqslant m^{n-1}$ et si $p > m$, on a $\{ m, n-1 \}_p =1$. De plus, comme $\sum\limits_{\mathfrak{p} \mid p} d_{\mathfrak{p}}=\delta$, on a $\gf{d_{\mathfrak{p}}}{\delta} \leqslant 1$. Par suite, dans tous les cas, si $m \leqslant s$, 
$$\log^{+} \left\vert \gf{G_m}{m!} \right\vert_{\mathfrak{p}, \mathrm{Gauss}} \leqslant s \log^{+} \left( \gf{1}{R_{\mathfrak{p}}} \right)+ (n-1) \log s +C_{\mathfrak{p}},$$ d'où, si $p \leqslant s$,

\begin{equation} \label{conclusionetape1ba}
\gf{1}{s} h(s, \mathfrak{p}) \leqslant \log^{+} \left( \gf{1}{R_{\mathfrak{p}}} \right)+(n-1) \gf{\log s}{s} + \gf{C_{\mathfrak{p}}}{s}
\end{equation} 

et si $p >s$, on a $\gf{1}{s} h(s, \mathfrak{p}) \leqslant  \log^{+} \left( \gf{1}{R_{\mathfrak{p}}} \right) + \gf{C_{\mathfrak{p}}}{s}$. Ainsi, $$\limsup\limits_{s \rightarrow + \infty} \gf{1}{s} h(s, \mathfrak{p}) \leqslant \log^{+} \left( \gf{1}{R_{\mathfrak{p}}} \right).$$ Ce résultat conjugué à \eqref{inegalite1etape1ba} donne $$\limsup\limits_{s \rightarrow + \infty} \gf{1}{s} h(s, \mathfrak{p}) = \log^{+} \left( \gf{1}{R_{\mathfrak{p}}} \right).$$
 
\bigskip
On a $\forall j \in \{0, \dots, n-1\}, |G_j|_{\mathfrak{p}, \mathrm{Gauss}} \leqslant 1$ pour tous les premiers $\mathfrak{p}$ sauf un nombre fini, de sorte que la constante $C_{\mathfrak{p}}$ est nulle sauf pour $\mathfrak{p}$ dans cet ensemble fini de premiers. Ainsi, en utilisant \eqref{conclusionetape1ba}, on peut trouver une constante $C>0$ telle que 

\begin{align*}
\sum\limits_{\mathfrak{p} \in \Spec(\Oal_\K)} \gf{1}{s} h(s, \mathfrak{p}) &\leqslant  \sum\limits_{\mathfrak{p} \in \Spec(\Oal_\K)} \log^{+}\left(\gf{1}{R_{\mathfrak{p}}}\right) + \sum\limits_{p \leqslant s} \sum\limits_{\mathfrak{p} \mid p} \gf{d_{\mathfrak{p}}}{\delta} (n-1) \gf{\log s}{s}+\gf{C}{s} \\
&\leqslant  \sum\limits_{\mathfrak{p} \in \Spec(\Oal_\K)} \log^{+}\left(\gf{1}{R_{\mathfrak{p}}}\right) +(n-1) \pi(s) \gf{\log s}{s}+\gf{C}{s},
\end{align*}
puisque pour tout premier $p$, $\sum\limits_{\mathfrak{p} \mid p} d_{\mathfrak{p}}=\delta$. Ici, $\pi$ est la fonction de comptage des nombres premiers. L'équivalent, fourni par le théorème des nombres premiers, $\pi(s) \sim \gf{s}{\log s}$ donne par passage à la limite de part et d'autre, 

$$\sigma(G) \leqslant \rho(G) + n-1.$$

\item \textbf{Étape 2} : soit $\mathfrak{p} \in \Spec(\Oal_\K)$, montrons que $$\log^{+} \left(\gf{1}{R_{\mathfrak{p}}} \right) = \liminf\limits_{s \rightarrow + \infty} \gf{1}{s} h(s, \mathfrak{p}).$$
Soit $y$ vecteur colonne tel que $Dy=Gy$. Soient $s, m \in \N$, on a $D^s y=G_s y$ et 

\begin{align*}
\gf{G_{s+m}}{(s+m)!} y &= \gf{D^m}{(s+m)!} D^s y=\gf{D^m}{(s+m)!}(G_s y) \overset{\mathrm{Leibniz}}{=} \gf{1}{(s+m)!} \sum_{i=0}^m \binom{m}{i} (D^i G_s) D^{m-i} y \\
&= \gf{1}{(s+m)!} \sum_{i=0}^m \binom{m}{i} (D^i G_s)G_{m-i} y.
\end{align*}

Donc 

$$ \gf{G_{s+m}}{(s+m)!}=\sum_{i=0}^{m} \gf{s! m!}{(s+m)!} \left(\gf{D^i G_s}{i! s!} \right) \gf{G_{m-i}}{(m-i)!}. $$

Or, pour tout $i \in \{0, \dots, m \}$, selon la proposition \ref{effet||0deriv}, $\left\vert \gf{D^i}{i!} \left(\gf{G_s}{s!} \right) \right\vert_{\mathfrak{p}, \mathrm{Gauss}} \leqslant \left\vert \gf{G_s}{s!}  \right\vert_{\mathfrak{p}, \mathrm{Gauss}}.$

Donc$$ \left\vert \gf{G_{s+m}}{(s+m)!} \right\vert_{\mathfrak{p}, \mathrm{Gauss}} \leqslant \left\vert \gf{G_s}{s!} \right\vert_{\mathfrak{p}, \mathrm{Gauss}} \left\vert \gf{s! m!}{(s+m)!}   \right\vert_{\mathfrak{p}}  \sup_{j \leqslant m} \left\vert \gf{G_j}{j!} \right\vert_{\mathfrak{p}, \mathrm{Gauss}},$$
de sorte que 
$$ \log^{+} \left\vert \gf{G_{s+m}}{(s+m)!} \right\vert_{\mathfrak{p}, \mathrm{Gauss}} \leqslant \log^{+} \left\vert \gf{G_s}{s!} \right\vert_{\mathfrak{p}, \mathrm{Gauss}} + h(m, \mathfrak{p})-\gf{d_{\mathfrak{p}}}{\delta} \log \left\vert \binom{s+m}{s}  \right\vert_{p}. $$

En effet, comme $| \cdot |_{\mathfrak{p}}$ est non archimédienne, $\left\vert \binom{s+m}{s}^{-1}  \right\vert_{p} \geqslant 1$, donc $$\log^{+} \left\vert \binom{s+m}{s}^{-1}  \right\vert_{p} = \log \left\vert \binom{s+m}{s}^{-1}  \right\vert_{p} =-\log\left\vert \binom{s+m}{s} \right\vert_{p}.$$ On obtient donc par récurrence que pour tout $k \in \N$,
\begin{align*}
\log^{+} \left\vert \gf{G_{s+km}}{(s+km)!} \right\vert_{\mathfrak{p}, \mathrm{Gauss}} &\leqslant \log^{+} \left\vert \gf{G_s}{s!} \right\vert_{\mathfrak{p}, \mathrm{Gauss}} + k h(m, \mathfrak{p})-\gf{d_{\mathfrak{p}}}{\delta}\log \left\vert \binom{s+m}{s}  \binom{s+2m}{s+m} \dots \binom{s+km}{s+(k-1)m} \right\vert_{p} \\
& \leqslant  \log^{+} \left\vert \gf{G_s}{s!} \right\vert_{\mathfrak{p}, \mathrm{Gauss}} + k h(m, \mathfrak{p})-\gf{d_{\mathfrak{p}}}{\delta} \log \left\vert \gf{(s+km)!}{s! (m!)^k}  \right\vert_{p}.
\end{align*}

Or, selon le lemme \ref{ordrefactorielle}, on a pour tout $h \in \N$, $\gf{-\log |h!|_{p}}{\log p} = \gf{h-S_h}{p-1}$, avec $h=a_0+a_1 p+ \dots + a_{\ell-1} p^{\ell-1}$ écrit en base $p$, $S_h=a_0+\dots+a_{\ell}$.

On remarque que comme $\forall i \in \{ 0, \dots, \ell-1 \}, 0 \leqslant a_i \leqslant p-1$, on a $\gf{S_h}{p-1} \leqslant \ell$, et puisque $h \geqslant p^{\ell-1}$, $\ell \leqslant 1+\gf{\log h}{\log p}$. Ici,

\begin{align*}
-\log \left\vert \gf{(s+km)!}{s! (m!)^k}  \right\vert_{p} &=  \left( \gf{(s+km)-S_{s+km}}{p-1} - \gf{s-S_s}{p-1} - k \gf{m-S_m}{p-1} \right) \log p \\
&=  \log p (S_s + k S_m-S_{s+km}) \leqslant  \log p \left( 1+\gf{\log s}{\log p} + k + k \gf{\log m}{\log p} \right)  \\
& \leqslant  (k+1) \log p+ k \log m + \log s  \leqslant  (k+1)(\log p + \log m) \;\; \text{si} \;\; s \leqslant m.
\end{align*}

Donc comme pour $s \leqslant m$, $\log^{+} \left\vert \gf{G_s}{s!} \right\vert_{\mathfrak{p}, \mathrm{Gauss}} \leqslant h(m, \mathfrak{p})$, on a si $s \leqslant m$, 
$$ \log^{+} \left\vert \gf{G_{s+km}}{(s+km)!} \right\vert_{\mathfrak{p}, \mathrm{Gauss}} \leqslant (k+1) \left( h(m, \mathfrak{p})+\gf{d_{\mathfrak{p}}}{\delta} \left(\log p + \log m \right)\right).$$

\medskip

Soient $N \geqslant 1$ et $m \geqslant 1$, écrivons $N=qm+s$, $0 \leqslant s < m$ la division euclidienne de $N$ par $m$. On a $q \leqslant \gf{N}{m}$, donc l'inégalité précédente devient 
$$ \log^{+} \left\vert \gf{G_N}{N!} \right\vert_{\mathfrak{p}, \mathrm{Gauss}} \leqslant (q+1)\left( h(m, \mathfrak{p})+\gf{d_{\mathfrak{p}}}{\delta} \left(\log p + \log m \right)\right) \leqslant \gf{N+m}{m} \left( h(m, \mathfrak{p})+\gf{d_{\mathfrak{p}}}{\delta} \left(\log p + \log m \right)\right),$$
inégalité qui est vérifiée \textit{a fortiori} pour $N' < N$, d'où 
$$ \gf{1}{N} h(N, \mathfrak{p}) \leqslant \left(\gf{1}{m}+ \gf{1}{N} \right) \left( h(m, \mathfrak{p})+\gf{d_{\mathfrak{p}}}{\delta} \left(\log p + \log m \right)\right).$$
Donc selon l'étape 1, 
\begin{equation} \label{lienliminflimsupba}
\forall m \in \N, \quad \log^{+} \left(\gf{1}{R_{\mathfrak{p}}} \right) = \limsup_{N \rightarrow + \infty} h(N, \mathfrak{p}) \leqslant \gf{1}{m} \left( h(m, \mathfrak{p})+\gf{d_{\mathfrak{p}}}{\delta} \left(\log p + \log m \right)\right).
\end{equation}

Considérons une extraction $\phi : \N \rightarrow \N$ telle que $$\gf{1}{\phi(r)} h(\phi(r), \mathfrak{p}) \xrightarrow[r \rightarrow + \infty]{} \liminf\limits_{s \rightarrow + \infty} \gf{1}{s} h(s, \mathfrak{p}).$$ En prenant $m=\phi(r)$ dans \eqref{lienliminflimsupba} et en passant à la limite, on obtient
$$\log^{+} \left(\gf{1}{R_{\mathfrak{p}}} \right) \leqslant \liminf\limits_{s \rightarrow +\infty} \gf{1}{s} h(s, \mathfrak{p}),$$
d'où $\log^{+} \left(\gf{1}{R_{\mathfrak{p}}} \right)=\lim\limits_{s \rightarrow +\infty} \gf{1}{s} h(s, \mathfrak{p})$. L'inégalité $\rho(G) \leqslant \sigma(G)$ découle finalement du lemme de Fatou puisque 
\begin{align*}
 \sum\limits_{ \mathfrak{p} \in \Spec(\Oal_\K)} \log^{+} \left(\gf{1}{R_{\mathfrak{p}}} \right) &= \sum\limits_{ \mathfrak{p} \in \Spec(\Oal_\K)} \liminf_{s \rightarrow +\infty} \gf{1}{s} h(s, \mathfrak{p}) \\
 &\leqslant \liminf_{s \rightarrow + \infty} \sum\limits_{ \mathfrak{p} \in \Spec(\Oal_\K)} \gf{1}{s} h(s, \mathfrak{p}) \\ 
&\leqslant \limsup_{s \rightarrow + \infty} \sum\limits_{ \mathfrak{p} \in \Spec(\Oal_\K)} \gf{1}{s} h(s, \mathfrak{p})=\sigma(G).
\end{align*}
\end{itemize}
Ceci conclut la démonstration.\end{dem}

\subsection{Démonstration du théorème \ref{katzchudandre}}

Pour finir, synthétisons les différents résultats obtenus en prouvant le théorème d'André-Chudnovsky-Katz (théorème \ref{katzchudandre}) cité en introduction.

Soit $f \in \Qbar\llbracket z\rrbracket $ une $G$-fonction, $L \in \Qbar(z)\left[\gf{\mathrm{d}}{\mathrm{d}z}\right]$ un opérateur différentiel d'ordre minimal pour $f$ tel que $L(f(z))=0$ et $\K$ un corps de nombres tel que $L \in \K(z)\left[\gf{\mathrm{d}}{\mathrm{d}z}\right]$. Alors

\begin{itemize}
\item Selon la remarque suivant le théorème des Chudnovsky (théorème \ref{chudnovsky}), la matrice compagnon associée $A_L$ associée à $L$ vérifie la condition de Galochkin.

\item Donc, selon le théorème d'André-Bombieri (théorème \ref{bombieriandre}), $y'=A_L y$ est un système différentiel globalement nilpotent.

\item La dernière remarque de la partie \ref{subsubsec:nilpotenceglobale} nous assure alors que $L \in \K(z)\left[\gf{\mathrm{d}}{\mathrm{d}z}\right]$ est un opérateur différentiel globalement nilpotent.

\item  Selon le théorème de Katz (théorème \ref{katz}) et la remarque qui le suit, $L$ est donc un opérateur singulier régulier en tout point de $\mathbb{P}^{1}(\Qbar)$ et ses exposants en tout point sont rationnels.

\item Les points de $\C \setminus \Qbar$ sont des points ordinaires puisque, en écrivant $L$ sous la forme $$L=P_{\mu}(z) \left(\gf{\mathrm{d}}{\mathrm{d}z}\right)^{\mu}+P_{\mu-1}(z) \left(\gf{\mathrm{d}}{\mathrm{d}z}\right)^{\mu-1}+\dots+P_{0}(z),$$ avec $P_i(z) \in \Qbar[z]$, ce qui est toujours possible quitte à multiplier $L$ à gauche par un polynôme convenable, toute singularité de $L$ est une racine de $P_{\mu}$, donc un élément de $\Qbar$.
\end{itemize} 

\newpage

\section*{Annexe : le théorème de Chebotarev}

\addcontentsline{toc}{section}{Annexe : le théorème de Chebotarev}

Intéressons-nous de plus près au théorème de Chebotarev dont on une conséquence est utilisée de manière cruciale dans la démonstration du théorème de Katz (sous-section \ref{subsec:thkatz}).

On peut définir trois notions de densité sur l'ensemble $\Spec(\Oal_\K)$ des premiers d'un corps de nombres $\K$. 

\begin{defi}
Soit $T \subset \mathrm{Spec}(\Oal_\K)$. On définit $$\zeta_{\K,T}(s) := \prod\limits_{\mathfrak{p} \in T} \gf{1}{1-N(\mathfrak{p})^{-s}},$$ où $N(\mathfrak{p})$ désigne $\mathrm{Card}\left( \Oal_\K/\mathfrak{p} \right)$, la norme de $\mathfrak{p}$.

S'il existe $n$ tel que $\zeta_{\K,T}^n$ se prolonge en une fonction méromorphe au voisinage de $1$ avec un pôle d'ordre $m$ en $1$ (en adoptant la convention qu'un zéro d'ordre $m$ est un pôle d'ordre $-m$) alors $T$ a une \emph{densité polaire} $d(T)=\gf{m}{n}$.
\end{defi}

On a quelques propriétés correspondant à ce que l'on est en droit d'attendre d'une densité, qui sont valables également pour les autres notions de densité. 

\begin{itemize}
\item L'ensemble de tous les idéaux premiers de $\K$ a une densité polaire de $1$. 
\item Si $T$ admet une densité polaire, $d(T) \geqslant 0$.
\item Un ensemble fini a une densité nulle.
\item Si $T$ est l'union disjointe de $T_1$ et $T_2$, et si deux d'entre eux ont une densité polaire, le troisième en a une et $d(T)=d(T_1)+d(T_2)$.
\item Si $T_1$ et $T_2$ ont des densités polaires et $T_1 \subset T_2$, alors $d(T_1) \leqslant d(T_2)$.
\end{itemize}

La deuxième notion, qui est celle utilisée pour définir la notion d'opérateur globalement nilpotent dans la section \ref{sec:opdiffnilpotents}, est la densité \emph{de Dirichlet} ou \emph{analytique} (cf \cite[pp. 255-257]{Descombes})
\begin{defi} \label{densitedirichlet}
Soit $T \subset \mathrm{Spec}(\Oal_\K)$. On dit que $T$ admet $d \geqslant 0$ pour \emph{densité de Dirichlet} lorsque $$ \gf{-1}{\log(s-1)} \sum_{\mathfrak{p} \in T} \gf{1}{N(\mathfrak{p})^s} \longrightarrow d$$ quand $s$ tend vers $1$ pour $s$ réel, $s >1$. 
\end{defi}

\begin{defi}
Soit $T \subset \mathrm{Spec}(\Oal_\K)$. On dit que $T$ admet $d \geqslant 0$ pour \emph{densité naturelle} si $$ \lim\limits_{n \rightarrow +\infty} \gf{\mathrm{Card}\left\{ \mathfrak{p} \in T : N(\mathfrak{p}) \leqslant n \right\}}{\mathrm{Card}\left\{ \mathfrak{p} \in \Spec(\Oal_\K) : N(\mathfrak{p}) \leqslant n \right\}}=d.$$
\end{defi}

\begin{prop} \label{prop:liensdensites}
Si la densité polaire de $T$ existe, alors sa densité de Dirichlet aussi et les deux quantités sont égales.
Si la densité naturelle de $T$ existe, alors sa densité polaire aussi et les deux quantités sont égales.
\end{prop}

En revanche, la réciproque de cette proposition est fausse : Serre (\cite[p. 126]{SerreCours}) cite un exemple d'ensemble de premiers ayant une densité de Dirichlet mais pas de densité naturelle.

\medskip

Une fois choisie une notion de densité parmi les trois ci-dessus, on peut donner un sens à l'expression \og{} presque tout \fg{}. 
\begin{defi}
On dit qu'une propriété est valable pour \emph{presque tout} premier $\mathfrak{p}$ d'un corps de nombres (relativement à la densité choisie) si elle vaut pour un ensemble de premiers de densité~$1$. Ceci est vrai en particulier, si elle vaut pour tous les premiers sauf un nombre fini.
\end{defi}

Le théorème suivant, dont la preuve dépasse largement le cadre de ce mémoire, a été démontré par Chebotarev en 1922. Pour plus de détails, on peut se référer à \cite[p. 131]{Serre81}, et \cite{Lenstra}. Une preuve de ce théorème pour la densité de Dirichlet peut être trouvée dans \cite[pp. 179--197]{JMCFT}.

\begin{Th}[Chebotarev] \label{Chebotarev}
Soit $\mathbb{L}/\K$ une extension de corps de nombres galoisienne. Soit $C \subset G=\Gal(\mathbb{L}/\K)$ une classe de conjugaison. Alors l'ensemble des idéaux premiers $\mathfrak{p}$ de $\Oal_\K$ qui sont non ramifiés et tels que $(\mathfrak{p}, \mathbb{L}/\K)=C$ a une densité naturelle de $\gf{|C|}{|G|}$, où $|X|:=\mathrm{Card}(X)$.
\end{Th}

Selon la proposition \ref{prop:liensdensites}, le théorème \ref{Chebotarev} est vrai en remplaçant \og{} densité naturelle \fg{} par \og{} densité polaire \fg{} ou \og{} densité de Dirichlet \fg{}. Dans toute la suite, on fixe une notion de densité. Les résultats énoncés vaudront pour toutes les notions de densité.

\medskip

Définissons les objets utilisés dans le théorème. Soit $\mathfrak{p}$ premier de $\Oal_\K$. Soit $\mathfrak{q}$ premier de $\Oal_\mathbb{L}$ au-dessus de $\mathfrak{p}$ et $D(\mathfrak{q})=\left\{ \sigma \in \Gal(\mathbb{L}/\K) : \sigma(\mathfrak{q})=\mathfrak{q} \right\}$ le groupe de décomposition de $\mathfrak{q}$.

Notons $k_{\mathfrak{p}}=\Oal_\K/\mathfrak{p}$ et $l_{\mathfrak{q}}=\Oal_\mathbb{L}/\mathfrak{q}$ les corps résiduels respectifs de $\mathfrak{p}$ et $\mathfrak{q}$. Notons $f$ l'indice d'inertie de $\mathfrak{q}/\mathfrak{p}$, c'est-à-dire le degré de l'extension $l_{\mathfrak{q}}/k_{\mathfrak{p}}$. On a un morphisme surjectif
$$\fonction{\varphi_{\mathfrak{q}}}{D(\mathfrak{q})}{\Gal(l_{\mathfrak{q}}/k_{\mathfrak{p}})}{\sigma}{\overline{\sigma} : x \mod \mathfrak{q} \mapsto \sigma(x) \mod \mathfrak{q}}.$$
En remarquant que $D(\mathfrak{q})$ est un stabilisateur de l'action transitive de $\Gal(\mathbb{L}/\K)$ sur les premiers de $\Oal_\mathbb{L}$ au-dessus de $\mathfrak{p}$, et en utilisant la formule $e(\mathfrak{q}/\mathfrak{p})fr=[\mathbb{L}:\K]$, où $r$ est le nombre de premiers au-dessus de $\mathfrak{p}$ dans $\Oal_\mathbb{L}$, on voit que si $\mathfrak{p}$ est non ramifié,  $D(\mathfrak{q})$ a pour cardinal $f$. Donc $\varphi_{\mathfrak{q}}$ est en réalité un isomorphisme de groupes.

En tant que groupe de Galois d'une extension de corps finis, $\Gal(l_{\mathfrak{q}}/k_{\mathfrak{p}})$ est engendré par le morphisme de Frobenius $\mathrm{Frob}_{\mathfrak{p}} : x \mapsto x^{N(\mathfrak{p})}$. On note $(\mathfrak{p}, \mathbb{L}/\K)_{\mathfrak{q}}=\varphi_{\mathfrak{q}}^{-1}(\mathrm{Frob}_{\mathfrak{p}})$, générateur de $D(\mathfrak{q})$. 

On vérifie que si $\mathfrak{q}'$ est un autre premier au-dessus de $\mathfrak{p}$ et $\sigma \in \Gal(\mathbb{L}/\K)$ est tel que $\mathfrak{q}'=\sigma(\mathfrak{q})$, alors $(\mathfrak{p},\mathbb{L}/\K)_{q'}=\sigma (\mathfrak{p}, \mathbb{L}/\K)_{\mathfrak{q}} \sigma^{-1}$. 

\begin{defi} \label{defifrobenius}
Soit $\mathfrak{p}$ un idéal premier non ramifié de $\K$ dans l'extension galoisienne $\mathbb{L}/\K$. Le \emph{Frobenius} de $\mathfrak{p}$ est la classe de conjugaison dans $\Gal(\mathbb{L}/\K)$ 

$$(\mathfrak{p}, \mathbb{L}/\K)=\left\{ (\mathfrak{p},\mathbb{L}/\K)_{q}, \mathfrak{q} \cap \Oal_\K=\mathfrak{p} \right\}.$$ 
\end{defi}

En particulier, si $\mathbb{L}/\K$ est une extension abélienne, le Frobenius de $\mathfrak{p}$ est un élément de $\Gal(\mathbb{L}/\K)$.

\begin{rqu} \label{reducmodpracines}
Si $f \in \Z[z]$ est \textbf{unitaire} et $p$ ne divise pas le discriminant $\Delta(f)$ de $f$, alors $f \mod p$ a même degré $n$ que $f$ et a autant de racines, simples, dans $\Oal_{\K}/\mathfrak{p}$, que $f$ dans $\K$. 

En effet, supposons que $f$ a $n$ racines distinctes $\alpha_1, \dots, \alpha_n$ dans $\K$.  Soit $x \in \K$ est tel que $f(x)=0$. Alors, comme $f$ est unitaire, on a $x \in \Oal_{\K}$ et en réduisant modulo $\mathfrak{p}$ pour $\mathfrak{p}$ premier de $\K$ au-dessus de $p$, on a $$(f \mod p)(x \mod \mathfrak{p})=0=f(x) \mod \mathfrak{p}.$$ Donc $\alpha_1 \mod \mathfrak{p}, \dots, \alpha_n \mod \mathfrak{p}$ sont des racines de $f \mod p$, qui sont deux à deux distinctes. En effet, $\Delta(f) \in \Z$, donc $\Delta(f) \mod \mathfrak{p} =\Delta(f) \mod p \neq 0$ car $p$ ne divise pas $\Delta(f)$ et $$ \Delta(f) \mod \mathfrak{p}=\prod_{1 \leqslant r < s \leqslant n} (\alpha_r - \alpha_s)^2 \mod p=\prod_{1 \leqslant r < s \leqslant n} (\alpha_r \mod \mathfrak{p} - \alpha_s \mod \mathfrak{p})^2 \neq 0.$$ Ainsi, via l'isomorphisme $\fonction{\varphi_{\mathfrak{p}}}{D(\mathfrak{p})}{\Gal(k_{\mathfrak{p}}/\mathbb{F}_p)}{\sigma}{\overline{\sigma}}$, on voit que l'action du morphisme de Frobenius $\mathrm{Frob}_{p} \in \Gal(k_{\mathfrak{p}}/\mathbb{F}_p)$ en tant que permutation de $\left\{ \alpha_1 \mod \mathfrak{p}, \dots, \alpha_n \mod \mathfrak{p} \right\}$ est la même que celle de $(p, \K/\Q)_{\mathfrak{p}}$ sur $\left\{ \alpha_1 , \dots, \alpha_n\right\}$.
\end{rqu}

\begin{prop}
Soit $f \in \Z[z]$ un polynôme unitaire non nul tel que $f \mod p \in \mathbb{F}_p[z]$ est scindé dans $\mathbb{F}_p[z]$ pour presque tout premier $p$. Alors $f$ est scindé dans $\Q[z]$.
\end{prop}

\begin{dem}
On note $\K$ un corps de décomposition de $f$. On prend $p$ un premier ne divisant pas $\Delta(f)$. Si $\mathfrak{p}$ est un premier de $\K$ au-dessus de $p$, on sait donc que $\mathfrak{p}$ est non ramifié. De plus, pour toute racine $\alpha$ de $f \mod p$, on a $\mathrm{Frob}_{p}(\alpha)=\alpha \ssi \alpha^{p}=\alpha \ssi \alpha \in \mathbb{F}_p$. Donc $f \mod p$ est scindé sur $\mathbb{F}_p$ si et seulement si $\mathrm{Frob}_{p}=\mathrm{\mathrm{id}}$. Via l'isomorphisme $\varphi_{\mathfrak{p}}$, cela équivaut à dire que $(p, \K/\Q)=\left\{\mathrm{id} \right\}$. 

Selon le théorème \ref{Chebotarev}, l'ensemble des premiers $p$ de $\Z$ tels que $(p, \K/\Q)=\left\{\mathrm{id} \right\}$ a pour densité $1/|\Gal(\K/\Q)|$. Or, comme $f \mod p$ est scindé sur $\mathbb{F}_p$ pour presque tout $p$, cette densité vaut $1$. Donc $|\Gal(\K/\Q)|=[\K : \Q]=1$, si bien que $f$ est scindé sur $\Q$.\end{dem}

\begin{coro} \label{corochebo1}
Soit $f \in \Z[z]$ un polynôme irréductible unitaire tel que $f \mod p \in \mathbb{F}_p[z]$ a un zéro dans $\mathbb{F}_p$ pour tout premier $p$ dans un sous-ensemble de densité $d > \gf{1}{2}$ de $\Spec(\Z)$. Alors $f$ est de degré $1$.
\end{coro}

On utilise le lemme suivant :
\begin{lem} \label{lemmecorochebo1}
Soit $G$ un groupe agissant transitivement sur un ensemble fini $\Omega$ de cardinal au moins 2. Alors il existe $\sigma \in G$ tel que $\forall \omega \in \Omega, \sigma(\omega) \neq \omega$.
\end{lem}

\begin{dem}
On utilise la formule de Burnside : si $\mathrm{Orb}(\Omega)$ est l'ensemble des orbites de l'action de $X$ sur $\Omega$, et pour tout $g \in G$, $\mathrm{Fix}_g=\left\{ \omega \in \Omega / g(\omega)=\omega \right\}$, alors

$$ |\mathrm{Orb}(\Omega)| \times |G|=\sum\limits_{g \in G} |\mathrm{Fix}_g|=|G| $$ car l'action est transitive. Mais $|\mathrm{Fix}_{\mathrm{id}}|=|\Omega|>1$, donc il existe $\sigma \in G$ tel que $|\mathrm{Fix}_{\sigma}|=0$, \emph{i.e.} $\sigma$ est sans point fixe.\end{dem}

\begin{dem}[du corollaire \ref{corochebo1}]
Par l'absurde, supposons que $\deg(f) >1$. Soit $\K$ un corps de décomposition de $f$ et $p$ un premier ne divisant pas $\Delta(f)$. Une racine de $f \mod p$ dans $\mathbb{F}_p$ correspond à un point fixe du morphisme de Frobenius par action sur les racines de $f \mod p$.

Selon la remarque suivant la définition \ref{defifrobenius}, si $f \mod p$ a une racine dans $\mathbb{F}_p$, tout élément de $(p, \K/\Q)$ a donc un point fixe par action sur les racines de $f$. 

Soit $G_0$ le sous groupe de $G=\Gal(\K/\Q)$ composé des $\sigma \in G$ tels que $\sigma$ a un point fixe. C'est aussi une classe de conjugaison dans $G$. Comme $f$ est irréductible, l'action de $G$ sur l'ensemble à au moins deux éléments de ses racines dans $\K$ est transitive. Selon le lemme \ref{lemmecorochebo1}, $|G_0| < |G|$. Donc l'ensemble des premiers $p$ de $\Z$ tels que $p$ est non ramifié et $(p, \K/\Q)=G_0$ a pour densité $\gf{|G_0|}{|G|} \leqslant \gf{1}{2}$ (car $G_0$ est un sous-groupe strict de $G$). Or, cet ensemble est en fait l'ensemble des premiers $p$ non ramifiés tels que $f \mod p$ a une racine dans $\mathbb{F}_p$, qui est de densité strictement plus grande que $\gf{1}{2}$ par hypothèse. On a donc une contradiction, par conséquent $\deg(f)=1$.\end{dem}

Nous pouvons maintenant prouver la proposition \ref{conseqchebokatz} utilisée pour démontrer le théorème de Katz :

\begin{coro}\label{conseqchebokatzannexe}
Soit $\K$ un corps de nombres et $\alpha \in \K$. Si pour presque tout premier $\mathfrak{p}$ de $\K$, $|\alpha|_{\mathfrak{p}} \leqslant 1$ et $\alpha \mod \mathfrak{p} \in \mathbb{F}_p$, où $(p)=\mathfrak{p} \cap \Z$, alors $\alpha \in \Q$. 
\end{coro}

\begin{dem}
Soit $d \in \Z$ tel que $d \alpha \in \Oal_\K$. On a pour presque tout premier $\mathfrak{p}$ de $\K$, si $(p)=\mathfrak{p} \cap \Z$, $d \mod \mathfrak{p}=d \mod p  \in \mathbb{F}_p$ et $(d \alpha) \mod \mathfrak{p}=(d \mod \mathfrak{p}) (\alpha \mod \mathfrak{p}) \in \mathbb{F}_p$ car $\alpha \mod \mathfrak{p} \in \mathbb{F}_p$. Donc le polynôme minimal unitaire $f \in \Z[z]$ de $d \alpha$ a une racine dans $\mathbb{F}_p$ pour presque tout $p$.  Selon le corollaire \ref{corochebo1}, $f$ est de degré $1$, donc $d \alpha \in \Z$, de sorte que $\alpha \in \Q$.\end{dem}

\newpage
\addcontentsline{toc}{section}{Références}
\nocite{*}
\printbibliography

\newpage

\end{document}